\documentclass[leqno,12pt]{amsart}

\usepackage{amssymb}
\usepackage{amscd}              
\usepackage[all]{xy}            
\xyoption{curve}

\CompileMatrices

\newcommand{\Cc}{{\mathcal{C}}}

\newcommand{\coker}{\text{coker}}
\newcommand{\comp}{\circ}
                  
\newcommand{\dual}{\vee}

\newcommand{\extensor}{\boxtimes}
\newcommand{\E}{{\mathcal{E}}}
\newcommand{\Eb}{\text{$\overline{\E}$}}
\newcommand{\Ext}{\text{Ext}}
\newcommand{\F}{{\mathcal{F}}}
\newcommand{\G}{{\mathcal{G}}}
\newcommand{\Gieseker}{\text{Gieseker}}

\newcommand{\Gln}{\text{Gl}_n}
\newcommand{\Grass}{\text{Grass}}
\newcommand{\GVB}{\text{GVB}}
\newcommand{\GVBD}{\text{GVBD}}

\newcommand{\Hc}{{\mathcal{H}}}
\newcommand{\Hilb}{\text{Hilb}}
\newcommand{\Hom}{\text{Hom}}
\newcommand{\id}{\text{id}}
\newcommand{\im}{\text{im}}
\newcommand{\injto}{\hookrightarrow}
\newcommand{\isomorph}{\cong}           
\newcommand{\isomto}{\overset{\sim}{\rightarrow}}  
\newcommand{\I}{{\mathcal{I}}}

\newcommand{\Isom}{\text{Isom}}
\newcommand{\Isomto}{\overset{\sim}{\longrightarrow}}
\newcommand{\J}{{\mathcal{J}}}

\newcommand{\KGl}{\text{KGl}}
\newcommand{\KGln}{\text{KGl}_n}
\newcommand{\Ll}{{\mathcal{L}}}                  
\newcommand{\m}{{\mathfrak{m}}}                  
\newcommand{\M}{{\mathcal{M}}}
\newcommand{\Md}{\M^{\dual}}
\newcommand{\N}{{\mathbb{N}}}                    
\newcommand{\Nc}{{\mathcal{N}}}
\newcommand{\Oo}{{\mathcal{O}}}                  
\newcommand{\Oplus}{\bigoplus}
\newcommand{\ot}{\leftarrow}
\newcommand{\p}{{\mathfrak{p}}}                  
\newcommand{\Pc}{{\mathcal{P}}}

\newcommand{\plim}{\underset{\leftarrow}{\text{lim}}\, }
\newcommand{\Pp}{{\mathbb{P}}}                   
\newcommand{\pr}{\text{pr}}
\newcommand{\Proj}{\text{Proj}\,}
\newcommand{\Qc}{{\mathcal{Q}}}
\newcommand{\Quot}{\text{Quot}}
\newcommand{\red}{\text{red}}
\newcommand{\rk}{\text{rk}\, }   
\newcommand{\Rb}{\text{$\overline{R}$}}
\newcommand{\Rr}{{\mathcal{R}}}
\newcommand{\Spec}{\text{Spec}\, }               
\newcommand{\Sc}{{\mathcal{S}}}
\newcommand{\Sym}{\text{Sym}\, }
\newcommand{\tensor}{\otimes}
\newcommand{\tfs}{\text{tfs}}

\newcommand{\TFS}{\text{TFS}}
\newcommand{\To}{\longrightarrow}

\newcommand{\univ}{\text{univ}}

\newcommand{\Ub}{\text{$\overline{U}$}}
 
\newcommand{\Vb}{\overline{V}}
\newcommand{\VB}{\text{VB}}
\newcommand{\vers}{\text{vers}}

\newcommand{\X}{{\mathcal{X}}}

\newcommand{\Y}{{\mathcal{Y}}}

\newtheorem{theorem}{Theorem}[section]
\newtheorem{proposition}[theorem]{Proposition}
\newtheorem{lemma}[theorem]{Lemma}

\theoremstyle{definition}
\newtheorem{definition}[theorem]{Definition}
\newtheorem{remark}[theorem]{Remark}

\newtheorem{construction}[theorem]{Construction}

\begin{document}

\title[A Degeneration of Moduli Stacks]
{A Gieseker Type Degeneration of Moduli Stacks of Vector Bundles 
       on Curves}
\author[Ivan Kausz]{Ivan Kausz}
\date{October 22, 2001.}
\address{NWF I - Mathematik, Universit\"{a}t Regensburg, 93040 Regensburg, 
Germany}
\email{ivan.kausz@mathematik.uni-regensburg.de}

\maketitle

\tableofcontents


\section{Introduction}

In this paper we construct a degeneration with nice properties
of the moduli stack of vector bundles on a smooth curve 
when the curve degenerates to a singular curve which is irreducible
with one double point.

Degeneration is a well known technique in the study of moduli 
spaces of bundles on curves. The technique can be described
as follows:
Suppose one is interested in some invariant $\Nc$ of the moduli
space $U(X)$ of (semistable) vector bundles (say of fixed rank and
degree) on a smooth curve $X$. Suppose furthermore that the invariant
does not really depend on the curve, but only on its genus $g$.
Then it is often useful to consider a proper flat family of genus $g$
curves over a one-dimensional connected base, say $\X\to Y$,
which is smooth over the complement of a point $y_0\in Y$ and 
whose fibre $\X_{y_0}$ over $y_0$ is singular. In addition one
requires that the genus of the normalization $\widetilde{\X}_{y_0}$
of $\X_{y_0}$ is strictly smaller than $g$. If $g\geq 1$, one
can allway find such a family. 
More specifically, there exists a family such that $\X_{y_0}$ is irreducible
with one double point (and thus the genus of $\widetilde{\X}_{y_0}$ 
is $g-1$).
Now one tries to construct a proper flat family $U(\X/Y)\to Y$
of varieties such that 
\begin{enumerate}
\item
for every $y\in Y\setminus\{y_0\}$ the
fibre $U(\X/Y)_y$ is isomorphic to $U(\X_y)$. 
\item
the invariant $\Nc$ makes sense and takes the same value
for {\em all} fibres $U(\X/Y)_y$.
\item
The varieties $U(\X/Y)_{y_0}$ and $U(\widetilde{\X}_{y_0})$  are
related in a way which enables one to compute the invariant 
for $U(\X/Y)_{y_0}$ in terms of the invariant for $U(\widetilde{\X}_{y_0})$.
\end{enumerate}
By this strategy one is finally reduced to the genus zero case 
where often one can determine the invariant directly.

I know of two examples where this technique has been employed
successfully:

One example is Sun's proof of the so called ``factorization rule''.
In this example, $U(X)$ is the space of semistable (parabolic)
bundles (of given rank and degree) on a smooth (pointed) curve $X$
 and the invariant $\Nc$ is the dimension of the space of global
sections of (some power of) the generalized theta line bundle
$\Theta$ on $U(X)$. In his proof Sun extends ideas of 
Narasimhan and Ramadas (cf. \cite{Narasimhan-Ramadas}) 
who did the rank two case.
For the degeneration technique he uses the family 
$U^{\tfs}(\X/Y)\to Y$ of moduli spaces of torsion free sheaves
as studied e.g. in \cite{Newstead} and \cite{Seshadri1}.
A variant of Bhosle's concept of generalized parabolic sheaves
(cf. \cite{Bhosle}) allows him to relate $U^{\tfs}(\X/Y)_{y_0}$
with $U(\widetilde{\X}_{y_0})$. 
Unfortunately his proof for that $\Nc$ is independent of
$y\in Y$ only works for $g\geq 3$. Therefore he cannot 
quite reduce to the genus zero case and his result does not
lead to a formula for $\Nc$.

The other example for a successful employment of the degeneration
technique is Gieseker's proof of a conjecture of Newstead and Ramanan.
In this example $U(X)$ is the space of stable vector bundles of rank
2 and odd degree on a smooth curve $X$ of genus $g$ and $\Nc$ is
the $m$-th Chern class of the tangent bundle of $U(X)$ where 
$m>2g-2$. The conjecture of Newstead and Ramanan says that $\Nc=0$.
For the degeneration technique, Gieseker uses a different family
of moduli spaces than Sun: The fibre over $y_0$ of Gieseker's
family $U^{\Gieseker}(\X/Y)\to Y$ is a moduli space 
$U^{\Gieseker}(\X_{y_0})$ for certain rank 2 vector bundles.
These vector bundles live either on $\X_{y_0}$ itself or on modifications
$\X'_{y_0}$ or $\X''_{y_0}$ of $\X_{y_0}$ which may be depicted as
follows:
$$
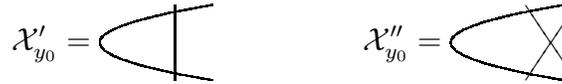
\X'_{y_0}=
\vcenter{
\xy
0;<1cm,0cm>:
(3,0.5); (3,-0.5) 
**\crv{(0,0)};
(2.5,0.5); (2.5,-0.5)
**@{-};
\endxy
}
\hspace{2cm}
\X''_{y_0}=
\vcenter{
\xy
0;<1cm,0cm>:
(3,0.5); (3,-0.5) 
**\crv{(0,0)};
(2.5,0.5); (3,-0.25)
**@{-};
(2.5,-0.5); (3,0.25)
**@{-};
\endxy
}
$$
In these figures the straight lines stand for copies of the
projective line $\Pp^1$ and the crooked ones stand for the
normalization $\widetilde{\X}_{y_0}$ of $\X_{y_0}$.
(As explained in the introduction of \cite{Kausz} the
space of vector bundles on $\X_{y_0}$ alone cannot be proper,
so a compactification of that space must parametrize additional objects).
The relationship between $U^{\Gieseker}(\X_{y_0})$ and 
$U(\widetilde{\X}_{y_0})$ is given by a diagram as follows:
$$
\vcenter{
\xymatrix{
& \Sc \ar[dl]_{\text{blowing up}} \ar[dr]^{\text{blowing up}} & \\
\Sc_{I} \ar[d]_{\text{locally trivial $\KGl_2$-fibration}}
& & \Sc_{II} \ar[d]^{\text{normalization}} \\
\text{$U(\widetilde{\X}_{y_0})$} & & U^{\Gieseker}(\X_{y_0})
}}
\eqno(*)
$$
Here, $\KGl_2$ is a certain compactification of the general
linear group $\text{Gl}_2$.
It is crucial for Gieseker's purpose that the variety 
$U^{\Gieseker}(\X/Y)$ is regular and that the fibre
 of $U^{\Gieseker}(\X/Y)$ over $y_0$ is a divisor with
normal crossings, properties which are not shared by the 
family $U^{tfs}(\X/Y)\to Y$.

Although expected by the experts (cf. introduction of \cite{Gieseker}
and \cite{Teixidor i Bigas}), a higher rank generalization of
Gieseker's family $U^{\Gieseker}(\X/Y)\to Y$ was missing for a long time. 
Only quite recently Nagaraj and Seshadri constructed such
a family in the arbitrary rank $n$ and coprime degree $d$ case
(cf. \cite{NS} and \cite{Seshadri2}). The fibre over $y_0$ of that
family is a moduli space $U^{\Gieseker}(\X_{y_0})$ for what we call
{\em stable Gieseker vector bundles} on $\X_{y_0}$. 
A Gieseker vector bundle on $\X_{y_0}$ is a vector bundle
$\E$ on a modification $\X^{(r)}_{y_0}$ of $\X_{y_0}$ of the
form 
$$
\X^{(r)}_{y_0}=
\vcenter{
\xy
0;<1cm,0cm>:
(3,0.5); (3,-0.5) **\crv{(0,0)};
(2.5,0.5); (3.6, 0.25) **@{-};
(3.4, 0.25); (4.6, 0.5) **@{-};
(5.5, 0.3) *{\dots};
(6.4, 0.5); (7.6, 0.25) **@{-};
(7.5, 0.4); (7.5, -0.4) **@{-};
(6.4, -0.5); (7.6, -0.25) **@{-};
(5.5, -0.3) *{\dots};
(3.4, -0.25); (4.6, -0.5) **@{-};
(2.5, -0.5); (3.6, -0.25)
**@{-};
\endxy
}
$$
(the number $r\geq 0$ is the length of the inserted chain
of projective lines, for $r=0$ we set $\X^{(0)}_{y_0}:=\X_{y_0}$)
with the following properties:
\begin{enumerate}
\item
The restriction of $\E$ to any of the inserted projective lines is
{\em strictly standard}, i.e. of the form 
$\Oo(1)^\delta\oplus\Oo^{n-\delta}$ for some $\delta>0$
which may depend on the projective line.
\item
The push forward $f_*\E$ of $\E$ by the canonical morphism 
$f:\X^{(r)}_{y_0}\to\X_{y_0}$ (which contracts the chain of projective
lines to the singular point of $\X_{y_0}$) is torsion-free.
\end{enumerate}
Such a Gieseker vector bundle $\E$ is called {\em stable}, if
in 2. the sheaf $f_*\E$ is a {\em stable} torsion-free sheaf.

Nagaraj and Seshadri show that the total space $U^{\Gieseker}(\X/Y)$
of their family is a regular projective variety and that the fibre
over $y_0$ is a divisor with normal crossings. Furthermore they
construct a canonical proper birational $Y$-morphism 
$U^{\Gieseker}(\X/Y)\to U^{\tfs}(\X/Y)$. On the other hand
they do not investigate the relationship between 
$U^{\Gieseker}(\X_{y_0})$ and $U(\widetilde{\X}_{y_0})$, i.e.
the higher rank analogue of diagram $(*)$. 
In particular they do not give a candidate for the compactification
of $\Gln$, which should be the higher rank generalization for $\KGl_2$.

In the quest for a geometric proof of the Verlinde formula
(cf. \cite{Sorger}) and about the same time as Nagaraj and
Seshadri, I arrived independently at a notion of (families of)
rank $n$ Gieseker vector bundles on an irreducible stable 
curve $C_0$ (over some field $k$) with one singularity 
(parmetrized by a $k$-scheme $S$). 
and announced
(in \cite{Kausz}) a result which assures that the associated
moduli functor $\GVB_n(C_0/k)$ fits into a diagram 
$$
\vcenter{
\xymatrix{
& \KGl(p_1^*\Pc,p_2^*\Pc) \ar[dl] \ar[dr] &\\
\text{$\VB_n(\tilde{C_0})$} & & \GVB_n(C_0/k)
}}
\eqno(**)
$$
where $\VB_n(\tilde{C_0})$ is the moduli functor of all 
- not just semistable - 
rank $n$ vector bundles on $\tilde{C_0}$ and
$\KGl(p_1^*\Pc,p_2^*\Pc)\to \VB_n(\tilde{C_0})$ is a locally trivial fibre
bundle with standard fibre the canonical compactification $\KGln$
of $\Gln$ as defined in \cite{Kausz}. 

In the present paper we prove more than that. Here is the precise statement.
Let $B$ be the spectrum of a discrete valuation ring,
$\Spec(k)\isomorph B_0\injto B$ its special and $B_1\injto B$
its generic point. We
assume that $k$ is algebraically closed and of characteristic
zero.
Let $C\to B$ be a stable curve over $B$, whose
generic fibre $C_1$  is smooth and whose special fibre
$C_0$ is irreducible with one singularity $p\in C_0$. 
We assume furthermore that $C$ is regular. 
Let $n\geq 1$ and let $\VB_n(C_1/B_1)$ be the algebraic 
$B_1$-stack of vector bundles
of rank $n$ over $C_1$.
We define the notion of (families of) Gieseker vector
bundles on $C$ and prove that the corresponding moduli
problem is representable by an algebraic $B$-stack
$\GVB_n(C/B)$ with the following properties:
\begin{enumerate}
\item
The stack $\GVB_n(C/B)$ is regular, locally of finite type and
flat over $B$. Its special fibre 
$\GVB_n(C_0/B_0)$ is a divisor with
 normal crossings (Theorem \ref{GVB is stack}).
\item
The generic fibre of $\GVB_n(C/B)$ is the moduli stack $\VB_n(C_1/B_1)$.
\item
In the diagram $(**)$ the arrows are morphisms of $B_0$-stacks.
$\KGl(p_1^*\Pc,p_2^*\Pc)$ is the fibre bundle over $\VB_n(\tilde{C_0})$
associated as in \cite{Kausz} Theorem 9.1 to the two vector bundles
$p_i^*\Pc$ ($i=1,2$) which are the pull-back of the Poincar{\'e} bundle on 
$\tilde{C_0}\times\VB_n(\tilde{C_0})$ by the sections given by the
points $p_1,p_2\in\tilde{C_0}$ lying above $p$. 
The morphism $\KGl(p_1^*\Pc,p_2^*\Pc)\to\GVB_n(C_0/k)$ identifies
$\KGl(p_1^*\Pc,p_2^*\Pc)$ with the normalization of 
$\GVB_n(C_0/k)$.
(Theorems \ref{normalization} and \ref{GVBD isomto KGl}).
\end{enumerate}

Returning to Nagaraj's and Seshadri's family $U^{\Gieseker}(\X/Y)\to Y$,
this result gives strong evidence for a relationship (by some blowing
up - blowing down procedure) between the normalization of
$U^{\Gieseker}(\X_{y_0})$ and a locally trivial  $\KGln$-fibration
over the variety $U(\widetilde{\X}_{y_0})$.
It also shows that the situation for stacks
is easier in this respect: No blowing up - blowing down is necessary;
the normalization of $\GVB_n(C_0/B_0)$ {\em is} already a locally
trivial $\KGln$-fibration over $\VB_n(\tilde{C_0})$.

There is also a stack version of the torsion-free sheaves approach:
For an arbitrary base scheme $B$ and any semistable curve
$C\to B$ there is an algebraic $B$-stack $\TFS_n(C/B)$ parametrizing
relatively torsion free sheaves of rank $n$ on $C/B$. However the stack
$\TFS_n(C/B)$ is not regular if $C\to B$ is not smooth. The
singularities of $\TFS_n(C/B)$ have been studied by Faltings in 
\cite{Faltings}. (That paper contains also an analogous treatment
of the stacks of principal $G$-bundles for $G$ the symplectic or
the orthogonal group). Returning to our case (i.e. where
$B$ is the spectrum of a discrete valuation ring, $C\to B$ generically
smooth with irreducible special fibre with only one singularity),
we have a morphism $\GVB_n(C/B)\to\TFS_n(C/B)$ which is defined in
the same way as the morphism 
$U^{\Gieseker}(\X/Y)\to U^{\tfs}(\X/Y)$ of Nagaraj and Seshadri
and which may be considered as a resolution of singularities of
$\TFS_n(C/B)$. The analogue in the torsion-free sheaves approach 
of the morphism $\KGl(p_1^*\Pc,p_2^*\Pc)\to\GVB_n(C_0/B_0)$ is a
morphism $\Grass_n(p_1^*\Pc\oplus p_2^*\Pc)\to\TFS_n(C_0/B_0)$.
We show (cf. Proposition \ref{GVB and TFS})     
that there exists a commutative diagram 
$$
\xymatrix{
\KGl(p_1^*\Pc,p_2^*\Pc) \ar[d] \ar[r] & \GVB_n(C_0/B_0) \ar[d] \\
\Grass_n(p_1^*\Pc\oplus p_2^*\Pc)\ar[r] & \TFS_n(C_0/B_0)
}
$$
The left vertical arrow has been defined in \cite{Kausz}
where we have also computed the fibres of that morphism.

Our paper is organized as follows:
In section \ref{section GVB} 
we define the notion of a Gieseker vector bundle on $C$
(resp. $C_0$)
over a $B$-scheme (resp. a $k$-scheme) $S$.
This is a pair $(\Cc,\E)$, where $\Cc$ is a nodal curve
over S which differs from the family $C\times_BS$ (resp. $C_0\times_kS$) 
only by the replacement of some of its singular fibres by modifications
$C_0^{(r)}$ of $C_0$ for some $r\leq n$, and where
the restriction of $\E$ to any fibre of the form $C_0^{(r)}$
is ``admissible'' in the following sense:
\begin{enumerate}
\item
It is strictly standard.
\item
Its push forward by $C_0^{(r)}\to C_0$ is torsion free.
\end{enumerate}
(Actually our definition (cf. \ref{admissible/X}) is different,
but equivalent to the above 
(cf. \ref{torsion-free})).
We then show that
the associated moduli problem is represented by a regular algebraic $B$-stack
$\GVB_n(C/B)$ and that its special fibre $\GVB_n(C_0/B_0)$
is a divisor with normal crossings.
This section relies heavily on results in \cite{NS}.

In section \ref{section GVBD} an 
auxiliary notion is introduced: that of a Gieseker
vector bundle data on $(\tilde{C_0},p_1,p_2)$ 
over a $k$-scheme $S$. For $S=\Spec(k)$ this
is just an admissible bundle on $C_0^{(r)}$ for some $r\geq 0$ plus the
specification of a singular point in $C_0^{(r)}$. 
The associated moduli problem
is  represented by an algebraic $B_0$-stack $\GVBD_n(\tilde{C_0},p_1,p_2)$
which we show to be canonically
isomorphic to the normalization of $\GVB_n(C_0/B_0)$.
Also in section \ref{section GVBD} we use results from \cite{NS}.

Sections \ref{section simple modifications} through 
\ref{section GVBD and KGl} are concerned with proving that
$\GVBD_n(\tilde{C_0},p_1,p_2)$ is a locally trivial $\KGln$-fibration over 
$\VB_n(\tilde{C_0}/B_0)$. For this, we  show that a 
Gieseker vector bundle data on $(\tilde{C_0},p_1,p_2)$ over 
a $k$-scheme $S$ uniquely corresponds to a vector bundle $\E$ over
$\tilde{C_0}\times_kS$ together with a generalized isomorphism 
(in the sense of \cite{Kausz} 5.2) from $p_1^*\E$ to $p_2^*\E$.
By definition, a generalized isomorphism is made up of what we
call bf-morphisms. A bf-morphism between vector bundles $E$ and $F$
over $S$ consists of
\begin{itemize}
\item
a pair $(L,\lambda)$ where $L$ is an invertible $\Oo_S$-module
and $\lambda$ is a global section of $L$.
\item
a diagram 
$
\xymatrix@C=2ex{
E
\ar@<-.5ex>[rr]
&&
F
\ar@<-.5ex>|{\tensor}[ll]
}
$
with certain properties,
where the arrow 
$
\xymatrix{
F \ar|{\tensor}[r] & E
}
$
means a morphism $F\to L\tensor E$ (cf. \cite{Kausz} 5.1).
\end{itemize}
Our strategy is to break up Gieseker vector bundle data into
simpler constituents and then to show a one to one correspondence
between these constituents and bf-morphisms.

In section \ref{section simple modifications} we show that
for a given curve $\pi:\Cc\to S$ together with a section $s$ of $\pi$
meeting $\Cc$ in the smooth locus,
giving data $(L,\lambda)$ is equivalent to giving a new
curve $\Cc'\to S$ (together with a section $s'$) which differs
from $\Cc$ only by that over points $z\in S$ where $\lambda$ vanishes,
a projective line is inserted at $s(z)$. We call $(\Cc',s')$ the
simple modification of $(\Cc,s)$ associated to $(L,\lambda)$.

In section \ref{pol} we prove a technical lemma needed in
section \ref{section simple modifications}.

In section \ref{section admissible bundles} we define the notion of
an admissible bundle over a simple modification. Then we show that
given a curve $\pi:\Cc\to S$, a section $s$ of $\pi$ and
a vector bundle $\E$ over $\Cc$ the data of a bf-morphism 
$
(L,\lambda,
\xymatrix{
E'
\ar@<-.5ex>[r]
&
E=s^*\E
\ar@<-.5ex>|(.6){\tensor}[l]
}
)
$
is equivalent to the data consisting of a simple modification
$(\Cc',s')$ of $(\Cc,s)$ together with an admissible bundle $\E'$
over $\Cc'$ which satisfies certain properties relative to $\E$ and $E'$.

In section \ref{section contractions} we prove that a Gieseker
vector bundle data gives rise to a sequence of simple modifications
together with admissible bundles on these.

Finally, in section \ref{section GVBD and KGl} we put together the 
results from section 
\ref{section simple modifications}-\ref{section contractions},
to prove that indeed Gieseker vector bundle data correspond to
generalized isomorphisms.

In the last section we establish the commutative diagram which 
clarifies the relationship of our construction with the
torsion-free sheaves approach.

A remark on the restrictions on $k$: Probably the results go through
without assuming algebraic closedness and characteristic zero. Indeed
the results in sections 
\ref{section simple modifications}-\ref{section tfs}
hold for general $k$. I imposed the restrictions on 
$k$ for the main results because 
I wanted to be able to cite
theorems from \cite{NS}
in sections \ref{section GVB} and \ref{section GVBD}.

During the preparation of this paper I had the opportunity to spend
six weeks at the Chennai Mathematical Institute. I would like to
thank V. Balaji, D.S. Nagaraj and C.S. Seshadri for their interest in my work
and all the other members of CMI for providing a friendly atmosphere.

\section{Notation}
\label{notation}
Here we collect some notations which will be used freely in this paper.
\begin{itemize}
\item
Throughout the article, $B$ denotes the spectrum of a discrete valuation ring,
$\Spec(k)\isomorph B_0\injto B$ its special and $B_1\injto B$
its generic point. 
Furthermore, we fix once and for all a stable curve $C/B$ 
with smooth generic fibre $C_1$ and irreducible special fibre
$C_0$ with only one singularity $p\in C_0$. We assume 
that $C$ is regular and that $\widehat{\Oo_{C_0,p}}\isomorph k[[u,v]]/(u.v)$.
The symbol $\tilde{C_0}$ denotes the normalization of $C_0$ and 
$p_1, p_2\in \tilde{C_0}$ are the points which map to the singular point 
$p\in C_0$.
\item
Let $X$ be a scheme, $x$ a point of $X$ and $\F$ a coherent $\Oo_X$-module.
Then we denote by $\F[x]$ the fibre $\F_x\tensor_{\Oo_{X,x}}\kappa(x)$
of $\F$ at $x$.
\item
For two integers $a,b$ the symbol $[a,b]$ denotes the set of integers
$c$ with $a\leq c\leq b$.
\item
Let $S$ be a scheme. Let $\pi:X\to S$ be a prestable curve over $S$ 
(i.e. $\pi$ is flat, proper, such that the geometric fibres are reduced
curves with at most ordinary double points).
Let $Y$ be the disjoint union of finitely many copies of $S$
and for $\nu=1,2$ let $i_\nu:Y\to X$ be closed immersions such that
\begin{enumerate}
\item
The morphisms $i_\nu$ are $S$-morphisms.
\item
The subschemes $i_1(Y)$ and $i_2(Y)$ of $X$ intersect trivially.
\item
For $\nu=1,2$ the subscheme $i_\nu(Y)$ is contained in the
smooth locus of $\pi$.
\end{enumerate}
By \cite{Knudsen}, Theorem 3.4 there exists a prestable curve
$X'/S$ and an $S$-morphism $g:X\to X'$, uniquely defined by the
following properties:
\begin{enumerate}
\item
$g\comp i_1=g\comp i_2$ and $g$ is universal with this property.
\item
as a topological space, $X'$ is the quotient of $X$ under the
equivalence relation $i_1(y)\sim i_2(y)$ for $y\in Y$.
\item
for $U'$ an open subset of $X'$ and $U:=g^{-1}(U')$ the equality
$\Gamma(U',\Oo_{X'})=\{f\in\Gamma(U,\Oo_X)\ |\ i_1^*f=i_2^*f\}$
holds.
\end{enumerate}
We refer to $X'$ as the {\em cokernel} of the double arrow 
$
\xymatrix{
Y \ar@<.6ex>[r]^{i_1} \ar@<-.6ex>[r]_{i_2} & X
}
$.
\item
An algebraic stack will always be understood in the sense of
\cite{LM}, (4.1).   
\end{itemize}

\section{Gieseker vector bundles}
\label{section GVB}

\begin{definition}
Let $K$ be a field. A {\em chain of rational curves over $K$} is
a connected projective prestable curve whose associated graph is 
linear and whose irreducible commponents intersect in $K$-rational
points and are each isomorphic to the rational line $\Pp_K^1$. 
The {\em length} of a chain of rational curves is the number of
its irreducible components. We define a chain
of length $0$ of rational curves over $K$ to be just $\Spec(K)$.
\end{definition}

Let $R$ be  a chain of rational curves over $K$ of length $r\geq 1$ 
and let 
$R_i\subseteq R$ $(1\leq i\leq r)$ be its irreducible components
such that $R_i\cap R_j\neq\emptyset$ if and only if $|i-j|\leq 1$.
For $1\leq i\leq r-1$ let $x_i\in R(K)$ be the ordinary double point,
in which the components $R_i$ and $R_{i+1}$ intersect. Furthermore,
let $x_0\in R_1(K)$ and $x_r\in R_r(K)$ be smooth points with $x_0\neq x_1$. 
Thus the situation may be
depicted es follows:

\vspace{3mm}

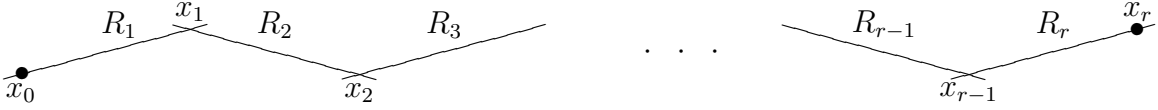
\begin{center}
\parbox{8cm}{
\xy <0mm,-4mm>; <.9mm,-4mm> :
(2.5,-6) *{\text{$x_0$}};
(0,-4);(30,4)**@{-} ?<<<*{\bullet}; (17,4) *{\text{$R_1$}};
(27.5,6) *{\text{$x_1$}};
(25,4);(55,-4)**@{-}; (40,4) *{\text{$R_2$}};
(52.5,-6) *{\text{$x_2$}};
(50,-4);(80,4)**@{-}; (65,4) *{\text{$R_3$}};
(95,0)*{\cdot};
(100,0)*{\cdot};
(105,0)*{\cdot};
(115,4);(145,-4) **@{-}; (130,4) *{\text{$R_{r-1}$}};
(142.5,-6) *{\text{$x_{r-1}$}};
(140,-4);(170,4) **@{-} ?>>>*{\bullet}; (155,4) *{\text{$R_r$}};
(167.5,6) *{\text{$x_r$}};
\endxy
}
\end{center}

\vspace{3mm}

Let $\E$ be a locally free $\Oo_R$-module of rank $n$ such that
for each $i$ one has
$$
\E|_{R_i}\isomorph\Oplus^{d_i}\Oo(1)\oplus\Oplus^{n-d_i}\Oo
$$ 
with certain $d_i\in[1,n]$. Note that although the direct sum 
decomposition of $\E|_{R_i}$ is not canonical, we have a canonical
exact sequence 
$$
0\to\F_i\to\E|_{R_i}\to\G_i\to 0
$$
such that $\F_i\isomorph\Oplus^{d_i}\Oo(1)$
and $\G_i\isomorph\Oplus^{n-d_i}\Oo$.
Indeed, $(\F_i\subseteq\E|_{R_i})$ is the Harder-Narasimhan filtration
of $\E|_{R_i}$. Therefore $\E$ induces a diagram of 
$K$-vector spaces

\vspace{3mm}
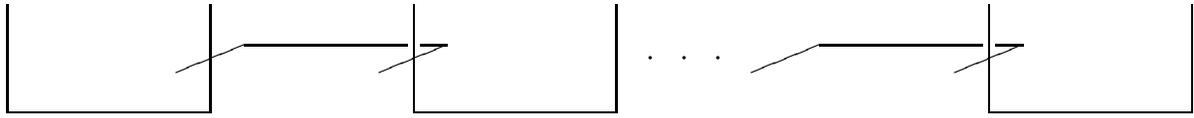
\begin{center}
\parbox{8cm}{
\xy <0mm,-4mm>; <.9mm,-4mm> :
(0,8);(0,-8)**@{-}; (0,-8);(30,-8)**@{-}; (30,-8);(30,8)**@{-};
(25,-2);(35,2)**@{-}; (35,2);(59,2)**@{-}, (61,2);(65,2)**@{-}; 
(65,2);(55,-2)**@{-};
(60,8);(60,-8)**@{-}; (60,-8);(90,-8)**@{-}; (90,-8);(90,8)**@{-};
(95,0)*{\cdot};
(100,0)*{\cdot};
(105,0)*{\cdot};
(110,-2);(120,2) **@{-}; (120,2);(144,2)**@{-}; (146,2);(150,2)**@{-}; 
(150,2);(140,-2)**@{-};
(145,8);(145,-8)**@{-}; (145,-8);(175,-8)**@{-}; (175,-8);(175,8)**@{-};
\endxy
}
\end{center}
\vspace{3mm}
made up of $r$ horseshoes,
where the $i$-th horseshoe stands for the diagram
$$
\xymatrix@R=2ex{
0\ar[d] &  & 0\ar[d] \\
\F_i[x_{i-1}]\ar[d] &  & \F_i[x_i]\ar[d] \\
\E[x_{i-1}]\ar[d] &  & \E[x_i]\ar[d] \\
\G_i[x_{i-1}]\ar[d] \ar@{=}[r] & 
H^0(R_i,\G_i) \ar@{=}[r] & 
\G_i[x_i]\ar[d] \\
0 &  & 0
}
$$
and the $i$-th and $(i+1)$-st horseshoe are linked together at $\E[x_i]$.
We define subspaces $V_i\subseteq \E[x_i]$ for $0\leq i\leq r$ 
inductively by setting
$V_0:=(0)$ and 
$V_i:=$
(preimage under $\E[x_i]\to\G_i[x_i]$ of the image of
 $V_{i-1}$ under the morphism $\E[x_{i-1}]\to\G_i[x_{i-1}]=\G_i[x_i]$).

\begin{definition}
\label{admissible}
Let $K$ be a field and let $R=\cup_{i=1}^r R_i$ be a chain of
rational curves over $K$. A locally free $\Oo_R$-module
of rank $n$ is called {\em admissible}, if either $R=\Spec(K)$ is
of length zero, or $R$ is of length $r\geq 1$ with successive 
irreducible components $R_1,\dots,R_r$ and the following 
two conditions hold:
\begin{enumerate}
\item
$\E$ is {\em strictly standard}, i.e.
for all $i$ 
the restriction of $\E$ to $R_i$ is isomorphic to 
$\Oplus^{d_i}\Oo(1)\oplus\Oplus^{n-d_i}\Oo$
for some $d_i\in[1,n]$
\item
With the above notation,
for all $i=1,\dots,r-1$  the two subspaces
$V_i$ and $\F_{i+1}[x_i]$ of $\E[x_i]$ intersect trivially.
\end{enumerate}
\end{definition}

\begin{lemma}
\label{equivalents for admissible}
Let $K$ be a field, let $R=\cup_{i=1}^r R_i$ be a chain of length $r\geq 1$ of
rational curves over $K$ and let $x_0\in R_1(K)$, $x_1\in R_r(K)$ be two 
different smooth points on the extremal components.
Let $\E$ be a strictly standard locally free $\Oo_R$-module
of rank $n$. Then the following conditions are equivalent:
\begin{enumerate}
\item
$\E$ is admissible
\item
$\dim H^0(R,\E(-x_0))= \sum_{i=1}^r \deg \E|_{R_i}$.
\item
$H^0(R,\E(-x_0-x_r))=(0)$.
\end{enumerate}
Furthermore, if one of these conditions holds, then 
for $i=1,\dots,r$ the subspace
$V_i\subseteq\E[x_i]$ which intervenes in \ref{admissible}.2.
coincides with the image of
$H^0(\E(-x_0)|_{\cup_{j=1}^iR_j})$ in $\E[x_i]$.
\end{lemma}

\begin{proof}
The equivalence of 2. and 3. has been shown in \cite{NS} lemma 2.
Consider the following commutative diagram of epimorphisms of 
$K$-vector spaces:
$$
\xymatrix@R=2.5ex{
& H^0(\E|_{R_{r}}) 
  \ar@{->>}[ld]_{\varphi_{r-1}} 
  \ar@{->>}[rd]^{\varphi_{r}} & 
\\
\E[x_{r-1}] 
  \ar@{->>}[d]^{\psi_{r-1}} & & 
\E[x_{r}] 
  \ar@{->>}[d]_{\psi_{r}} 
\\
\G_{r}[x_{r-1}] \ar@{=}[rr] & & \G_{r}[x_{r}]
}
$$
I claim that for any subspace $V\subseteq\E[x_{r-1}]$ we have
$\varphi_{r}(\varphi_{r-1}^{-1}(V))=\psi_{r}^{-1}(\psi_{r-1}(V))$.
Indeed, for this it is clearly sufficient to show
that the kernel of $\psi_{r}$ is contained in 
$\varphi_{r}(\varphi_{r-1}^{-1}(V))$.
But we have 
$
\ker(\varphi_{r-1})=H^0(\E|_{R_{r}}(-x_{r-1}))=
H^0(\F_{r}(-x_{r-1}))$
and therefore 
$
\ker(\psi_{r})=\F_{r}[x_{r}]=\varphi_{r}(H^0(\F_{r}(-x_{r-1})))
\subseteq \varphi_{r}(\varphi_{r-1}^{-1}(V))
$.
We now show by induction on the length $r$ of $R$ that conditions
1. and 2. are equivalent and that the equality
$
V_r=\im(H^0(\E(-x_0))\to\E[x_r])
$
holds.
For $r=1$ conditions 1. and 2. are trivially equivalent
and equality $V_1=\im(H^0(\E(-x_0))\to\E[x_1])$ follows
from the above claim applied to $V=(0)$.
If $r\geq 2$, we have
$$
H^0(\E(-x_0))=
H^0(\E(-x_0)|_{\cup_{i=1}^{r-1}R_i})
\times_{\E[x_{r-1}]}
H^0(\E|_{R_r})
\quad.
$$
Assuming condition 1. we clearly have 
$\dim(V_{r-1})=\sum_{i=1}^{r-1}\deg\E|_{R_i}$. 
By induction hypothesis we therefore have 
$\dim(H^0(\E(-x_0)|_{\cup_{i=1}^{r-1}R_i})=\dim(V_{r-1})$
and furthermore 
$V_{r-1}=\im(H^0(\E(-x_0)|_{\cup_{i=1}^{r-1}R_i})\to\E[x_{r-1}])$.
This implies in particular that the canonical map from
$
H^0(\E(-x_0)|_{\cup_{i=1}^{r-1}R_i})
$
to
$
\E[x_{r-1}]
$
is injective. 
Therefore $H^0(\E(-x_0))$ may be identified with the
subspace $\varphi_{r-1}^{-1}(V_{r-1})$ of $H^0(\E|_{R_r})$.
It follows that
$\dim(H^0(\E(-x_0)))=\dim(V_{r-1})+\dim(\ker\varphi_{r-1})
=\sum_{i=1}^{r}\deg\E|_{R_i}$, i.e. $\E$ satisfies condition 2.
and furthermore that we have
$$
\im(H^0(\E(-x_0))\to\E[x_r]) = 
\varphi_r(\varphi_{r-1}^{-1}(V_{r-1})) =
\psi_r^{-1}(\psi_{r-1}(V_{r-1})=V_r
\quad.
$$
Conversely, if $\E$ satisfies condition 2. and 3., then
the map 
$
H^0(\E(-x_0)|_{\cup_{i=1}^{r-1}R_i})
\to\E[x_{r-1}]
$
is injective, since its kernel is 
$H^0(\E(-x_0-x_{r-1})|_{\cup_{i=1}^{r-1}R_i})=(0)$.
As above, we conclude that
$\im(H^0(\E(-x_0))\to\E[x_r]) = V_r$.
Since $H^0(\E(-x_0))\to\E[x_r]$ is injective by condition 3.,
it follows that $\dim V_r = \sum_{i=1}^{r}\deg\E|_{R_i}$
which is clearly equivalent to $\E$ being admissible.
\end{proof}

\begin{remark}
\label{remark to admissibility}
Observe that admissibility of $\E$ implies that the
length of $R$ is at most the rank of $\E$.
From lemma \ref{equivalents for admissible}
it is easy to deduce the
following statement:
Let $R':=\cup_{i=1}^{r'}R_i$ and 
$R'':=\cup_{i=r'+1}^rR_i$ be subchains of $R$.
If $\E$ is admissible then $\E|_{R'}$ and $\E|_{R''}$ are
both admissible. Indeed, by \ref{equivalents for admissible}
it suffices to show that 
$H^0(R',\E|_{R'}(-x_0-x_{r'}))=(0)=H^0(R'',\E|_{R''}(-x_{r'}-x_r))$.
Suppose that there exists a nontrivial section
$\alpha'\in H^0(R',E|_{R'}(-x_0-x_{r'}))$, say. Then we can extend 
$\alpha'$ by zero to obtain a nontrivial section 
$\alpha\in H^0(R,\E(-x_0-x_r))$ which is impossible by assumption.
\end{remark}

\begin{definition}
\label{modification/K}
Let $K$ be a field and $X/K$ a curve. Let $x\in X(K)$ be a point with
$\widehat{\Oo_{X,x}}\isomorph K[[u,v]]/(uv)$ and let
$\tilde{X}\to X$ be the blowing-up of $X$ at $x$.
Let $x',x''\in\tilde{X}(K)$ be the two points above $x$.
A {\em modification of $X$ at $x$} is a morphism
$h:Y\to X$ such that 
\begin{enumerate}
\item
$Y$ is isomorphic to the cokernel of the double arrow 
(cf. section \ref{notation})
$$
\xymatrix{
\Spec(K)\amalg\Spec(K) 
\ar@<0.6ex>[r]^(.7){x'\amalg x''} \ar@<-0.6ex>[r]_(.7){x_0\amalg x_r} &
\text{$\tilde{X}$}\amalg R
}
$$
where $R$ is a chain of rational curves over $K$ of length $r$ and
$x_0, x_r$ are smooth $K$-rational points of $R$, which lie in different
extremal components if $r\geq 2$, are not equal if $r=1$, and are identical
if $r=0$.
\item
The composed morphism $\tilde{X}\to Y\to X$ (resp. $R\to Y\to X$)
is the blowing-up morphism 
(resp. the constant morphism which maps $R$ onto $x$).
\end{enumerate}  
\end{definition}
If $r=0$, then $Y=X$. For general $r$, a modification $Y$ looks like:

\vspace{3mm}
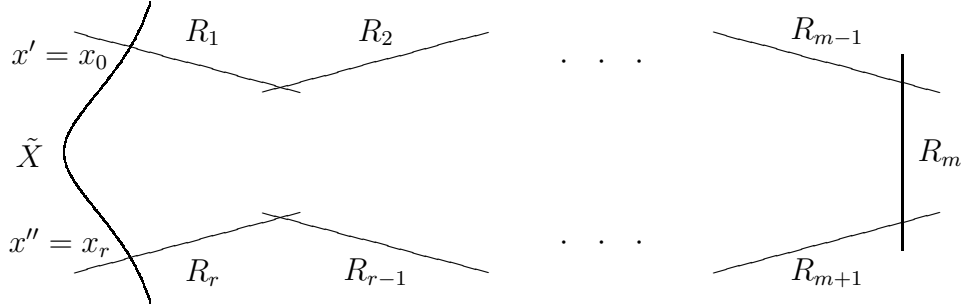
\begin{center}
\parbox{8cm}{
\xy <0mm,-20mm>; <1mm,-20mm>:
(10,20);(10,-20)**\crv{(7,10)&(-10,0)&(7,-10)};
(0,16);(30,8) **@{-}; (17,16) *{\text{$R_1$}};
(25,8);(55,16) **@{-}; (40,16) *{\text{$R_2$}};
(65,12)*{\cdot};
(70,12)*{\cdot};
(75,12)*{\cdot};
(85,16);(115,8) **@{-}; (100,16) *{\text{$R_{m-1}$}};
(110,13);(110,-13) **@{-}; (115,0) *{\text{$R_m$}};
(0,-16);(30,-8) **@{-}; (17,-16) *{\text{$R_r$}};
(25,-8);(55,-16) **@{-}; (40,-16) *{\text{$R_{r-1}$}};
(65,-12)*{\cdot};
(70,-12)*{\cdot};
(75,-12)*{\cdot};
(85,-16);(115,-8) **@{-}; (100,-16) *{\text{$R_{m+1}$}};
(-2,13) *{\text{$x'=x_0$}};
(-2,-12) *{\text{$x''=x_r$}};
(-6,0) *{\text{$\tilde{X}$}};
\endxy
}
\end{center}
\vspace{3mm}

\begin{definition}
\label{admissible/X}
Let $X/K$ and $x\in X(K)$ be as in definition \ref{modification/K}
and let $h:Y\to X$ be a modification of $X$ at $x$.
Let $R:=h^{-1}(x)$ be the chain of rational curves which is contracted
to $x$. A locally free $\Oo_Y$-module $\E$ is called 
{\em admissible for $Y\to X$},
if its restriction to $R$ is admissible.
\end{definition}

\begin{lemma}(\cite{NS}, Prop. 5)
\label{torsion-free}
With the notation of definition \ref{admissible/X} and assuming
$X$ to be irreducible with only one double point, we have:
A locally free $\Oo_Y$-module $\E$ is admissible for $Y\to X$ if and 
only if $\E$ is strictly positive and 
$h_*\E$ is torsion-free, where ``strictly positive'' means that
the restriction of $\E$ to any component of $R$ is of the
form $\oplus_i\Oo(a_i)$ for some $a_i\geq 0$ with
$\sum_ia_i\neq 0$.
\end{lemma}

\begin{definition}
\label{modification/S}
Let $C\to B$ be as in section \ref{notation}.
\begin{enumerate}
\item
Let $S$ be a $B$-scheme.
A {\em modification of $C$ over  $S$}
is a commutative diagram of schemes
$$
\xymatrix@R=2ex{
\Cc \ar[rr]^h \ar[dr]_{\pi} & & C\times_{B}S \ar[dl]^{\pr_2} \\
& S &
}
$$
where 
\begin{enumerate}
\item
$\pi$ is flat.
\item
$h$ is projective and finitely presented and induces an
isomorphism 
$\Cc\times_S S_1\isomto C_1\times_{B_1}S_1$,
where $S_1:=S\times_B B_1$.
\item
For all $z\in S\times_B B_0$ the induced morphism 
$\Cc_z\to C_0\times_{B_0}\Spec(\kappa(z))$ 
on the fibres over $z$ is a modification 
(in the sense of \ref{modification/K}) of the curve
$C_0\times_{B_0}\Spec(\kappa(z))$
in the point $p\times_{B_0}\Spec(\kappa(z))$.
\end{enumerate}
\item
Let $S$ be a $k$-scheme. A {\em modification of $C_0$ over $S$}
is a modification of $C$ over $S$, where $S$ is viewed as a
$B$-scheme via $S\to\Spec(k)=B_0\injto B$.
\end{enumerate}
\end{definition}

\begin{lemma}
\label{pi_*Oo}
Let $S$ be a $B$-scheme, let $h:\Cc\to C\times_BS$ be a modification of
$C$ over $S$ and let $\pi:\Cc\to S$ be the structure morphism.
Then the natural morphism $\Oo_S\to\pi_*\Oo_{\Cc}$ is an isomorphism.
\end{lemma}

\begin{proof}
First we show that $h_*\Oo_{\Cc}=\Oo_{X\times_BS}$.
By \cite{Knudsen}, Cor 1.5 the $\Oo_{X\times_BS}$-module
$h_*\Oo_{\Cc}$ is flat over $S$ and we have 
for all $z\in S$ the equality 
$h_*\Oo_{\Cc}\tensor_{\Oo_S}\kappa(z)=(h_z)_*\Oo_{\Cc_z}$,
where $h_z:\Cc_z\to X\times_B\Spec(\kappa(z))$ is the induced
morphism between the fibres over $z$.
Since obviously $(h_z)_*\Oo_{\Cc_z}=\Oo_{X\times_B\Spec(\kappa(z))}$,
it follows by Nakayama's lemma that we have indeed
$h_*\Oo_{\Cc}=\Oo_{X\times_BS}$.
Now consider the functor $T:M\mapsto H^0(C,\Oo_C\tensor_{\Oo_B}M)$
from the category of $\Oo_B$-modules into the category of $\Oo_B$-modules.
Since $C\to B$ is a stable curve, we have $H^0(C_b,\Oo_{C_b})=\kappa(b)$
for $b\in B$. Therefore $T$ is exact by \cite{EGA III}(7.8.4).
By \cite{EGA III}(7.7.5.3) it follows that 
$(\pr_2)_*\Oo_{C\times_BS}=H^0(C,\Oo_C)\tensor_{\Oo_B}\Oo_S$.
Thus it suffices to show that $H^0(C,\Oo_C)=\Oo_B$.
But this is clear: $H^0(C,\Oo_C)$ is a free $\Oo_B$-module
by \cite{EGA III}(7.8.4) and the morphism $\Oo_B\to H^0(C,\Oo_C)$
is an isomorphism by \cite{EGA III}(7.7.5)d).
\end{proof}

\begin{definition}
If  $\Cc\to C\times_B S$ is a modification of $C$ over
a $B$-scheme $S$ (cf. definition \ref{modification/S}), 
a locally free $\Oo_{\Cc}$-module $\E$ is called 
{\em admissible for $\Cc\to C\times_B S$},
if for all $z\in S\times_B B_0$ the restriction
of $\E$ to the fibre $\Cc_z$ over $z$ is admissible for 
$\Cc_z\to C_0\times_{B_0}\Spec(\kappa(z))$ in the sense
of definition \ref{admissible/X}.
\end{definition}

\begin{definition}
\label{definition of GVB}
Let $C\to B$ be as in section \ref{notation}.
Let $S$ be a $B$-scheme. 
A {\em Gieseker vector bundle (of rank $n$) on $C$ over $S$} is a modification 
$\Cc\to C\times_BS$ together with a locally free $\Oo_{\Cc}$-module 
$\E$ (of rank $n$), which is admissible for $\Cc\to C\times_BS$.
If $S$ is a $k$-scheme, a 
{\em Gieseker vector bundle on $C_0$ over $S$} is a
Gieseker vector bundle on $C$ over $S$, where $S$ is viewed as
a $B$-scheme via $S\to\Spec(k)=B_0\injto B$.
\end{definition}

There is an obvious notion of an isomorphism between two Gieseker
vector bundles on $C$ over a $B$-scheme $S$, and of pull-back
of a Gieseker vector bundle under a morphism $S'\to S$. Thus we have
a functor

\begin{eqnarray*}
\GVB_n(C/B) &:& \left\{
\begin{array}{ll}
\{\text{$B$-schemes}\}^o &\to\quad \{\text{groupoids}\} \\
\qquad S &
\mapsto\quad (\text{Gieseker vector bundles of rank $n$ on $C$ over $S$})
\end{array}
\right.
\end{eqnarray*}
Similarly, we have a functor

\begin{eqnarray*}
\GVB_n(C_0/k) &:& \left\{
\begin{array}{ll}
\{\text{$k$-schemes}\}^o &\to\quad \{\text{groupoids}\} \\
\qquad S &
\mapsto\quad (\text{Gieseker vector bundles of rank $n$ on $C_0$ over $S$})
\end{array}
\right.
\end{eqnarray*}
Since we will always work with bundles of fixed rank, 
we will often suppress the index $n$ in the notation
of these functors.
We may consider $\GVB(C_0/k)$ also as a functor on $\{\text{$B$-schemes}\}$
by setting $\GVB(C_0/k)(S):=\GVB(C_0/k)(S\to\Spec(k))$, if the structure
morphism $S\to B$ factorizes over $B_0=\Spec(k)$ and setting 
$\GVB(C_0/k)(S):=\emptyset$ else. 
In this sense we have the following cartesian diagram in the category
of $B$-groupoids:
$$
\xymatrix{
\GVB(C_0/k) \ar@{^{(}->}[r] \ar[d] & \GVB(C/B) \ar[d] \\
\Spec(k) \ar@{^{(}->}[r] & B
}
$$
i.e. $\GVB(C_0/k)$ is the special fibre of $\GVB(C/B)\to B$.

\begin{theorem}
\label{GVB is stack}
Assume that the residue field $k$ of $B$ is algebraically
closed and of characteristic zero. Then
the $B$-groupoid $\GVB(C/B)$ is an algebraic stack which is
locally of finite type and flat over $B$
and 
$\GVB(C_0/k)\injto\GVB(C/B)$ is a divisor with normal crossings. 
\end{theorem}

\begin{proof}
This follows from \ref{descent}, \ref{Isom},
\ref{thm H_N} and \ref{smooth domination} below.
\end{proof}

\begin{lemma}
\label{openness of being a closed immersion}
Let $S$ be a locally noetherian scheme and let $f:X\to Y$
an $S$-morphism  of proper $S$-schemes.
Let $s\in S$ be a point such that the morphism $f_s:X_s\to Y_s$
between the fibres over $s$ is a closed immersion.
Then there is an open neighbourhood  $S'\subseteq S$ of $s$
such that the induced morphism 
$X\times_SS'\to Y\times_SS'$ is a closed immersion. 
\end{lemma}

\begin{proof}
(I suppose this is well known but I couldn't find a suitable reference).
The assumptions imply that the morphism $f$ is proper.
By a theorem of Chevalley (cf. \cite{EGA IV}, (13.1.3)) there is an
open neighbourhood $Y'$ of the fibre $Y_s$ such that all fibres
of $X':=f^{-1}(Y')\to Y'$ are $0$-dimensional.
By Zariski's main theorem the morphism $X'\to Y'$ is finite
and by Nakayama's lemma the morphism $\Oo_{Y,y}\to(f_*\Oo_X)_y$
is surjective for all $y\in Y$. Therefore there is an open
neighbourhood $Y''\subseteq Y'$ of $Y_s$ such that 
$\Oo_{Y''}\to f''_*\Oo_{X''}$ is surjective, where
$X'':=f^{-1}(Y'')$ and $f'':=f|_{X''}$.
Let $q:Y\to S$ be the structure morphism and let 
$S':=S\setminus q(Y\setminus Y'')$.
Then $S'$ is an open neighbourhood of $s$ in $S$ with 
$Y\times_SS'\subseteq Y''$.
The morphism $f':X\times_SS'\to Y\times_SS'$ is finite
and $\Oo_{Y\times_SS'}\to f'_*\Oo_{X\times_SS'}$ is surjective.
Thus $f'$ is a closed immersion.
\end{proof}

\begin{lemma}
\label{open}
Let $C\to B$ be as in section \ref{notation}.
Let $S$ be a locally noetherian $B$-scheme, 
$h:\Cc\to C\times_BS$ a proper morphism
such that $\pi:=\text{pr}_2\comp h$ is flat,
and let $\E$ be a coherent $\Oo_{\Cc}$-module which is flat over $S$.
Let $z\in S$ be a point such that the pair
$(\Cc_z\to C\times_B\Spec(\kappa(z)); \E|_{\Cc_z})$
is a Gieseker vector bundle of rank $n$ on $C$ over 
$\Spec(\kappa(z))$. Then there is an open neighbourhood
$S^o$ of $z$ in $S$ such that
$(\Cc^o:=\Cc\times_SS^o\to C\times_BS^o; \E|_{\Cc^o})$
is a Gieseker vector bundle of rank $n$ on $C$ over $S^o$.
\end{lemma}

\begin{proof}
As has been shown in \cite{NS} (appendix), being a modification
is an open condition. Thus there exists an open neighbourhood
$U\subseteq S$ of $z$, such that  
$h_U:\Cc_U:=\Cc\times_SU
                  \to C\times_BU$
is a modification of $C$ over $U$
It is easy to see that after possibly shrinking $U$
we may also assume $\E_U:=\E|_U$
to be locally free of rank $n$.
We claim that if we choose $U$ small enough, we may even
assume  $\E_U$ to be admissible for $h_U$.
By defintion of a modification, $h_U$ is an isomorphism
over the generic fibre of $U\to B$, so admissibility of 
$\E_U$ over points in the generic fibre is trivial.
Thus we can restrict to the case
where $U\to B$ factorizes over the closed point $B_0$ of $B$.
Let $\F:=(h_U)_*\E_U$ and let $\E_z$ (resp. $\F_z$) be the
restriction of $\E$ (resp. of $\F$) to the fibre 
$\Cc_z:=\Cc\times_S\Spec(\kappa(z))$ 
(resp. to $C\times_B\Spec(\kappa(z))$).
By \cite{Knudsen} Cor. 1.5 we have
$\F_z=(h_z)_*\E_z$, where 
$h_z:\Cc_z:=\Cc\times_S\Spec(\kappa(z))\to C\times_B\Spec(\kappa(z))$
denotes the morphism between the fibres induced by $h$.
Therefore lemma \ref{torsion-free} implies that
$\E_z$ is strictly positive and that $\F_{z}$ is torsion free. 
We have to show that $\E_{z'}$ being strictly positive and
$\F_{z'}$ being torsion free are open conditions on $z'\in U$.
Let $\Ll$ be a very ample invertible sheaf on $C_0:=C\times_BB_0$. Then
the strict positivity of $\E_{z'}$ is equivalent to the following:
There exists an $N\in\N$ such that the sheaf 
$\E_{z'}(N):=\E_{z'}\tensor_{\Oo_{\Cc_{z'}}}
              h_{z'}^*(\Ll^{\tensor N}\tensor_{\Oo_B}\kappa({z'}))$
has the properties 
(i) $H^i(\Cc_{z'},\E_{z'}(N))=0$ for $i\geq 1$.
(ii) $\E_{z'}(N)$ is generated by its global sections.
(iii) the morpism from $\Cc_{z'}$ to the Grassmannian 
      $\Grass_n(H^0(\E_{z'}(N)))$
      induced by the surjection
      $H^0(\E_{z'}(N))\tensor_{\kappa({z'})}\Oo_{\Cc_{z'}}\to\E_{z'}$ 
      is a closed embedding.
With this characterization and lemma 
\ref{openness of being a closed immersion}    
it is easy to see that strict positivity
is indeed an open condition on $z'$. 
Since $C_0$ is irreducible, torsion freeness of $\F_{z'}$ is equivalent
to $\F_{z'}$ being of depth one (cf. \cite{Seshadri1}, p. 146, lemma 2).
By \cite{EGA IV} (12.2.2) this is an open condition on $z'$.
\end{proof}

The following lemma will allow us in some proofs to restrict to
the case of a noetherian base scheme.

\begin{lemma}
\label{noetherian}
Let $C\to B$ be as in section \ref{notation}.
Let $(S_\lambda)_\lambda$ 
be a projective system of affine $B$-schemes and let
$S:=\plim S_\lambda$.
Let
$(\Cc\to C\times_BS; \E)$ be a Gieseker vector bundle 
of rank $n$ on $C$ over $S$.
Then there exists an index $\alpha$, an open subscheme $S'\subset S_\alpha$
and a Gieseker vector bundle $(\Cc'\to C\times_BS'; \E')$ 
of rank $n$ on $C$ over $S'$ 
such that the morphism $S\to S_\alpha$ factorizes through $S'$ and
$(\Cc\to C\times_BS; \E)$
is induced by pull-back via $S\to S'$ 
from $(\Cc'\to C\times_BS'; \E')$.
\end{lemma}

\begin{proof}
By \cite{EGA IV} (8.8.2)(ii) there exists an index $\lambda$ 
and an $S_{\lambda}$-scheme of finite type $\Cc_{\lambda}$
with $\Cc\isomorph\Cc_{\lambda}\times_{S_{\lambda}}S$.
By \cite{EGA IV}, (8.10.5)(xii) and (11.2.6)(ii) we may assume
that $\Cc_{\lambda}\to S_{\lambda}$ is proper and flat.
By \cite{EGA IV} (8.8.2)(i), (8.5.2) and (11.2.6)(ii) 
we may assume furthermore that the morphism
$\Cc\to C\times_BS$ is induced by base-change
from a proper morphism $\Cc_{\lambda}\to C\times_BS_{\lambda}$ 
and that $\E$ is the pull-back of a coherent 
$\Oo_{\Cc_{\lambda}}$-module $\E_{\lambda}$, flat over $S_{\lambda}$.
It is easy to see that
if $K/\kappa$ is a field extension, $\Spec(\kappa)\to B$ a morphism,
$\Cc_{\kappa}\to C\times_B\Spec(\kappa)$ a proper morphism and $\E_{\kappa}$
a coherent $\Oo_{\Cc_{\kappa}}$-module such that the pull-back
$(\Cc_K\to C\times_B\Spec(K); \E_K)$ of the data
$(\Cc_\kappa\to C\times_B\Spec(\kappa); \E_\kappa)$ 
via $\Spec(K)\to\Spec(\kappa)$
is a Gieseker vector bundle
on $C$ over $\Spec(K)$, then there is a finite sub-extension 
$\kappa'/\kappa$ of $K/\kappa$,
such that already the pull-back of
$(\Cc_\kappa\to C\times_B\Spec(\kappa); \E_\kappa)$ 
via $\Spec(\kappa')\to\Spec(\kappa)$ is a Gieseker vector bundle
on $C$ over $\Spec(\kappa')$.
%
%
%
Therefore we can find for each $z\in S$ an index $\lambda_{z}\geq \lambda$,
such that the pull-back of $(\Cc_\lambda\to C\times_BS_\lambda;\E_\lambda)$
via $\Spec(\kappa(z))\to S_\lambda$ is a Gieseker vector bundle, where
$z'$ denotes the image of $z$ in $S_{\lambda_{z}}$. By lemma
\ref{open}, there is an open neighbourhood 
$U_{z}\subseteq S_{\lambda_{z}}$ of $z'$, such that the pull-back of 
$(\Cc_\lambda\to C\times_BS_\lambda;\E_\lambda)$
via $U_{z}\to S_\lambda$ is a Gieseker vector bundle 
on $C$ over $U_{z}$. Since $S$ is quasi-compact, there exist finitely
many points $z_1,\dots,z_m \in S$ such that the preimages of
the $U_{z_i}$ by $S\to S_{\lambda_{z_i}}$
cover $S$. We choose an index $\alpha$ which
dominates all the $\lambda_{z_i}$. Let 
$S'\subset S_{\alpha}$ be the union of all the preimages of 
the $U_{z_i}$ by $S_{\alpha}\to S_{\lambda_i}$. 
Then $S\to S_{\alpha}$ factorizes through 
$S$ and the pair 
$(\Cc\to C\times_BS'; \E):=
(\Cc_\lambda\times_{S_\lambda}S'\to C\times_B S'; 
 \E_\lambda\tensor_{\Oo_{S_\lambda}}\Oo_{S'})$
is a Gieseker vector bundle with the required property.
\end{proof}

\begin{proposition}
\label{descent}
Let $S$ be an affine $B$-scheme, let $g:S'\to S$ be an 
\'etale covering of $S$, let $S'':=S'\times_SS'$ and denote
by $p_i:S''\to S'$ the projection onto the $i$-th factor.
Let $(\Cc'\to C\times_BS'\ ,\ \E')$ be a Gieseker vector bundle
on $C$ over $S'$. Then every descent data
$$
f: p_1^*(\Cc'\to C\times_BS'\ ,\ \E')\isomto
   p_2^*(\Cc'\to C\times_BS'\ ,\ \E')
$$
for $(\Cc'\to C\times_BS'\ ,\ \E')$ 
with respect to $S'\to S$ is effective.
\end{proposition}

\begin{proof}
First of all we claim that it suffices to prove the proposition in the case
where $S$ is a locally noetherian scheme. Indeed, by 
\cite{EGA IV} (8.8.2)(ii), (8.10.5)(vi), (17.7.8)(ii)
there exists an \'etale covering $S'_0\to S_0$ 
of noetherian schemes which pulls back to the given covering
via a morphism $S\to S_0$. 
Lemma \ref{noetherian} implies that we may assume that there
exists a Gieseker vector bundle $(\Cc'_0\to C\times_BS_0', \E_0')$
which pulls back via $S'\to S_0'$ to the given one.
By \cite{EGA IV} (8.8.2)(i), (8.5.2)(i), (8.8.2.4), (8.5.2.4) 
we may assume that the
morphism $f$ comes from a morphism 
$f_0: p_1^*(\Cc'_0\to C\times_BS_0', \E_0')\isomto
      p_2^*(\Cc'_0\to C\times_BS_0', \E_0')$
and that $f_0$ is itself a descent data. This proves our claim.
So we can assume $S$ to be a locally noetherian scheme.
Let $T:=C\times_BS$, $T':C\times_BS'$, $T'':C\times_BS''$.
Then $T'\to T$ is an \'etale covering and $f$ induces
a descent data 
$p_1^*(\Cc',\det\E')\isomto p_2^*(\Cc',\det\E')$,
where $p_i:T''=T'\times_TT'\to T'$ is the $i$-th projection.
Since $\det\E'$ is relatively ample with respect to the morphism
$\Cc'\to T'$ (this follows e.g. by \cite{EGA III}, (4.7.1)),
this descent data is effective (\cite{SGA1}, exp. VII, Prop 7.8).
By \cite{SGA1}, exp. VII, Cor 1.3 and Prop. 1.10, also $\E'$
descends.
\end{proof}

\begin{proposition}
\label{Isom}
Let $S$ be an affine $B$-scheme and let 
$(\Cc_1\to C\times_BS\ ,\ \E_1)$,
$(\Cc_2\to C\times_BS\ ,\ \E_2)$
be two Gieseker vector bundles on $C$ over $S$.
Then the contravariant functor from the category of $S$-schemes
to the category of sets, which to every $S$-scheme $T$ associates
the set of isomorphisms from 
$(\Cc_1\to C\times_BS\ ,\ \E_1)_T$ to
$(\Cc_2\to C\times_BS\ ,\ \E_2)_T$,
is representable by a quasi-compact separable $S$-scheme.
\end{proposition}

The proof of proposition \ref{Isom} will be given after
lemma \ref{additivity of degree} below.

\begin{definition}
\label{definiton of degree}
Let $X/K$ be a projective curve over a field $K$ and let $\Ll$
be an invertible $\Oo_X$-module. We define 
$\deg_{X/K}(\Ll):=\chi(\Ll)-\chi(\Oo_X)$ where
for a coherent $\Oo_X$-module $\F$ we set 
$\chi(\F):=\dim_KH^0(X,\F)-\dim_KH^1(X,\F)$.
\end{definition}

\begin{lemma}
\label{additivity of degree}
Let $X$ be a prestable projective curve over an algebraically closed
field and let $\Ll$, $\M$ be two invertible $\Oo_X$-modules.
Then we have 
$\deg_{X/K}(\Ll\tensor\M)=\deg_{X/K}(\Ll)+\deg_{X/K}(\M)$.
\end{lemma}

\begin{proof}
Let $x\in X$ be an ordinary double point and let $f:\tilde{X}\to X$
be the blowing up of $X$ at $x$. Then there is an exact sequence 
of $\Oo_X$-modules as follows:
$$
0 \To \Ll \To f_*f^*\Ll \To \Ll[x] \To 0
\quad.
$$
From this and the analogous exact sequence for $\Oo_X$ 
one deduces easily the equality 
$deg_X(\Ll)=deg_{\tilde{X}}(f^*\Ll)$.
This shows that we may assume $X$ to be the disjoint union of smooth
projective curves or for that matter to be a single connected smooth
projective curve. The Riemann-Roch theorem shows that in this case
our definition of the degree of a line bundle coincides with 
the usual one, so the required formula holds.  
\end{proof}

\begin{proof}
(of proposition \ref{Isom})
Employing techniques used in the proof of lemma \ref{noetherian}
one sees 
that it suffices to show that in the case of 
noetherian $S$ the functor 
$T\mapsto\Isom((\Cc_1\to C\times_BS,\E_1)_T,(\Cc_2\to C\times_BS,\E_2)_T)$
from the category of noetherian $S$-schemes to the category of sets
is representable by a quasi-compact separable $S$-scheme.
Let $\Hilb_{\Cc_1\times_S\Cc_2/S}$ be the Hilbert scheme which 
parametrizes closed subschemes of $\Cc_1\times_S\Cc_2/S$ which are flat
over $S$. Let 
$Z\injto(\Cc_1\times_S\Cc_2)\times_S\Hilb_{\Cc_1\times_S\Cc_2/S}$
be the universal object.
It is easy to see that there is a locally closed 
subscheme $Y\injto \Hilb_{\Cc_1\times_S\Cc_2/S}$ with the property
that a morphism $T\to\Hilb_{\Cc_1\times_S\Cc_2/S}$ 
factorizes over $Y$ if and only if the following holds:
\begin{enumerate}
\item
The $i$-th projection $\pr_i:Z_T\to(\Cc_i)_T$ is an isomorphism for $i=1,2$.
\item
The diagram
$$
\xymatrix@R=1.5ex{
(\Cc_1)_T \ar[dr] & Z_T \ar[l]_{\pr_1} \ar[r]^{\pr_2} & (\Cc_2)_T \ar[dl] \\
& C\times_B T &
}
$$
is commutative.
\end{enumerate}
In particular, the pair
$(Y,\pr_2\comp\pr_1^{-1}:(\Cc_1)_Y\isomto Z_Y\isomto(\Cc_2)_Y)$
represents the (contravariant) functor
$$
T\mapsto 
\{\text{all $(C\times_B T)$-isomorphisms $(\Cc_1)_T\isomto(\Cc_2)_T$}\}
$$
from $\{\text{$S$-schemes}\}$ to $\{\text{sets}\}$.
I claim that $Y/S$ is of quasi-projective.
We may assume without loss of generality that $S$ is connected.
Let $\Ll_i$ be a very ample invertible sheaf on $\Cc_i$ 
relatively to $S$ $(i=1,2)$. Then the sheaf 
$\Ll:=\Ll_1\extensor_{\Oo_S}\Ll_2$ is very ample on 
$\Cc_1\times_X\Cc_2$ over $S$.
For a polynomial $P$ let $\Hilb_{\Cc_1\times_S\Cc_2/S}^P$
be the component of $\Hilb_{\Cc_1\times_S\Cc_2/S}$ parametrizing
objects whith Hilbert polynomial $P$ with respect to the polarization
$\Ll$. Since $\Hilb_{\Cc_1\times_S\Cc_2/S}^P$ is projective over $S$,
it is clearly enough for the proof of the claim to show that 
$Y$ is contained in $\Hilb_{\Cc_1\times_S\Cc_2/S}^P$ for some $P$.
Let $T$ be the spectrum of an algebraically closed field and
let $T\to S$ be a morphism. Let $Z_T\injto (\Cc_1\times_S\Cc_2)_T$
be a closed immersion with the above properties 1. and 2..
The Hilbert polynomial $P$ of the coherent 
$\Oo_{(\Cc_1\times_S\Cc_2)_T}$-module $\Oo_{Z_T}$  with respect to $\Ll_T$
is given by
$$
P(m)=\chi_{(\Cc_1)_T}((\Ll_1)_T^{\tensor m}\tensor f_T^*(\Ll_2)_T^{\tensor m})
\quad,
$$
where $f_T:=\pr_2\comp\pr_1^{-1}:(\Cc_1)_T\isomto(\Cc_2)_T$.
We have to show that $P$ is independent of $T\to S$ and of
$Z_T\injto(\Cc_1\times_S\Cc_2)_T$.
By lemma \ref{additivity of degree} we have
$$
P(m)=\deg_{(\Cc_1)_T}((\Ll_1)_T^{\tensor m})
    +\deg_{(\Cc_2)_T}((\Ll_2)_T^{\tensor m})
    +\chi_{(\Cc_1)_T}(\Oo_{(\Cc_1)_T})
$$
which shows the required independence.
This proves our claim and in particular that 
$Y$ is a quasi-compact separable $S$-scheme.

Let $f:(\Cc_1)_Y\isomto(\Cc_2)_Y$ be the universal isomorphism.
To prove the proposition it suffices now to show that the functor
$T\mapsto \Isom_{\Oo_{(\Cc_1)_T}}((\E_1)_T,f^*(\E_2)_T)$ 
from the category of $Y$-schemes to the category of sets
is representable by a quasi-compact separable $Y$-scheme.
But this is proven in \cite{LM} (proof of theorem (4.6.2.1)).
\end{proof}

\begin{definition}
\label{def H_N}
Let $C\to B$ be as in section \ref{notation}.
For every $B$-scheme $S$ let $H_N(S)$ be the set of all closed
subschemes $\Cc\injto C\times_B\Grass_n(\Oo_B^N)\times_BS$
with the following properties:
\begin{enumerate}
\item
The induced morphism $\Cc\to C\times_BS$ is a modification of
$C$ over $S$.
\item
The induced morphism $\Cc\to\Grass_n(\Oo_B^N)\times_BS$
is a closed immersion.
\item
Let 
$
\xymatrix{
\Oo_{\Cc}^N\ar@{->>}[r] & \E
}
$ 
be the pull back by 
$\Cc\injto\Grass_n(\Oo_B^N)\times_BS$
of the universal quotient. Then for every point $z\in S$
the induced morphism
$H^0(\Cc_z,\Oo_{\Cc_z}^N)\to H^0(\Cc_z,\E_z)$ 
is an isomorphism and the cohomology groups $H^i(\Cc_z,\E_z)$ vanish
for $i\geq 1$.
\item
The $\Oo_{\Cc}$-module
$\E$ is admissible for $\Cc \to C\times_BS$.
\end{enumerate}
The association $S\mapsto H_N(S)$ is in an obvious way a contravariant
functor $H_N$ from the category of $B$-schemes to the category of sets.
\end{definition}

\begin{theorem}
\label{thm H_N}
(\cite{NS})
Assume that the residue field $k$ of $B$ is algebraically closed
and of characteristic zero. Then
the functor $H_N$ defined in \ref{def H_N} 
is represented by a subscheme (also denoted $H_N$)
of the Hilbert scheme 
$\Hilb_{C\times_B\Grass_n(\Oo_B^N)/B}$.
Furthermore, the scheme $H_N$ is regular and flat over $B$, 
its generic fibre (over
$B$) is smooth and its special fibre is a divisor 
with normal crossings in $H_N$.
\end{theorem}

\begin{proof}
By \cite{NS} Prop. 8 there exists an open subscheme $\Y$ of
the Hilbert scheme
\linebreak[4]
$\Hilb_{C\times_B\Grass_n(\Oo_B^N)/B}$
which parametrizes closed subschemes satisfying properties 1-3
in definition \ref{def H_N}.
Furthermore, $\Y$ is regular, its generic fibre over $B$ is smooth
and its special fibre is a divisor with normal crossings.
Let $S$ be a $B$-scheme and let 
$\Cc\injto C\times_B\Grass_n(\Oo_B^N)\times_BS$
be an $S$-valued point of $\Y$.
Let $h:\Cc\to C\times_BS$ be the projection and let
$
\xymatrix{
\Oo_{\Cc}^N\ar@{->>}[r] & \E
}
$ 
be the pull-back by $\Cc\injto\Grass_n(\Oo_B^N)\times_BS$
of the universal quotient.
In \cite{NS} Prop. 9 it is shown that 
application of $h_*$ gives a surjection 
$
\xymatrix{
\Oo_{C\times_BS}^N\ar@{->>}[r] & h_*\E
}
$ 
such that $h_*\E$ is flat over $S$.
In particular, push-forward induces a morphism 
from $\Y$ to the $\Quot$-scheme $\Quot_{\Oo_{\Cc}^N/\Cc/B}$ 
which parametrizes  coherent quotients of $\Oo_{C}^N$ which
are flat over the base. Let $\Qc\subset\Quot_{\Oo_{\Cc}^N/\Cc/B}$
be the open subscheme parametrizing torsion-free quotients.
I claim that the preimage of $\Qc$ by the morphism 
$\Y\to\Quot_{\Oo_{\Cc}^N/\Cc/B}$ represents the functor $H_N$.
Indeed, let $S$ be a $B$-scheme which is the spectrum of a field,
let $\Cc\injto C\times_B\Grass_n(\Oo_B^N)\times_BS$
be an $S$-valued point of $\Y$, let $h:\Cc\to C\times_BS$ be the projection
and let $\E$ be the pull-back of the universal quotient.
Property 2. in definition \ref{def H_N} insures that $\E$ is strictly positive.
Therefore lemma \ref{torsion-free} tells us that $\E$ is admissible for
$h$ if and only if $h_*\E$ is torsion-free.
\end{proof}

\begin{definition}
We fix a very ample invertible sheaf $\Oo_C(1)$ on the curve $C$.
For a Gieseker vector bundle $(f:\Cc\to C\times_BS\ ,\ \E)$  on $C$ 
over a $B$-scheme $S$  and for an integer $N'$ 
we denote by $\E(-N')$ the (admissible) $\Oo_{\Cc}$-module
$\E\tensor_{\Oo_{\Cc}}f^*(\Oo_C(1)^{\tensor (-N')}\extensor_{\Oo_B}\Oo_S)$.
For every pair of integers $N\geq n$, $N'\geq 0$ we have a
morphism of $B$-groupoids
$$
P_{N,N'}: H_N\to\GVB_n(C/B)
$$
which to an $S$-valued point 
$(\Cc\injto C\times_B\times\Grass_n(\Oo_B^N)\times_BS)$ of $H_N$ associates
the $S$-valued point $(\Cc\to C\times_BS\ ,\ \E(-N'))$ of $\GVB_n(C/B)$,
where $\E$ is the pull back by $\Cc\injto\Grass_n(\Oo_B^N)\times_BS$ of
the universal quotient sheaf on the Grassmannian.
\end{definition}

\begin{lemma}
\label{S_NN'}
Let $S$ be an affine $B$-scheme and let 
$(\Cc\to C\times_BS\ ,\ \E)$ be Gieseker vector bundle of rank $n$
on $C$ over $S$. Let $N\geq n$ and $N'\geq 0$.
Then there is an open subscheme $S_{N,N'}\subseteq S$ such that
a morphism $T\to S$ factorizes over $S_{N,N'}$ if and only if
for every point $t\in T$ the following holds:
\begin{enumerate}
\item
$\dim H^i(\Cc_t,\E(N')_t) =
\left\{
\begin{array}{ll}
N &\text{for $i=0$} \\
0 &\text{else}
\end{array}
\right.
$
\item
The canonical morphism
$H^0(\Cc_t,\E(N')_t)\tensor_{\kappa(t)}\Oo_{\Cc_t}
\longrightarrow \E(N')_t$
is surjective.
\item
The morphism $\Cc_t\to\Grass_n(H^0(\Cc_t,\E(N')_t))$
induced by the surjection in 2. is a closed immersion.
\end{enumerate}
\end{lemma}

\begin{proof}
By \ref{noetherian} one is easily reduced to the case where
$S$ is noetherian. We define subsets 
$S'''\subseteq S''\subseteq S'\subseteq S$
as follows:
$S'$ is the set of points $s\in S$ with the property that
$H^i(\Cc_s,\E(N')_s)$ vanishes if $i>0$ and is of dimension
$N$ if $i=0$.
$S''$ is the set of points in $s\in S'$ such that
the canonical morphism 
$H^0(\Cc_s,\E(N')_s)\tensor_{\kappa(s)}\Oo_{\Cc_s}\to \E(N')_s$
is surjective.
Finally, $S'''$ is the set of points $s\in S''$ with the property
that the morphism $\Cc_s\to\Grass_n(H^0(\Cc_s,\E(N')))$ 
induced by the surjection 
$H^0(\Cc_s,\E(N')_s)\tensor_{\kappa(s)}\Oo_{\Cc_s}\to \E(N')_s$
is a closed
immersion. I claim that $S'''$ is open in $S$.
Indeed $S'$ is open in $S$ by the semicontinuity theorem 
(cf. \cite{EGA III} (7.7.5.1)) and by \cite{EGA III} (7.9.9).
Since $\Cc\to S$ is proper, an application of Nakayama's
lemma shows that $S''$ is open in $S'$. 
The openness of $S'''$ follows from lemma 
\ref{openness of being a closed immersion}.
It is clear that $S_{N,N'}:=S'''$ has all the required properties.
\end{proof}

\begin{proposition}
\label{smooth domination}
The following morphism of $B$-groupoids 
$$
P:=\underset{\stackrel{\scriptstyle N\geq n}{N'\geq 0}}{\amalg}P_{N,N'}:
\quad
\underset{\stackrel{\scriptstyle N\geq n}{N'\geq 0}}{\amalg} 
H_N \longrightarrow \GVB_n(C/B)
$$
is smooth and surjective.
\end{proposition}

\begin{proof}
Let $S$ be an affine $B$-scheme and let $S\to\GVB(C/B)$ be an $S$-valued
point of $\GVB(C/B)$, given by a Gieseker vector bundle 
$(\Cc\to C\times_BS\ ,\ \E)$ on $C$ over $S$. Let $Y$ be the $B$-groupoid
defined by the cartesian diagram
$$
\xymatrix{
Y \ar[d] \ar[r] & H_N \ar[d]^{P_{N,N'}} \\
S \ar[r] & \GVB(C/B)
}
$$
I claim that 
$Y=\Isom(\Oo^N_{S_{N,N'}},\pi_*\E(N')|_{S_{N,N'}})$,
where $\pi:\Cc\to S$ is the structure morphism and
where $S_{N,N'}\subseteq S$ is the open subscheme defined with respect to
$(\Cc\to C\times_BS\ ,\ \E)$ in lemma \ref{S_NN'}.

Let $U$ be an $S$-scheme. By definition, a $U$-valued point of $Y$
is given by the following data:
\begin{enumerate}
\item
A closed subscheme $\Cc'\injto C\times_B\Grass_n(\Oo_B^N)\times_BU$
with the properties 1.-4. listed in definition \ref{def H_N}.
\item
An isomorphism $\varphi:\Cc|_U\isomto \Cc'$ of $U$-schemes.
\item
An isomorphism $\E(N')|_U\isomto\varphi^*\E'$,
where 
$
\xymatrix{
\Oo^N_{\Cc'} \ar@{->>}[r] & \E'
}
$
is the pull-back of the canonical quotient on the Grassmannian.
\end{enumerate}
Observe that $U\to S$ factorizes over $S_{N,N'}$. Indeed this
follows easily from the characterization of $S_{N,N'}$ given in
lemma \ref{S_NN'} and the fact  (cf. \ref{pi_*Oo})
that for any $z\in S$ we have
$H^0(\Cc_z,\Oo_{\Cc_z})=\kappa(z)$.
By \ref{pi_*Oo} and \ref{def H_N}.3 we have 
$\Oo_U^N=\pi'_*\Oo_{\Cc'}^N|_U\isomto\pi'_*\E'|_U$.
Therefore the above data induce an isomorphism
$\Oo_U^N\isomto\pi_*\E(N')|_U$, i.e. a $U$-valued point
of $\Isom(\Oo^N_{S_{N,N'}},\pi_*\E(N')|_{S_{N,N'}})$.

Conversely, let $U$ be an $S_{N,N'}$-scheme and 
let $\Oo_U^N\isomto\pi_*\E(N')|_U$ be a $U$-valued point
of $\Isom(\Oo^N_{S_{N,N'}},\pi_*\E(N')|_{S_{N,N'}})$.
Then we have:
\begin{enumerate}
\item
Let $(\Cc_U\to C\times_BU\ ,\ \E_U)$ be the pull back of 
$(\Cc\to C\times_BS\ ,\ \E)$ to $U$
and let $\pi_U:\Cc_U\to U$ be the projection.
Then we have a canonical morphism
$\Oo_{\Cc_U}^N\isomto\pi_U^*(\pi_U)_*\E_U(N')\to\E_U(N')$, which by 
definition of $S_{N,N'}$ (cf. lemma \ref{S_NN'}) is surjective 
and induces a closed immersion $\Cc_U\injto\Grass_n(\Oo_B^N)\times_BU$.
Thus we have a closed immersion
$\Cc_U\injto C\times_B\Grass_n(\Oo_B^N)\times_BU$ which obviously has
the properties 1.-4. listed in \ref{def H_N}.
\item
The identity isomorphism $\Cc|_U\isomto\Cc_U$.
\item
The natural isomorphism $\E(N')|_U=\E_U(N')\isomto\E'$ where 
$\Oo^N_{\Cc_U}\to\E'$ is the pull back of the canonical quotient of
the Grassmannian.
\end{enumerate}
These are the data needed to define a morphism $U\to Y$.
Thus we have established a correspondence between $U$-valued
points of $Y$ and of $\Isom(\Oo^N_{S_{N,N'}},\pi_*\E(N')|_{S_{N,N'}})$
which is easily seen to be a bijection functorial in $U$.
This proves our claim.

Since $\Isom(\Oo^N_{S_{N,N'}},\pi_*\E(N')|_{S_{N,N'}})$
is smooth and surjective over $S_{N,N'}$, all that remains to
be shown, is that the open subschemes $(S_{N,N'})_{N\geq n, N'\geq 0}$
cover $S$. But this is an immediate consequence of 
\cite{NS}, proposition 4.
\end{proof}

\section{Gieseker vector bundle data}
\label{section GVBD}

\begin{definition}
\label{modification of X,x',x''}
Let $K$ be a field, $\tilde{X}/K$ a curve and $x',x''$ two
different smooth $K$-rational points of $\tilde{X}$.
A {\em modification of $(\tilde{X},x',x'')$}
is a morphism $(Y,y',y'')\to(\tilde{X},x',x'')$ of two-pointed
curves, such that
\begin{enumerate}
\item
$Y$ is isomorphic to the cokernel of the double arrow 
(cf. section \ref{notation})
$$
\xymatrix{
\Spec(K)\amalg\Spec(K)
\ar@<0.6ex>[r]^(.6){x'\amalg x''} \ar@<-0.6ex>[r]_(.6){x'_0\amalg x''_0} &
\text{$\tilde{X}$}\amalg R'\amalg R''
}
$$
where $R'$ (resp. $R''$) is a chain of rational curves over $K$ of 
length $r'$ (resp. $r''$)
and $x'_0, x'_{r'}$ (resp. $x''_0, x''_{r''}$) are smooth $K$ rational
points of $R'$ (resp. of $R''$), which lie on different extremal components,
if $r'\geq 2$ (resp. if $r''\geq 2$), are not equal, if $r'=1$ 
(resp. if $r''=1$), and are identical, if $r'=0$ (resp. if $r''=0$).
\item
The composed morphism $\tilde{X}\to Y\to\tilde{X}$ is the identity,
the composed morphism $R'\to Y\to\tilde{X}$ 
(resp. $R''\to Y\to\tilde{X}$) is the constant morphism onto $x'$
(resp. onto $x''$) and $y'$ (resp. $y''$) is the $K$-rational point 
induced by $x'_{r'}$ (resp. by $x''_{r''}$).
\end{enumerate}
\end{definition}
The situation may be visualized as follows:

\vspace{3mm}
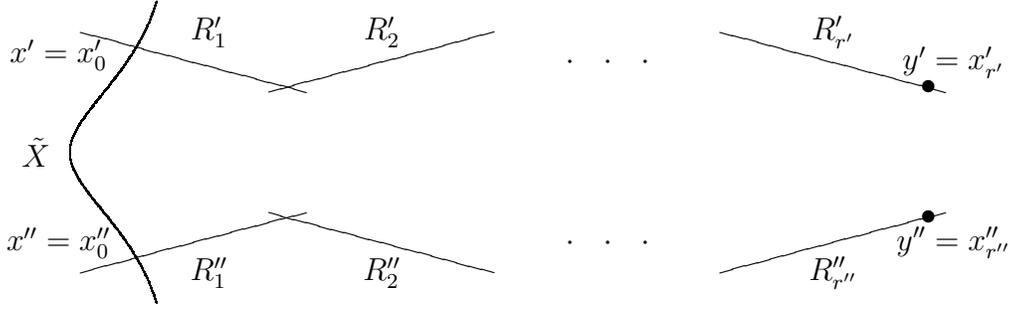
\begin{center}
\parbox{8cm}{
\xy <0mm,-20mm>; <1mm,-20mm>:
(10,20);(10,-20)**\crv{(7,10)&(-10,0)&(7,-10)};
(0,16);(30,8) **@{-}; (17,16) *{\text{$R'_1$}};
(25,8);(55,16) **@{-}; (40,16) *{\text{$R'_2$}};
(65,12)*{\cdot};
(70,12)*{\cdot};
(75,12)*{\cdot};
(85,16);(115,8) **@{-} ?>>>*{\bullet}; (100,16) *{\text{$R'_{r'}$}};
(0,-16);(30,-8) **@{-}; (17,-16) *{\text{$R''_1$}};
(25,-8);(55,-16) **@{-}; (40,-16) *{\text{$R''_{2}$}};
(65,-12)*{\cdot};
(70,-12)*{\cdot};
(75,-12)*{\cdot};
(85,-16);(115,-8) **@{-} ?>>>*{\bullet}; (100,-16) *{\text{$R''_{r''}$}};
(-6,0) *{\text{$\tilde{X}$}};
(-3,13) *{\text{$x'=x'_0$}};
(-3,-12) *{\text{$x''=x''_0$}};
(116,12) *{\text{$y'=x'_{r'}$}};
(116,-12) *{\text{$y''=x''_{r''}$}};
\endxy
}
\end{center}
\vspace{3mm}

\begin{definition}
\label{admissible data}
Let $h:(Y,y',y'')\to(\tilde{X},x',x'')$ be a modification of 
$(\tilde{X},x',x'')$ as in definition \ref{modification of X,x',x''}.
A locally free $\Oo_Y$-module $\E$ is called {\em admissible for
$(Y,y',y'')\to(\tilde{X},x',x'')$}, if the restrictions of $\E$
to the two chains of rational curves $R'=h^{-1}(x')$ and 
$R''=h^{-1}(x'')$ are both admissible in the sense of \ref{admissible}.
The {\em extremal degree} of an admissible $\Oo_Y$-module 
$\E$ of rank $n$ is the pair of numbers $(d',d'')\in[1,n+1]^2$,
where for $\iota\in\{',''\}$ the number
$d^{\iota}$  is the degree of $\E|_{R^{\iota}_{r^{\iota}}}$
if $R^{\iota}$ is of length $r^{\iota}\geq 1$, where 
$R^{\iota}_{r^{\iota}}$ 
is the component of $R^{\iota}$ which contains the point $y^{\iota}$, 
and $d^{\iota}:=\infty$ if $R^{\iota}$ is of lenght $0$.
A pair consisting of a locally free $\Oo_Y$-module $\E$ and
an isomorphism $\varphi:\E[y']\isomto\E[y'']$ is called 
{\em admissible for $(Y,y',y'')\to(\tilde{X},x',x'')$}, if
the the module induced by $(\E,\varphi)$ on the chain
$R:=(R'\amalg R'')/(y'=y'')$ of rational curves of
length $r'+r''$ 
(by identifying via $\varphi$ the fibres at the points $y'$ and $y''$) 
is admissible in the sense of definition \ref{admissible}.
\end{definition}

\begin{remark}
With the notation in the above definition, admissibility of
the pair $(\E,\varphi)$ implies admissibility of $\E$ for
$(Y,y',y'')\to(\tilde{X},x',x'')$.
Cf. the end of remark \ref{remark to admissibility}.
\end{remark}

\begin{definition}
\label{modification/S of C0,p1,p2}
Let $C_0\to \Spec(k)$ be the  curve from section \ref{notation}
and let $(\tilde{C_0},p_1,p_2)$ be the normalization together
with its two distinguished points. Let $S$ be a $k$-scheme.
A {\em modification of $(\tilde{C_0},p_1,p_2)$ over $S$} is a 
diagram 
$$
\xymatrix{
\Cc \ar@{->}[dr]^{\pi} \ar@{->}[rr]^(.4){h}    &   & 
\text{$\tilde{C_0}\times_{\Spec(k)}S$} \ar[dl]_{\pr_2}  \\
& S \ar@/^/[ul] \ar@<1ex>@/^/[ul]^{s_1,s_2} 
   \ar@/_/[ur] \ar@<-1ex>@/_/[ur]_{p_2\times\id,p_1\times\id} &   
}
$$
which is commutative in the sense that $\pr_2\comp h=\pi$ and 
$h\comp s_i=p_i\times\id$ for $i=1,2$,
where $\pi$ is flat, $h$ is proper and finitely presented, 
the $s_i$ are sections of $\pi$,
and for every 
point $z\in S$ the induced morphism on the fibre 
$(\Cc_{z},s_1(z),s_2(z))\to
(\tilde{C_0},p_1,p_2)\times_{Spec(k)}\Spec(\kappa(z))$ 
is a modification
of $(\tilde{C_0},p_1,p_2)\times_{Spec(k)}\Spec(\kappa(z))$ 
in the sense of definition \ref{modification of X,x',x''}.
\end{definition}

\begin{definition}
\label{admissible for Cc,pi,s1,s2,f}
Let $(\Cc,\pi,s_1,s_2,h)$ be modification of 
$(\tilde{C_0},p_1,p_2)$ over $S$ as in definition
\ref{modification/S of C0,p1,p2}. A locally free $\Oo_{\Cc}$-module
$\E$ is called 
{\em admissible for $(\Cc,\pi,s_1,s_2,h)$},
if for all $z\in S$ the restriction of $\E$ to the fibre $\Cc_z$ is
admissible for 
$(\Cc_{z},s_1(z),s_2(z))\to
(\tilde{C_0},p_1,p_2)\times_{Spec(k)}\Spec(\kappa(z))$
in the sense of definition \ref{admissible data}.
We define $\E$ to be of {\em extremal degree $\geq (d_1,d_2)$},
if for all $z\in S$ we have $d'_z\geq d_1$ and $d''_z\geq d_2$,
where $(d'_z,d''_z)$ is the extremal degree of $\E|_{\Cc_z}$ 
in the sense of \ref{admissible data}.
A pair consisting of a locally free $\Oo_{\Cc}$-module $\E$ and
an isomorphism $\varphi:s_1^*\E\isomto s_2^*\E$ is called
{\em admissible for $(\Cc,\pi,s_1,s_2,h)$},
if it is fibrewise admissible in the sense of \ref{admissible data}.
\end{definition}

\begin{lemma}
\label{GVBD/noetherian}
Let $(\tilde{C_0},p_1,p_2)$ be the two-pointed $k$-curve from 
section \ref{notation}
and let $S'$ be an affine $k$-scheme. Let 
$(\Cc',\pi',s'_1,s'_2,h')$ be a modification of $(\tilde{C_0},p_1,p_2)$ over 
$S'$
and let $\E'$ be a locally free $\Oo_{\Cc'}$-module which is admissible
for $(\Cc',\pi',s'_1,s'_2,h')$.
Then there exists a $k$-scheme of finite type $S$, a modification 
$(\Cc,\pi,s_1,s_2,h)$ of $(\tilde{C_0},p_1,p_2)$ over $S$,
a locally free $\Oo_{\Cc}$-module admissible for $(\Cc,\pi,s_1,s_2,h)$
and a $k$-morphism $S'\to S$ such that $(\Cc',\pi',s'_1,s'_2,h';\E')$
is induced by pull-back via $S'\to S$ from $(\Cc,\pi,s_1,s_2,h;\E)$.
\end{lemma}

\begin{proof}
We omit the proof which is similar to the proof of lemma \ref{noetherian}.   
\end{proof}

\begin{definition}
\label{GVBD}
Let $(\tilde{C_0},p_1,p_2)$ be the two-pointed $k$-curve from 
section \ref{notation}
and let $S$ be a $k$-scheme. A 
{\em Gieseker vector bundle data (of rank $n$) on 
$(\tilde{C_0},p_1,p_2)$ over $S$}
is a modification $(\Cc,\pi,s_1,s_2,h)$ of $(\tilde{C_0},p_1,p_2)$
over $S$ together with a locally free $\Oo_{\Cc}$-module $\E$ 
(of rank $n$) and an isomorphism $\varphi:s_1^*\E\isomto s_2^*\E$,
such that $(\E,\varphi)$ is admissible for $(\Cc,\pi,s_1,s_2,h)$.
\end{definition}

As all Gieseker vector bundle data of rank $n$ on $(C_0,p_1,p_2)$ over
a fixed $k$-scheme $S$ form in an obvious way a groupoid, and as there is
a natural notion of pull-back under a morphism $S'\to S$ of such data,
we have the following functor
\begin{eqnarray*}
\GVBD_n(\tilde{C_0},p_1,p_2) &:& \left\{
\begin{array}{ll}
\{\text{$k$-schemes}\}^o &\to\quad \{\text{groupoids}\} \\
\qquad S &
\mapsto\quad
\left( 
\begin{array}{ll}
\text{Gieseker vector bundle data} \\ 
\text{of rank $n$ on $(C_0,p_1,p_2)$  over $S$}
\end{array}
\right)
\end{array}
\right.
\end{eqnarray*}
Again, we often suppress the index $n$ in the notation of this $k$-groupoid.

Let $S$ be a $k$ scheme and let 
$(\tilde{\Cc},\tilde{\pi},s_1,s_2,\tilde{h},\tilde{\E},\varphi)$ 
be an object in $\GVBD(\tilde{C_0},p_1,p_2)(S)$.
Let $\Cc$ be the cokernel of the double arrow (cf. section \ref{notation})
$$
\xymatrix{
S \ar@<0.6ex>[r]^{s_1} \ar@<-0.6ex>[r]_{s_2} & \text{$\tilde{\Cc}$}
}
$$ 
and let $\pi:\Cc\to S$ the morphism induced by $\tilde{\pi}$.
By the universal property of $\Cc$ (cf. \cite{Knudsen}, Theorem 3.4),
$\tilde{h}$ induces an $S$-morphism $h:\Cc \to C_0\times_{\Spec(k)}S$.
Furthermore, the data $(\tilde{\E},\varphi)$ induces a bundle $\E$ on
$\Cc$ which clearly is admissible for $h$.
Thus we have constructed an object $(\Cc,\pi,h,\E)$ in $\GVB(C_0/k)(S)$.
Since the construction
$$
(\tilde{\Cc},\tilde{\pi},s_1,s_2,\tilde{h},\tilde{\E},\varphi)
\mapsto
(\Cc,\pi,h,\E)
$$
is functorial with respect to isomorphisms and 
commutes with base-change by a morphism $S'\to S$, we obtain a
1-morphism $\GVBD(\tilde{C_0},p_1,p_2)\to\GVB(C_0/k)$ of $k$-groupoids.

\begin{definition}
(cf. \cite{Deligne})   
Let $S$ be a scheme and let $\Y$ be a reduced algebraic $S$-stack.
An algebraic $S$-stack $\X$ together with a representable 1-morphism
$\X\to\Y$ is called the {\em normalization} of $\Y$, if 
there exists a cartesian diagram
$$
\xymatrix{
X \ar[r] \ar[d] & Y \ar[d] \\
\X \ar[r] & \Y 
}
$$
where $X$ and $Y$ are algebraic spaces, 
$Y\to\Y$ is a smooth surjective 1-morphism and
$X$ is the normalization of $Y$. 
\end{definition}

\begin{theorem}
\label{normalization}
Assume that $k$ is algebraically closed and of characteristic zero.
Then the natural morphism of $k$-groupoids 
$$
\GVBD(\tilde{C_0},p_1,p_2)\to
\GVB(C_0/k)
$$
identifies $\GVBD(\tilde{C_0},p_1,p_2)$ with the normalization
of  $\GVB(C_0/k)$.
\end{theorem}

The proof of the theorem will be given after \ref{versal} below.

\begin{lemma}
\label{normality criterion}
Let $k$ be an algebraically closed field and let $X'\to X$ be 
a proper morphism of $k$-schemes which are of finite 
type over $k$. Assume that for every closed point $x\in X$
there exist numbers $1\leq m\leq n$ and a commutative diagram
$$
\xymatrix{
(t_j,\dots,t_j) 
&
\Spec\prod_{i=1}^{m}k[[t_1,\dots,t_n]]/(t_i) \ar[r]^(.6){\isomorph} \ar[d]
&
X'\times_X\Spec\widehat{\Oo_{X,x}} \ar[d] 
\\
t_j \ar@{|->}[u]
&
\Spec k[[t_1,\dots,t_n]]/(t_1\cdot\dots\cdot t_m) \ar[r]^(.65){\isomorph}
&
\Spec\widehat{\Oo_{X,x}}
}
$$
where the horizontal arrows are isomorphisms and the left vertical arrow
is defined as indicated. Assume furthermore that $X'\to X$ is an isomorphism
outside the singular locus of $X$. Then $X$ is reduced, $X'$ is smooth and 
$X'\to X$ is the normalization of $X$.
\end{lemma}

\begin{proof}
The assumptions clearly imply that $X'$ is smooth and that $X$ is reduced. 
Let $(X'_i)_{i\in I}$ and $(X_j)_{j\in J}$ be the irreducible components
of $X'$ and $X$ respectively. Since each $X_j$ is generically smooth
and by assumption $X'\to X$ is an isomorphism outside the singular locus
of $X$, it follows that we may identify the index sets $I$ and $J$
such that $X'\to X$ induces a birational morphism $X'_i\to X_i$ for
all $i\in I$. By \cite{EGA II} (6.3.9) the morphism $X'_i\to X_i$ 
factorizes uniquely over the normalization $\tilde{X_i}\to X_i$ of $X_i$.
The induced morphism $X'_i\to\tilde{X_i}$ is 
proper, birational and quasifinite
and hence an isomorphism by \cite{EGA III} (4.4.9). Therefore 
$X'=\amalg_i X'_i\isomto \amalg_i \tilde{X_i} =\tilde{X}$ is the
normalization of $X$ (cf. \cite{EGA II} (6.3.8)).
\end{proof}

\begin{construction}
\label{versal}
Let $r\geq 0$. For $l\in[0,r]$ let 
$V_l:=\Spec(k[[t_0,\dots,t_{l-1},t_{l+1},\dots,t_r]])$.
We define modifications 
$(\tilde{Z}_l^{(i)}, \pi^{(i)}, s_1^{(i)}, s_2^{(i)}, h^{(i)})$ 
of
$(\tilde{C_0},p_1,p_2)$ over $V_l$ 
for  $i=0,\dots,l-1$ inductively as follows.
Let $\tilde{Z}_l^{(0)}:=\tilde{C_0}\times_kV_l$, the sections 
$s_\nu^{(0)}$ be induced by the points $p_\nu$, the morphism
$\pi^{(0)}$ be the projection
onto $V_l$ and $h^{(0)}$ the identity morphism.
Assume that 
$(\tilde{Z}_l^{(i)}, \pi^{(i)}, s_1^{(i)}, s_2^{(i)}, h^{(i)})$ 
has already been constructed for some $i\in[0,l-2]$.
Then $\tilde{Z}_l^{(i+1)}$ is defined to be the blowing up of 
$\tilde{Z}_l^{(i)}$
along the closed subscheme $s_1^{(i)}(\{t_i=0\})$.
The morphisms $\pi^{(i+1)}$ and $h^{(i+1)}$ are by definition
the composition of the morphism 
$\tilde{Z}_l^{(i+1)}\to\tilde{Z}_l^{(i)}$
with the morphisms $\pi^{(i)}$ and $h{(i)}$ respectively, and 
the sections $s^{(i+1)}_\nu$ are the proper transforms of the sections
$s^{(i)}_\nu$.
Now let 
$
(\tilde{Z}_l^{[0]}, \pi^{[0]}, s_1^{[0]}, s_2^{[0]}, h^{[0]}):=
(\tilde{Z}_l^{(l-1)}, \pi^{(l-1)}, s_1^{(l-1)}, s_2^{(l-1)}, h^{(l-1)})
$ 
and define
$(\tilde{Z}_l^{[i]}, \pi^{[i]}, s_1^{[i]}, s_2^{[i]}, h^{[i]})$ 
for $i:=1,\dots,r-l$ inductively as follows: For $i\geq 0$
the scheme $\tilde{Z}_l^{[i+1]}$ is the blowing up of the scheme
$\tilde{Z}_l^{[i]}$ along the closed subscheme
$s_2^{[i]}(\{t_{r-i}=0\})$. The morphisms $\pi^{[i+1]}$ and $h^{[i+1]}$
are again the composition of the blowing-up morphism with the
morphisms $\pi^{[i]}$ and $h^{[i]}$ respectively and the sections
$s^{[i+1]}_\nu$ are the proper transforms of the sections
$s^{[i]}_\nu$.
Let $Z_l$ be the cokernel of the double arrow (cf. section \ref{notation})
$$
\xymatrix{
V_l \ar@<0.6ex>[r]^{s_1} \ar@<-0.6ex>[r]_{s_2} 
& \text{$\tilde{Z}_l^{[r-l]}$}
}
$$
where $s_\nu:=s_\nu^{[r-l]}$.
Then the morphism $h^{[r-l]}$ factorizes over $Z_l$ and
the induced morphism $Z_l\to V_l$ is a modification of 
$C_0$ over $V_l$. Furthermore, the sections $s_1$ and $s_2$
induce a section $\sigma_l$ of $Z_l\to V_l$.
Observe that the special fibre of $Z_l\to V_l$ contains a chain
of projective lines of length $r$ and that $\sigma_l$ cuts this
chain into two chains of length $l$ and $r-l$ respectively.

Let $V:=\Spec(k[[t_0,\dots,t_r]]/(t_1\dots t_r))$ and let 
$V_l\injto V$ be the closed immersion defined by 
$$
\xymatrix@R=1ex{
k[[t_0,\dots,t_r]]/(t_1\dots t_r)) \ar[r] &
k[[t_0,\dots,t_{l-1},t_{l+1},\dots,t_r]] 
\\
t_i \ar@{|->}[r] & t_i \quad \text{for $i\neq l$} 
\\
t_i \ar@{|->}[r] & 0 \quad \text{for $i= l$} 
}
$$
Thus we have a morphism 
$
\tilde{V}:=\amalg_{l=0}^{r}V_l \to V
$.
Let $Z':=\amalg_{l=0}^r Z_l$ and let $Z'\to \tilde{V}$
be the morphism induced by the morphisms $Z_l\to V_l$.
Let $V'':=\tilde{V}\times_V\tilde{V}$ and let $\pr_\nu:V''\to\tilde{V}$
be the projection onto the $\nu$-th component.
Let $Z''_\nu:=Z'\times_{\tilde{V},\pr_\nu}V''$ 
be the pull-back of $Z'$ by $\pr_\nu$. Then there is an isomorphism
$Z''_1\isomto Z''_2$ over $V''$ which provides an effective descent datum
for $Z'$ relative to the morphism $\tilde{V}\to V$.
Thus there exists a $V$-scheme $Z$ such that 
$Z'\isomorph Z\times_V\tilde{V}$. The morphisms $Z_l\to C_0\times_k V_l$ 
induce a morphism $Z\to C_0\times_k V$ and it is clear that this
defines a modification of 
$C_0$ over $V$. 
Furthermore, $Z\to C_0\times_k V$ is universal in the following sense:
Let $S$ be the spectrum of a complete local $k$-algebra with residue
field $k$ and let $\Cc\to C_0\times_kS$ be a modification of $C_0$ over
$S$ such that its special fibre contains a chain of projective
lines of length $r$ (and none of length $r+1$). Then there exists
a morphism $S\to V$ such that $\Cc\isomorph Z\times_VS$ 
(cf. appendix in \cite{NS}).

Finally we remark that the composed morphisms 
$$
\xymatrix{
V_l \ar[r]^{\sigma_l} & Z_l \ar[r] & Z
}
$$
induce a closed immersion $\tilde{V}\injto Z$ which identifies
$\tilde{V}$ with the subscheme $\Sigma$ of $Z$ defined by the first
Fitting ideal of $\Omega^1_{Z/V}$.
\end{construction}

\begin{proof}(of theorem \ref{normalization})
Let $N\geq n$ and
let $H_{N,0}/k$ be the special fibre of the $B$-scheme $H_N$ defined
in \ref{def H_N}.
Let 
$\Cc_{H_{N,0}}\injto 
C_0\times_{\Spec(k)}\Grass_n(k^N)\times_{\Spec(k)}H_{N,0}$
be the universal object over $H_{N,0}$.
Let $\widetilde{H_{N,0}}\injto\Cc_{H_{N,0}}$ be the locus of points where the
morphism $\Cc_{H_{N,0}}\to H_{N,0}$ fails to be smooth.
More precisely, $\widetilde{H_{N,0}}$ is the closed subscheme of
$\Cc_{H_{N,0}}$, defined by the first Fitting ideal of the sheaf 
$\Omega^1_{\Cc_{H_{N,0}}/H_{N,0}}$.
I claim that
$\widetilde{H_{N,0}}\to H_{N,0}$ is the normalization of $H_{N,0}$.
Indeed, let $x$ be a closed point of $H_{N,0}$ and let $H_x$ be
the spectrum of the completion of the local ring of $H_{N,0}$ at $x$.
Let $\Cc\to H_x$ be the pull-back of the universal curve 
$\Cc_{H_{N,0}}$ by the natural map $H_x\to H_{N,0}$.
Let $\Sigma\injto Z\to V$ be as in 
\ref{versal} with $r\geq 0$ chosen such that
the special fibres of $\Cc\to H_x$ and $Z\to V$ be isomorphic.
Then there is a morphism $H_x\to V$ such that $\Cc=Z\times_VH_x$
and it follows that 
$\widetilde{H_{N,0}}\times_{H_{N,0}}H_x=\Sigma\times_VH_x$
Since furthermore the morphism $H_x\to V$ is formally smooth 
(cf. \cite{NS}, appendix), the claim now follows from lemma
\ref{normality criterion} (observe that $H_{N,0}$ is the disjoint union
of schemes of finite type over $k$).

It suffices now to show that there exists a morphism 
$\widetilde{H_{N,0}}\to\GVBD(\tilde{C_0},p_1,p_2)$ such that the following
diagram of $k$-groupoids is cartesian:
$$
\xymatrix{
\text{$\widetilde{H_{N,0}}$} \ar[d]\ar[r] & H_{N,0} \ar[d] \\
\text{$\GVBD(\tilde{C_0},p_1,p_2)$} \ar[r] & \GVB(C_0/k)
}
$$
But this is easy: Let $S$ be an affine $k$-scheme. An $S$-valued point
of $\widetilde{H_{N,0}}$ is an $S$-valued point 
$\Cc\injto C_0\times_k\times\Grass_n(k^n)\times_kS$
of $H_{N,0}$ plus a section $S\to\Cc$ of the projection $\pi:\Cc\to S$,
meeting $\Cc$ in the nonsmooth locus of $\pi$. Observe that 
this is exactly the description
of an $S$-valued point of the fibre product
$\GVBD(\tilde{C_0},p_1,p_2)\times_{\GVB(C_0/k)}H_{N,0}$.
\end{proof}

\section{Simple modifications of a relative curve with a smooth section}
\label{section simple modifications}

\begin{definition}
\label{def mpc}
Let $S$ be a scheme and 
$\pi:\Cc\to S$ a flat morphism whose geometric fibers
are connected and of dimension one. 
Let $s: S\to \Cc$ be a section of $\pi$ whose image
lies in the subset of $\Cc$ where $\pi$ is smooth.
A {\em simple modification $(\Cc',f,\pi',s')$ of $(\Cc,\pi,s)$} 
is a diagram
$$
\xymatrix{
\Cc' \ar[rr]^f \ar[dr]^{\pi'} & & \Cc \ar[dl]_{\pi} \\
   & S \ar@/^/[ul]^{s'} \ar@/_/[ur]_s
}
$$
with the following properties:
\begin{enumerate}
\item
The morphism $\pi'$ is again flat with connected 
geometric fibers of dimension one and $s'$ is a section of $\pi'$ whose 
image consists of points where $\pi'$ is smooth.
\item
The diagram is commutative in the sense that $\pi\comp f=\pi'$
and $f\comp s'=s$.
\item 
The morphism $f$ is proper and finitely presented.
\item
Let $z\in S$ be a point. Then there are two possibilities for
the induced morphism $f_z:\Cc'_z\to\Cc_z$ of fibres over $z$:
Either $f_z$ is an isomorphism, or $\Cc'_z$ is the quotient
of $R\isomorph\Pp_{\kappa(z)}^1$ and $\Cc_z$ by the identification
of a point in $R(\kappa(z))$ with $s(z)\in\Cc_z$
and $f_z$ contracts $R$ to the point $s(z)$.
\end{enumerate}
\end{definition}

\begin{construction}
\label{constr1}
Let $(\Cc',f,\pi',s')$ be a simple modification of $(\Cc,\pi,s)$ 
as in definition \ref{def mpc} and let
$\M:=\Oo_{\Cc'}(-s')\tensor f^*\Oo_{\Cc}(s)$.
We have the following exact diagram of $\Oo_{\Cc'}$-modules :
$$
\xymatrix@R=3ex{
0 \ar[r] &
\Oo_{\Cc'}(-s') \ar[r] \ar[d] &
\Oo_{\Cc'} \ar[r] \ar@{.>}[dl] \ar[d] &
s'_*s^*\Oo_{\Cc} \ar[r] \ar[d]^0 &
0 \\
0 \ar[r]&
\M \ar[r] &
f^*\Oo_{\Cc}(s) \ar[r] &
s'_*s^*\Oo_{\Cc}(s) \ar[r] & 
0 
}
$$
where the vertical arrows are induced by the
morphism $\Oo_{\Cc'}\to f^*\Oo_{\Cc}(s)$ obtained
by applying $f^*$ to the natural injection $\Oo_{\Cc}\injto\Oo_{\Cc}(s)$.
Since the right vertical arrow obviously vanishes,
the middle vertical arrow factorizes as indicated by the dotted
arrow. 
Applying $(s')^*$ to the morphism 
$m:\Oo_{\Cc'}\to \M$ thus obtained, yields a morphism
$$
\mu_1:\Oo_S\to M_1:=(s')^*\Oo_{\Cc'}(-s')\tensor s^*\Oo_{\Cc}(s)
\quad.
$$
\end{construction}

\begin{proposition}
\label{mod to Mmu}
Let $S$ be an arbitrary scheme, $\pi:\Cc\to S$ a proper, finitely 
presented, flat morphism 
with connected geometric fibers of dimension one 
and let 
$s:S\to\Cc$ be a section of $\pi$ whose image consists of smooth points of 
$\pi$. 
Assume that for every point $z\in S$ we have 
$H^0(\Cc_z,\Oo_{\Cc_z})=\kappa(z)$ and
$H^0(\Cc_z,\Oo_{\Cc_z}(-s(z)))=(0)$ where $\Cc_z$ denotes
the fibre of $\pi$ at $z$.
Then the above construction \ref{constr1} 
yields an isomorphism of groupoids:
$$
\left\{
\begin{array}{ll}
\text{Simple modifications of $(\Cc,\pi,s)$} \\
\text{in the sense of definition \ref{def mpc}}
\end{array}
\right\}
\isomto
\left\{
\begin{array}{llll}
\text{Pairs $( M,\mu)$, where $ M$ is } \\
\text{an invertible $\Oo_S$-module and}\\
\text{$\mu:\Oo_S\to  M$ is a global}\\
\text{section of $ M$}
\end{array}
\right\}
$$
\end{proposition}

In the remaining of this section we prove the above proposition.
For this, we fix a scheme $S$ and data $(\Cc,\pi,s)$
as in the proposition.

\begin{construction}
\label{constr2}
To a pair $( M,\mu)$, consisting of an invertible $\Oo_S$-module $ M$
together with a global section $\mu$ we associate a diagram
$$
\xymatrix{
\Cc'_1 \ar[rr]^{f_1} \ar[dr]^{\pi'_1} & & \Cc \ar[dl]_{\pi} \\
   & S \ar@/^/[ul]^{s'_1} \ar@/_/[ur]_s
}
$$
as follows:

Let $\I$ be the $\Oo_{\Cc}$-module defined by the exact sequence
$$
0 \to
\pi^* M^{\dual}\tensor\Oo_{\Cc}(-s) \to
\pi^* M^{\dual}\oplus\Oo_{\Cc}(-s) \to
\I \to
0
\quad,
$$
where 
$\pi^* M^{\dual}\tensor\Oo_{\Cc}(-s) \to \pi^* M^{\dual}$
is induced by the natural injection 
$\Oo_{\Cc}(-s) \to \Oo_{\Cc}$
and 
$\pi^* M^{\dual}\tensor\Oo_{\Cc}(-s) \to \Oo_{\Cc}(-s)$
is induced by $-\mu$.
We define 
$$
\Cc'_1:=\Proj(\Sym\I) \quad.
$$
Let $f_1: \Cc'_1\to\Cc$ be the projection and $\pi'_1:=\pi\comp f_1$.
The morphism $s'_1:S\to\Cc'_1$ is the one induced by the invertible quotient
$s^*\I\to M^{\dual}$ of $s^*\I$
which is adjoint to the arrow $\I\to s_* M^{\dual}$ defined
by the exact diagram
$$
\xymatrix@R=3ex{
& & 0\ar[d] & 0\ar[d] & \\
& &
\Oo_{\Cc}(-s) \ar[r]^{\text{id}} \ar[d]&
\Oo_{\Cc}(-s) \ar[r] \ar[d] & 0 \\
0 \ar[r] &
\pi^* M^{\dual}\tensor\Oo_{\Cc}(-s) \ar[r] \ar[d]_{\text{id}} &
\pi^* M^{\dual}\oplus\Oo_{\Cc}(-s) \ar[r] \ar[d] &
\I \ar[r] \ar[d] & 0 \\
0 \ar[r] &
\pi^* M^{\dual}\tensor\Oo_{\Cc}(-s) \ar[r] \ar[d] &
\pi^* M^{\dual} \ar[r] \ar[d] &
s_* M^{\dual} \ar[r] \ar[d] & 0 \\
& 0 & 0 & 0 &
}
$$
\end{construction}

\begin{lemma}
\label{local}
Let $(\Cc,\pi,s)$ be as above.
Assume that $S=\Spec A$ and $\Cc=\Spec B$ for a ring
$A$ and an $A$-algebra $B$. Let $b\in B$ be a regular element that defines
the section $s:S\to \Cc$. In particular, 
we have $B/(b)=A$ as $A$-algebras.
Let $a\in A$ and let $\I$ be the $B$-module associated 
as in construction \ref{constr2} to the pair $(M,\mu):=(A,a)$.
Then we have an isomorphism of graded $B$-algebras
$$
\Sym(\I)\isomorph B[X,Y]/(bX-aY)
$$
and $\Sym(\I)$ is flat over $A$.
\end{lemma}

\begin{proof}
By definition we have an exact sequence of $B$-modules
$$
\xymatrix@R=.2ex{
0 \ar[r] &
B \ar[r] &
B\oplus B \ar[r] &
\I \ar[r] &
0 \\
&
1 \ar@{|->}[r] &
(b,-a) & & &
}
$$
By \cite{Algebre}, III.2, Prop.4, this induces an exact sequence
of graded $B$-modules
$$
\xymatrix@R=.2ex{
\Sym(B\oplus B)(-1) \ar[r] &
\Sym(B\oplus B) \ar[r] &
\Sym(\I) \ar[r] &
0 \\
f \ar@{|->}[r] &
(b,-a)\cdot f & & &
}
$$
This proves the required isomorphism.
The flatness of $\Sym(\I)$ now follows from the fact that the image 
of $b$ in $B/(\p)$ and
thus the image of $bX-aY$ in $B/(\p)[X,Y]$ is regular
for any prime ideal $\p\subset A$ (cf. \cite{Matsumura}, 22.6).
\end{proof}

\begin{lemma}
The diagram 
$$
\xymatrix{
\Cc'_1 \ar[rr]^{f_1} \ar[dr]^{\pi'_1} & & \Cc \ar[dl]_{\pi} \\
   & S \ar@/^/[ul]^{s'_1} \ar@/_/[ur]_s
}
$$
constructed in \ref{constr2} is a simple modification 
of $(\Cc,\pi,s)$. Furthermore, construction \ref{constr2} commutes
with base change $S'\to S$.
\end{lemma}

\begin{proof}
It is clear 
that $f_1$ is proper and that
$\pi\comp f_1=\pi'_1$,\ \ $f_1\comp s'_1=s$,\ \ $\pi'_1\comp s'_1=\id_S$.
From the exact sequence
$
0\to\Oo_{\Cc}(-s)\to\I\to s_*M^{\dual}\to 0
$
it follows that $\I$ is invertible outside $s(S)$, thus $f_1$ is an 
isomorphism over the complement of the section $s$.
The flatness of $\pi'_1$ now follows easily from lemma \ref{local}.
Let $x\in s(S)$. By lemma \ref{local} we have
$\I_x\isomorph\Oo_{\Cc,x}[X,Y]/(bX-aY)$ 
for certain elements $a\in\Oo_{S,\pi(x)}$,
$b\in\Oo_{\Cc,x}$. Therefore, if $f_1$ is not an isomorphism in any
neighbourhood of $x$ then both $a$ and $b$ are in the maximal ideal
of $\Oo_{\Cc,x}$. This implies $f^{-1}(x)\isomorph\Pp^1_{\kappa(x)}$.
It is immediate that construction \ref{constr2} commutes with base
change. An easy calculation in the case $S=\Spec(k)$, with $k$ algebraically
closed, shows that the geometric 
fibers of $\pi'_1$ are connected and
of dimension one and that the image of $s'_1$ consists of smooth
points of $\pi'_1$.
\end{proof}

We want to show that under the assumptions of proposition
\ref{mod to Mmu},
construction \ref{constr2} is the inverse to
construction \ref{constr1}.
For this, we start with a simple modification 
$(\Cc',f,\pi',s')$ of $(\Cc,\pi,s)$, apply construction
\ref{constr1} to get a pair $( M_1,\mu_1)$ and then apply construction
\ref{constr2} to $( M,\mu):=( M_1,\mu_1)$ to get  a 
sheaf $\I$
and a simple modification
$(\Cc'_1,f_1,\pi'_1,s'_1)$  of $(\Cc,\pi,s)$,
where $\Cc'_1=\Proj(\Sym(\I))$. 
We have to show the existence of a canonical isomorphism
$\Cc'\isomto\Cc'_1$ which is compatible with the respective projections
and sections. This is done in lemma \ref{C' isomto C} below.

We need some preliminary results.
Let $\M:=\Oo_{\Cc'}(-s')\tensor  f^*\Oo_{\Cc}(s)$.
Since $(s')^*\Md=M^{\dual}$, we have a canonical morphism
$\Md\to s'_*M^{\dual}$. Application of $\pi'_*$ yields
a morphism $\pi'_*\Md \to M^{\dual}$.

\begin{lemma}
\label{piMM to M}
Let $(\Cc,\pi,s)$ satisfy the assumptions of proposition
\ref{mod to Mmu} and let $(\Cc',f,\pi',s')$ be a simple
modification of $(\Cc,\pi,s)$. Then the 
canonical morphism $\pi'_*\Md \to M^{\dual}$
is an isomorphism.
\end{lemma}

\begin{proof}
The proof of this lemma will be given in section \ref{pol}.
\end{proof}

\begin{lemma}                
\label{exact}
The following sequence of $\Oo_{\Cc'}$-modules is exact
$$
0 \to
(\pi')^* M^{\dual}\tensor\Oo_{\Cc'}(-s') \to
(\pi')^* M^{\dual}\oplus f^*\Oo_{\Cc}(-s) \to
\Md\to 0
\quad ,
$$
where the involved morphisms are defined as follows:
\begin{itemize}
\item
$(\pi')^* M^{\dual}\tensor\Oo_{\Cc'}(-s') \to (\pi')^* M^{\dual}$
is induced by $\Oo_{\Cc'}(-s')\injto\Oo_{\Cc'}$.
\item
$(\pi')^* M^{\dual}\tensor\Oo_{\Cc'}(-s') \to f^*\Oo_{\Cc}(-s)$
is induced by the negative of the adjoint of $M^{\dual}\isomto\pi'_*\Md$
(cf. lemma \ref{piMM to M}).
\item
$(\pi')^*M^{\dual}\to \Md$
is the adjoint of $M^{\dual}\isomto\pi'_*\Md$
(cf. lemma \ref{piMM to M}).
\item
$f^*\Oo_{\Cc}(-s)\to\Md$
is induced by $\Oo_{\Cc'}\injto\Oo_{\Cc'}(s')$.
\end{itemize}
\end{lemma}

\begin{proof}
First of all, the sequence under consideration is a complex, as
its middle part is of the shape
$$
\xymatrix{
N\tensor L_1 \ar[r]^{(u,-v)} &
N\oplus L_2 \ar[r]^{u+v} &
L_1^{\dual}\tensor L_2
}
\quad,
$$
where $N$, $L_1$, $L_2$ are invertible sheaves and
$u:L_1\to\Oo$ and $v:N\to L_1^{\dual}\tensor L_2$
are morphisms.
Since all the sheaves occuring in the sequence are locally free
and the rank of the sheaf in the middle is the sum of the ranks
of the two other sheaves, it suffices to show for each point
$x\in\Cc'$ the injectivity of $\alpha_x$ and the surjectivity
of $\beta_x$, where
\begin{eqnarray*}
\alpha_x &:&
((\pi')^* M^{\dual}\tensor\Oo_{\Cc'}(-s'))[x] \to
((\pi')^* M^{\dual}\oplus f^*\Oo_{\Cc}(-s))[x] \\
\beta_x &:&
((\pi')^* M^{\dual}\oplus f^*\Oo_{\Cc}(-s))[x] \to
\Md[x]
\end{eqnarray*}
are the induced morphisms on the fibers of the sheaves at $x$.
In particular, we can assume that $S=\Spec(k)$ for a field $k$.
If $f$ is an isomorphism, the statement is easy to see.
Otherwise, $\Cc'=R\cup\Cc$ and $R\cap\Cc=\{s\}$, where $R\isomorph\Pp_k^1$
and three cases are possible:
Either
$x\in\Cc$,
or
$x\in R\setminus\{s,s'\}$,
or
$x=s'$.
It is not difficult to check the statement in each of these cases.
\end{proof}

\begin{lemma}
\label{f_*}
\begin{enumerate}
\item \label{f_*.one}
We have $R^1f_*\Oo_{\Cc'}(-s') = 0$.
\item \label{f_*.two}
The adjunction morphism
$\Oo_{\Cc} \to f_*\Oo_{\Cc'}$
is an isomorphism.
\item \label{f_*.three}
The morphism 
$\Oo_{\Cc}(-s) \to f_*\Oo_{\Cc'}(-s')$
adjoint to the morphism
$f^*\Oo_{\Cc}(-s) \to \Oo_{\Cc'}(-s')$
(cf. construction \ref{constr1})
is an isomorphism.
\item \label{f_*.four}
We have $\pi_*\Oo_{\Cc}=\pi'_*\Oo_{\Cc'}=\Oo_S$.
\end{enumerate}
\end{lemma}

\begin{proof}
Since obviously $H^1(f^{-1}(x),\Oo_{\Cc'}(-s')\tensor_{\Oo_{\Cc}}\kappa(x))=0$
for every point $x\in\Cc$, Corollary 1.5 in \cite{Knudsen} implies
the first assertion in the lemma.
Similarly, it follows from loc. cit. that
$f_*\Oo_{\Cc'}$ is flat over $S$ and that for every $z\in S$ we have
$f_*\Oo_{\Cc'}\tensor_{\Oo_S}\kappa(z)=(f_z)_*\Oo_{\Cc'_z}$,
where $f_z:\Cc'_z\to\Cc_z$ is the induced morphism between the fibers 
over $z$. With the help of Nakayama's lemma, this implies the second
assertion. For the third assertion it suffices to prove the
commutativity of the diagram
$$
\xymatrix{
0 \ar[r] &
\Oo_{\Cc}(-s) \ar[d] \ar[r] &
\Oo_{\Cc} \ar[d]_{\isomorph} \ar[r] &
s_*\Oo_S \ar[d]_{=} \ar[r] &
0 
\\
0 \ar[r] &
f_*\Oo_{\Cc'}(-s') \ar[r] &
f_*\Oo_{\Cc'} \ar[r] &
s_*\Oo_S \ar[r] &
0 
}
$$
where the middle vertical arrow is the isomorphism from 2. and
the lower sequence comes from applying $f_*$ to the exact sequence
$0\to\Oo_{\Cc'}(-s')\to\Oo_{\Cc'}\to s'_*\Oo_S\to 0$.
The right square in this diagram is obviously commutative.
The commutativity of the left square is equivalent to the
commutativity of the adjoint square
$$
\xymatrix{
f^*\Oo_{\Cc}(-s) \ar[d] \ar[r] &
f^*\Oo_{\Cc}=\Oo_{\Cc'} \ar[d]_{=} 
\\
\Oo_{\Cc'}(-s') \ar[r] &
\Oo_{\Cc'}
}
$$
But this square commutes by construction \ref{constr1}.
By \cite{EGA III} (7.8.8), 
the equality $\pi_*\Oo_{\Cc}=\Oo_S$ follows from the
assumption $H^0(\Cc_z,\Oo_{\Cc_z})=\kappa(z)$ for all $z$.
The assertion $\pi'_*\Oo_{\Cc'}=\Oo_S$ follows now from 2.
\end{proof}

\begin{lemma}
\label{I=f_*M^dual}
There exists a canonical isomorphism $\I\isomto f_*\Md$.
\end{lemma}

\begin{proof}
Consider the following diagram of $\Oo_{\Cc}$-modules:
$$
\xymatrix{
0 \ar[r] &
\pi^* M^{\dual}\tensor\Oo_{\Cc}(-s) 
\ar[d]_{\isomorph} \ar[r] &
\pi^* M^{\dual}\oplus\Oo_{\Cc}(-s) 
\ar[d]_{\isomorph} \ar[r] &
\I \ar@{.>}[d] \ar[r] &
0
\\
0 \ar[r] &
\pi^* M^{\dual}\tensor f_*\Oo_{\Cc'}(-s') \ar[r] &
(\pi^* M^{\dual}\oplus \Oo_{\Cc}(-s))\tensor f_*\Oo_{\Cc'} \ar[r] &
f_*\Md \ar[r] &
0
}
$$
where the first row is the defining exact sequence for $\I$
(cf. construction \ref{constr2}) and the second row  comes
from applying $f_*$ to the exact sequence of lemma \ref{exact}.
It is exact by lemma \ref{f_*}.\ref{f_*.one}.
The left vertical arrow is induced by the isomorphism in
\ref{f_*}.\ref{f_*.three}, the middle vertical arrow by the one in 
\ref{f_*}.\ref{f_*.two}.
The lemma follows, if we can show the commutativity of the left square
in this diagram. 
We give the details, since at some point we make use of the
equality $\pi_*\Oo_{\Cc}=\Oo_S$ and thus of the assumption
$H^0(\Cc_z,\Oo_{\Cc_z})=\kappa(z)$.
The proof of the commutativity of the left square in the above diagram
amounts to showing that the two squares
$$
\xymatrix{
\Oo_{\Cc}(-S) \ar[r] \ar[d]^{\isomorph} &
\Oo_{\Cc} \ar[d]^{\isomorph} &
\pi^*M^{\dual}\tensor\Oo_{\Cc}(-s) \ar[r] \ar[d]^{\isomorph} &
\Oo_{\Cc}(-s) \ar[d]^{\isomorph}
\\
f_*\Oo_{\Cc'}(-s') \ar[r] &
f_*\Oo_{\Cc'} &
f_*((\pi')^*M^{\dual}\tensor\Oo_{\Cc'}(-s')) \ar[r] &
f_*f^*\Oo_{\Cc}(-s)
}
$$
are both commutative. For the left square this has been done already
in the proof of lemma \ref{f_*}.
The commutativity of the right square is equivalent to the commutativity
of the adjoint square
$$
\xymatrix{
(\pi')^*M^{\dual}\tensor f^*\Oo_{\Cc}(-s) \ar[r] \ar[d] &
f^*\Oo_{\Cc}(-s) \ar[d]^{=}
\\
(\pi')^*M^{\dual}\tensor \Oo_{\Cc'}(-s') \ar[r]  &
f^*\Oo_{\Cc}(-s)
}
$$
which in turn is equivalent to the commutativity of
$$
\vcenter{
\xymatrix{
(\pi')^*M^{\dual}\tensor\Md \ar[r] \ar[d] &
\Md \ar[d]^=
\\
(\pi')^*M^{\dual} \ar[r] &
\Md
}}
\eqno(*)
$$
where the upper horizontal and the left vertical arrows are
induced by $\Oo_S\to M$ and $\Oo_{\Cc'}\to \M$ 
respectively, which were constructed in \ref{constr1},
and the lower horizontal arrow is the adjoint of the inverse of the
isomorphism 
$\pi'_*\Md\isomto M^{\dual}$ from
lemma \ref{piMM to M}.      
To prove the commutativity of $(*)$, observe first that
by general nonsense the following diagram commutes:
$$
\vcenter{
\xymatrix{
\Hom(M^{\dual},\pi'_*\Md) \ar[d] \ar[r]^{\isomorph} &
\Hom((\pi')^*M^{\dual},\Md) \ar[d] \ar[dr]
\\
\Hom(M^{\dual},\pi'_*\Oo_{\Cc'}) \ar[r]^{\isomorph} &
\Hom((\pi')^*M^{\dual},\Oo_{\Cc'}) \ar[r]^(.45){\isomorph} & 
\Hom((\pi')^*M^{\dual}\tensor\Md,\Md) 
}}
\eqno(**)
$$
where the vertical arrows and the oblique arrow 
are all induced by $\Oo_{\Cc'}\to \M$ from \ref{constr1}.
The commutativity of $(*)$ follows, if we can show that
the left vertical arrow of $(**)$ maps the inverse of the
isomorphism 
$\pi'_*\Md\isomto M^{\dual}$ from
lemma \ref{piMM to M}
to the morphism $M^{\dual}\to\pi'_*\Oo_{\Cc'}$ obtained by composing
$\mu:M^{\dual}\to\Oo_S$ from \ref{constr1} with the canonical morphism
$a:\Oo_S\to\pi'_*\Oo_{\Cc'}$. This amounts to showing the commutativity of
the diagram
$$
\vcenter{
\xymatrix{
\pi'_*\Md \ar[r]^{\ref{piMM to M}} \ar[d] &
M^{\dual} \ar[d] 
\\
\pi'_*\Oo_{\Cc'} &
\Oo_S \ar[l]_a
}}
\eqno(***)
$$
Now consider the two squares
$$
\xymatrix{
\Md \ar[r] \ar[d] &
s'_*M^{\dual} \ar[d] & &
\pi'_*\Md \ar[r]^{\ref{piMM to M}} \ar[d] &
M^{\dual} \ar[d]
\\
\Oo_{\Cc'} \ar[r] &
s'_*\Oo_S & &
\pi'_*\Oo_{\Cc'} \ar[r]^b &
\Oo_S
}
$$
The left square is obviously commutative, therefore so is the right
square, which is the image of the left one under the functor $\pi'_*$.
Since by lemma \ref{f_*}.\ref{f_*.four} the morphisms $a$ and $b$ are 
inverse to each other (this is the point, where the equality
$\pi_*(\Oo_{\Cc})=\Oo_S$ enters), the commutativity of $(***)$ follows and
we are done.
\end{proof}

By \cite{Knudsen} Cor. 1.5, the morphism $f^*\I\to\Md$
adjoint to the isomorphism $\I\isomto f_*\Md$
is an epimorphism.   

\begin{lemma}         
\label{C' isomto C}
The morphism $\Cc'\to\Cc'_1$  induced by the invertible quotient
$f^*\I\to\Md$
of $f^*\I$ defined above is an isomorphism
of $\Cc$-schemes. Furthermore, the following diagram commutes:
$$
\xymatrix@R=0.5ex{
\Cc' \ar[rr] & & \Cc'_1 \\
& S \ar[ul]^{s'} \ar[ur]_{s'_1} &
}
$$
\end{lemma}

\begin{proof}
It is clear that $\Cc'\to\Cc'_1$ is a morphism of $\Cc$-schemes.
To show that this is in fact an isomorphism, we first assume
that $S=\Spec(k)$ for some field $k$.
If $f$ is an isomorphism, then $\I$ is invertible by 
\ref{I=f_*M^dual} and therefore $f_1:\Cc'_1\to\Cc$ and
$\Cc'\to\Cc'_1$ are isomorphisms.
Otherwise, $\Cc'$ is the glueing together of $R\isomorph\Pp_k^1$
and $\Cc$ by identifying a point in $R(k)$ with the point $s\in\Cc(k)$.
Since both $f$ and $f_1$ are isomorphisms over the complement
of $s$, may assume that $\Cc$ is the Spectrum of a discrete valuation ring
with closed point $s$.
It is now easy to see that the surjection
$
f^*\I\isomto f^*f_*\Md\to\Md
$
induces an isomorphism
$$
\Cc'\isomto
\Proj(\Sym(f^*f_*\Md))
\isomto
\Proj(\Sym(\I))=\Cc'_1
\quad.
$$
Now let $S$ be arbitrary.
Since both $\Cc'$ and $\Cc'_1$ are flat over $S$ and all constructions
commute with base change, it follows from the special case
that  $\Cc'\to\Cc'_1$ is an isomorphism in general.
It remains to show that $s'_1$ is the composition of $s'$ with
$\Cc'\to\Cc'_1$. For this, it suffices to show that the morphism
$s^*\I\to M^{\dual}$ which one obtains by applying $(s')^*$ to
the morphism $f^*\I\to\Md$ coincides with the morphism
$s^*\I\to M^{\dual}$ defining $s'_1$ (cf. construction \ref{constr2}).
This follows formally by the yoga of pairs of adjoint morphisms. 
\end{proof}

Conversely, let a pair $(M,\mu)$ consisting of an invertible $\Oo_S$-module
$M$ and a global section $\mu$ of $M$ be given. Applying construction
\ref{constr2} to this pair, we get an $\Oo_{\Cc}$-module $\I$ and a 
simple modification 
$(\Cc'_1,f_1,\pi'_1,s'_1)$ of $(\Cc,\pi,s)$, where $\Cc'_1=\Proj(\Sym(\I))$. 
Construction \ref{constr1} applied to the simple modification
$(\Cc',f,\pi',s'):=(\Cc'_1,f_1,\pi'_1,s'_1)$ yields a pair $(M_1,\mu_1)$.
We have to establish a canonical isomorphism $M\isomto M_1$ which
makes the respective global sections correspond.

Let $\Sc$ be the graded $\Oo_{\Cc}$-module $\Sym(\I)$.
By definition, the section $s':S\to\Cc'$
is given by applying $\Proj$ to the surjection 
$\Sc\to\Sym(s_*M^{\dual})$ which is induced by the epimorphism 
$\I\to s_*M^{\dual}$ defined in construction \ref{constr2}.
Therefore, we have
$
(s')^*\Oo_{\Cc'}(1)
=\left(
\Sc(1)\tensor_{\Sc}\Sym(s_*M^{\dual})
\right)^{\sim}
$,
where $\Oo_{\Cc'}(1)$ denotes the tautological invertible sheaf
on $\Cc'=\Proj(\Sym(\I))$ (cf. \cite{EGA II} (3.2.5.1)) and
where $^\sim$ means taking the quasicoherent sheaf on 
$S=\Proj(\Sym(s_*M^{\dual}))$
associated to a graded $\Sym(s_*M^{\dual})$-module.
It is easy to see that the sheaf 
$
\left(
\Sc(1)\tensor_{\Sc}\Sym(s_*M^{\dual})
\right)^{\sim}
$
is just $M^{\dual}$.
Thus we have 
$$
(s')^*\Oo_{\Cc'}(1)=M^{\dual}
\quad.
$$

From the exact sequence of $\Oo_{\Cc}$-modules
$
0 \to \Oo_{\Cc}(-s) \to \I \to s_*M^{\dual} \to 0
$
(cf. construction \ref{constr2})
we get by \cite{Algebre} III.2 Prop.4
the exact sequence of graded $\Oo_{\Cc}$-modules
$$
\Sc(-1)\tensor_{\Oo_{\Cc}}\Oo_{\Cc}(-s) \to \Sc \to \Sym(s_*M^{\dual}) \to 0
\quad.
$$
From lemma \ref{local} it follows easily that the left arrow is
injective.
Since obviously  
$
\Oo_{\Cc'}(-s')
$
is the $\Oo_{\Cc'}$-module associated to the graded
$\Sc$-module $\ker(\Sc \to \Sym(s_*M^{\dual}))$,
we obtain a canonical isomorphism
$$ \label{O(-1)=M}
\Oo_{\Cc'}(s')\tensor f^*\Oo_{\Cc}(-s) \isomto \Oo_{\Cc'}(1)
\quad.
$$
Applying $(s')^*$ to it yields the required isomorphism
$M\isomto M_1$.

\begin{lemma}              
The isomorphism $M\isomto M_1$ maps
the global section $\mu$ of $M$ to the global section $\mu_1$ of $M_1$.
\end{lemma}

\begin{proof}
Since both construction \ref{constr1} and construction \ref{constr2}
commute with base change $S'\to S$ and since furthermore the 
definition of the isomorphism $M\isomto M_1$ depends only on
a neighbourhood of the section $s$ in $\Cc$, we may assume
$S=\Spec(A)$ and $\Cc=\Spec(B)$ for a local ring $A$ and local 
$A$-algebra $B$. Let $l: M\isomto A$ be an isomorphism.
Let $a:=l(\mu)$ and $b$ a generator of the kernel $N$
of the surjection $B\to A$ associated to the section $s: S\to \Cc$.
By lemma \ref{local} we may identify the graded $B$ module 
$\Sc$ with $B[X,Y]/(bX-aY)$.
Let $\J$ be the graded $\Sc$-module $\ker(\Sc\to \Sym(M^{\dual}))$.
It is generated by the element $Y\in \Sc$.
By definition, we have
$$
M_1=(\J\tensor_{\Sc}\Sym(M^{\dual}))_{(l)}\tensor_A N^{\dual}
\quad,
$$
where the index $(l)$ means taking the zero-degree part
of the localization by powers of $l\in \Sym(M^{\dual})$ 
(cf. \cite{EGA II} (2.2.1)).
It is easy to check that the isomorphism $M\to M_1$ is given
by
$$
\xi\mapsto \frac{Y\tensor l(\xi)}{l}\tensor b^{\dual}
\quad,
$$
where $b^{\dual}\in N^{\dual}$ denotes the generator dual to $b$.
On the other hand, it follows directly from the definitions that
$$
\mu_1=\frac{bX\tensor 1}{l}\tensor b^{\dual}
\quad.
$$
Therefore the statement
follows from the equality $bX=aY$ which holds in $\Sc$.

\end{proof}

This completes the proof of proposition \ref{mod to Mmu}.

\section{Proof of lemma \ref{piMM to M}}
\label{pol}

\begin{lemma}
\label{S=Spec k}
Lemma \ref{piMM to M} is true in the case where $S$ is the
spectrum of a field $k$.
\end{lemma}

\begin{proof}
If $f$ is an isomorphism, then $\Md=\Oo_{\Cc}$ and the
lemma follows directly from the assumption
$H^0(\Cc,\Oo_{\Cc})=k$.
Otherwise, $\Cc'$ is the union of a projective line $R\isomorph \Pp^1$
and the curve $\Cc$, and we have $R\cap\Cc=\{s\}$. 
Furthermore, the restriction of $f$ 
to $R$ is the constant map onto the point $s\in\Cc$ and the restriction
of $f$ to $\Cc\subset\Cc'$ is the identity on $\Cc$.
Now we have 
\begin{eqnarray*}
H^0(\Cc',\Md) &=& 
H^0(R,\Md|_R) \times_{\Md[s]} H^0(\Cc,\Md|_{\Cc})\\
&=& H^0(R,\Oo_R(s')\tensor_k s^*\Oo_{\Cc}(-s)) \times_{\Md[s]} (0)\\
&=& H^0(R,\Oo_R(s'-s))\tensor_k s^*\Oo_{\Cc}(-s) 
\end{eqnarray*}
and  the morphism 
$H^0(\Cc',\Md)\to M^{\dual}$
is induced by the isomorphism
$H^0(R,\Oo_R(s'-s))\isomto(s')^*\Oo_R(s')=(s')^*\Oo_{\Cc'}(s')$
(restriction to the point $s'$).
\end{proof}

\begin{lemma}
\label{universal}
Let $k$ be a field and let 
$\tilde{f}_{vers}:\tilde{\Cc}'_{\vers}\to\tilde{\Cc}_{\vers}$ be the 
blowing up of $\tilde{\Cc}_{\vers}:=\Spec(k[[t]][[z]])$
in the closed point $t=z=0$. Let 
$\tilde{\Cc}'_{\vers}\to\Spec(k[[t]])$ be the obvious 
projection morphism. Denote by $\tilde{\Cc'_0}$ the
special fibre of this morphism. Then $\tilde{\Cc'_0}$
has two irreducible components $R\isomorph\Pp^1_k$ 
and $\tilde{\Cc_0}\isomorph\Spec(k[[z]])$ which meet
in a $k$-rational ordinary double point of $\tilde{\Cc'_0}$.
Furthermore, $\tilde{\Cc}'_{\vers}$ over
the base $\Spec(k[[t]])$ is a versal deformation of
$\tilde{\Cc'_0}$.
\end{lemma}

\begin{proof}
It is clear that the special fibre $\tilde{\Cc}'_0$ of
$\tilde{\Cc}'_{\vers}\to\Spec(k[[t]])$ has the indicated shape.
To show that $\tilde{\Cc}'_{\vers}$ over $\Spec(k[[t]])$
is a versal deformation of its special fibre, it suffices
to prove that
$\Ext^2(\Omega^1_{\tilde{\Cc}'_0/k},\Oo_{\tilde{\Cc}'_0})=(0)$ and
$\Ext^1(\Omega^1_{\tilde{\Cc}'_0/k},\Oo_{\tilde{\Cc}'_0})\isomorph k$
and that the deformation $\tilde{\Cc}'_{\vers}$ has nonvanishing
Kodaira-Spencer class.
We identify $\tilde{\Cc}'_0$ with the closed subscheme of 
$\Pp_{k[[z]]}^1=\Proj(k[[z]][u,v])$ defined by the
equation $z.u$. Let $I\subset\Oo_{\Pp_{k[[z]]}^1}$ be the
defining sheaf of ideals. Then we have an exact sequence
$$
\xymatrix{
0\ar[r] & 
I/I^2 \ar[r] & 
\Omega^1_{\Pp^1_{k[[z]]}/k}|_{\tilde{\Cc}'_0} \ar[r] &
\Omega^1_{\tilde{\Cc}'_0/k}\ar[r] & 
0\quad.
}
$$
which is a resolution of $\Omega^1_{\tilde{\Cc}'_0/k}$ by locally
free $\Oo_{\tilde{\Cc}'_0}$-modules.
Applying the functor $\Hom(\cdot,\Oo_{\tilde{\Cc}'_0})$ to
this sequence, we obtain a long exact sequence of $\Ext$-groups.
It is easy to see that for 
$\F=\Omega^1_{\Pp^1_{k[[z]]}/k}|_{\tilde{\Cc}'_0}$
or
$\F=I/I^2$
we have 
$
\Ext^i(\F,\Oo_{\tilde{\Cc}'_0})=
H^i(\tilde{\Cc}'_0,\F^{\dual})=(0)
$ 
for every $i\geq 1$.
Therefore $\Ext^2(\Omega^1_{\tilde{\Cc}'_0/k},\Oo_{\tilde{\Cc}'_0})$ vanishes
and $\Ext^1(\Omega^1_{\tilde{\Cc}'_0/k},\Oo_{\tilde{\Cc}'_0})$ is the 
cokernel of the morphism
$
H^0(\tilde{\Cc}'_0,(\Omega^1_{\Pp^1_{k[[z]]}/k}|_{\tilde{\Cc}'_0})^{\dual})
\to
H^0(\tilde{\Cc}'_0,(I/I^2)^{\dual})
$
which in turn is easily seen to be one-dimensional.
By inspection, one also shows that the Kodaira-Spencer class of
the deformation $\tilde{\Cc}'_{\vers}$ is nonzero.
\end{proof}

\begin{lemma}
\label{adaptation}
Let $k$ be a field, $A$ a local artinian $k$-algebra with residue
field $k$ and let $S:=\Spec(A)$. Let $\tilde{\Cc}:=\Spec(A[[z]])$
and let $\tilde{s}:S\to\tilde{\Cc}$ be the section $z=0$ of the structure
morphism $\tilde{\pi}:\tilde{\Cc}\to S$. Let 
$(\tilde{\Cc}',\tilde{f},\tilde{\pi}',\tilde{s}')$ be an arbitrary simple 
modification of $(\tilde{\Cc}, \tilde{\pi}, \tilde{s})$ such
that $\tilde{f}$ is not an isomorphism.
Then there is a morphism $S\to\Spec(k[[t]])$ 
such that the diagram
$$
\xymatrix{
\text{$\tilde{\Cc}'$} \ar[rr]^{\tilde{f}}\ar[dr]^{\tilde{\pi}'}
& & 
\text{$\tilde{\Cc}$}\ar[dl]_{\tilde{\pi}} \\
& S \ar@/^/[ul]^{\tilde{s}'} \ar@/_/[ur]_{\tilde{s}} &
}
$$
is isomorphic to the diagram which is induced by the base change 
$S\to\Spec(k[[t]])$ from the diagram
$$
\xymatrix{
\text{$\tilde{\Cc}'_{\vers}$} 
\ar[rr]^{\tilde{f}_{\vers}} \ar[dr]
& & \text{$\tilde{\Cc}_{\vers}$} \ar[dl] \\
& \Spec(k[[t]])
\ar@/^/[ul]^{\tilde{s}'_{\vers}}
\ar@/_/[ur]_{\tilde{s}_{\vers}} &
}
$$
where $\tilde{s}_{\vers}$ is the section $z=0$ of 
$\tilde{\Cc}_{\vers} \to \Spec(k[[t]])$
and $\tilde{s}'_{\vers}$ is the proper transform of $\tilde{s}_{\vers}$.
\end{lemma}

\begin{proof}
By lemma \ref{universal}, there is a morphism $S\to\Spec(k[[t]])$ such
that $\tilde{\Cc}'$ is isomorphic to 
$\tilde{\Cc}_{\vers}\times_{\Spec(k[[t]])}S$.
Thus there exists an open affine covering 
$\tilde{\Cc}'=U_1\cup U_2$ with
$U_1:=\Spec(B_A[u]/(uz-a))$ and $U_2:=\Spec(B_A[v]/(z-va))$ where
$B_A:=A[[z]]$ and $a\in\m_A$ is the image of $t$ by the
morphism $k[[t]]\to A$. The glueing of $U_1$ and $U_2$ is
given by 
$$
\xymatrix@R=1ex{
1/u & &
v \ar@{|->}[ll] \\
\Spec(B_A[u]/(uz-a))_{(u)} \ar[rr]^{\isomorph} \ar[dr]_{\isomorph} & &
\Spec(B_A[v]/(z-va))_{(v)} \ar[dl]^{\isomorph} \\
& U_1\cap U_2 &
}
$$
The section $\tilde{s}':S\to\tilde{\Cc}'$ factorizes through $U_2$ and
is given by an $A$-morphism $\Oo_{U_2}\to A$, $v\mapsto\alpha$.
We want to show that there exists an automorphism 
$\varphi'$ of $\tilde{\Cc}'$ such that 
the composition $\varphi'\comp\tilde{s}'$ is given by the $A$-morphism
$\Oo_{U_2}\to A$, $v\mapsto 0$.
By a standard argument we may assume $\alpha\m_A=(0)$.
Let $\varphi'$ be defined on the affine
pieces $U_i$ by the $B_A$-morphisms
$$
\vcenter{
\xymatrix@R=.8ex{
\Oo_{U_1}\ar[r]^{\isomorph} & 
\Oo_{U_1} \\
u \ar@{|->}[r] &
u(1+\alpha u)
}}
\qquad\quad
\text{and}
\qquad\quad
\vcenter{
\xymatrix@R=1ex{
\Oo_{U_2}\ar[r]^{\isomorph} & 
\Oo_{U_2} \\
v \ar@{|->}[r] &
v-\alpha
}}
$$
The automorphism $\varphi'$ has the required property.
Therefore we may assume that $\alpha=0$, 
i.e. that $\tilde{s}'$ is induced by base change from $\tilde{s}'_{\vers}$.
The morphism $\tilde{f}$ is given by $A$-morphisms
$$
\vcenter{
\xymatrix@R=.8ex{
B_A\ar[r] & 
\Oo_{U_1} \\
z \ar@{|->}[r] &
z+p
}}
\qquad\quad
\text{and}
\qquad\quad
\vcenter{
\xymatrix@R=1ex{
B_A\ar[r] & 
\Oo_{U_2} \\
z \ar@{|->}[r] &
z+q
}}
$$
where $p\in\m_A\Oo_{U_1}$ and $q\in\m_A\Oo_{U_2}$.
Again we may assume that $\m_Ap=(0)=\m_Aq$.
The constraint $\tilde{s}=\tilde{f}\comp\tilde{s}'$ and
compatibility with the glueing morphism over $U_1\cap U_2$
has the consequence that in fact $q=0$ and $p=zp_1$ for
some $p_1\in\m_AB_A$ with $\m_Ap_1=(0)$.
We define an automorphism $\varphi$ of $\tilde{\Cc}=\Spec(B_A)$
by $z\mapsto(1-p_1)z$. 
Then the composed morphism 
$\varphi\comp\tilde{f}:\tilde{\Cc}\to\tilde{\Cc}'$ 
is given
by
$$
\vcenter{
\xymatrix@R=.8ex{
B_A\ar[r] & 
\Oo_{U_1} \\
z \ar@{|->}[r] &
z
}}
\qquad\quad
\text{and}
\qquad\quad
\vcenter{
\xymatrix@R=1ex{
B_A\ar[r] & 
\Oo_{U_2} \\
z \ar@{|->}[r] &
z
}}
$$
i.e. it is the one induced by base change from the morphism 
$\tilde{f}_{\vers}$. Since furthermore
we have $\varphi\comp\tilde{s}=\tilde{s}$, which is also
the morphism induced by base change from $\tilde{s}_{\vers}$,
this proves the lemma.
\end{proof}

\begin{lemma}
\label{presentation of N}
With the notation of lemma \ref{adaptation}, let 
$N:=\tilde{f}_*(\Oo_{\tilde{\Cc}'}(\tilde{s}'))
    \tensor\Oo_{\tilde{\Cc}}(-\tilde{s})$
and denote by
$\tilde{\Cc}'_0$, $\tilde{s}'_0$ etc. the objects induced by
$\tilde{\Cc}'$, $\tilde{s}'$ etc. via base change $\Spec(k)\to S$.
Then there exist canonical isomorphisms
$$
N\tensor_{\Oo_{\tilde{\Cc}}} Q_A\isomto Q_A
\quad
\text{and}
\quad
N\tensor_A k\isomto 
\m_{\tilde{\Cc}_0,\tilde{s}_0} \oplus 
((\m_{\tilde{\Cc}'_0,\tilde{s}'_0}/ \m_{\tilde{\Cc}'_0,\tilde{s}'_0}^2)^{\dual}
\tensor
\m_{\tilde{\Cc}_0,\tilde{s}_0}/ \m_{\tilde{\Cc}_0,\tilde{s}_0}^2)
$$
where $Q_A:=\Oo_{\tilde{\Cc}\setminus\{\tilde{s}\}}=A[[z]][1/z]$,
such that the following diagram commutes
$$
\vcenter{
\xymatrix@R=2ex{
N \ar[r] \ar[d] &
N\tensor_Ak \ar[r]^(.2){\isomorph} &
\m_{\tilde{\Cc}_0,\tilde{s}_0} \oplus 
((\m_{\tilde{\Cc}'_0,\tilde{s}'_0}/ \m_{\tilde{\Cc}'_0,\tilde{s}'_0}^2)^{\dual}
\tensor
\m_{\tilde{\Cc}_0,\tilde{s}_0}/ \m_{\tilde{\Cc}_0,\tilde{s}_0}^2)
\ar[d]^{\text{pr}_1}
\\
N\tensor_{\Oo_{\tilde{\Cc}}} Q_A \ar[d]^{\isomorph} & &
\m_{\tilde{\Cc}_0,\tilde{s}_0} \ar@{^{(}->}[d]
\\
Q_A \ar[r] &
Q_A\tensor_Ak \ar@{=}[r] &
\Oo_{\tilde{\Cc}_0\setminus\{\tilde{s}_0\}}
}}
\eqno(*)
$$
Furthermore, there is an element $a\in\m_A$ and an
exact sequence of $\Oo_{\tilde{\Cc}}$-modules
$$
\xymatrix@R=0.7ex{
0\ar[r] &
\Oo_{\tilde{\Cc}}\ar[r] &
\Oo_{\tilde{\Cc}}\oplus\Oo_{\tilde{\Cc}}\ar[r] & 
N\ar[r] &
0 \\
& 1\ar@{|->}[r] & (a,-z)
}
$$
which by tensoring with $Q_A$ and $k$ induces the exact sequences
$$
\xymatrix@R=0.7ex{
0\ar[r] &
Q_A\ar[r] &
Q_A\oplus Q_A\ar[r] & 
Q_A\ar[r] &
0 \\
& & (1,0)\ar@{|->}[r] & z &\\
& & (0,1)\ar@{|->}[r] & a &\\
}
$$
and
$$
\xymatrix@R=0.7ex{
0\ar[r] &
\Oo_{\tilde{\Cc}_0}\ar[r] &
\Oo_{\tilde{\Cc}_0}\oplus\Oo_{\tilde{\Cc}_0}\ar[r] & 
\m_{\tilde{\Cc}_0,\tilde{s}_0} \oplus 
((\m_{\tilde{\Cc}'_0,\tilde{s}'_0}/ \m_{\tilde{\Cc}'_0,\tilde{s}'_0}^2)^{\dual}
\tensor
\m_{\tilde{\Cc}_0,\tilde{s}_0}/ \m_{\tilde{\Cc}_0,\tilde{s}_0}^2)
\ar[r] &
0 \\
& & (1,0)\ar@{|->}[r] & (z,0) &\\
& & (0,1)\ar@{|->}[r] & (0,e) &\\
}
$$
where $0\neq e\in 
(\m_{\tilde{\Cc}'_0,\tilde{s}'_0}/ \m_{\tilde{\Cc}'_0,\tilde{s}'_0}^2)^{\dual}
\tensor
\m_{\tilde{\Cc}_0,\tilde{s}_0}/ \m_{\tilde{\Cc}_0,\tilde{s}_0}^2$.
\end{lemma}

\begin{proof}
Since both $\tilde{f}_*\Oo_{\tilde{\Cc}'}(\tilde{s}')$ and 
$\Oo_{\tilde{\Cc}}(-\tilde{s})$ are canonically trivialized 
over $\tilde{\Cc}\setminus\tilde{s}$, the isomorphism
$N\tensor_{\Oo_{\tilde{\Cc}}}Q_A\isomto Q_A$ is clear. 
By \cite{Knudsen}, Cor. 1.5 we have 
$N\tensor_A k=(f_0)_*(\Md_0)$. Let $R$ denote the component
isomorphic to $\Pp_k^1$ of $\tilde{\Cc}'_0$ and by abuse of
notation let $\tilde{\Cc}_0$ also denote the component of 
$\tilde{\Cc}'_0$, which is mapped isomorphically
to $\tilde{\Cc}_0$ by $\tilde{f}_0$.
Then we have, similarly as in the proof of \ref{S=Spec k},
\begin{eqnarray*}
(f_0)_*(\Md_0) &=&
\Md_0|_{\tilde{\Cc}_0}\times_{\Md_0[\tilde{s}_0]}H^0(R,\Md_0|_R) 
 \\
&=&
\Md_0|_{\tilde{\Cc}_0}\oplus
\ker(H^0(R,\Md_0|_R)\to \Md_0[\tilde{s}_0])
\\
&=&
\m_{\tilde{\Cc}_0,\tilde{s}_0} \oplus 
((\m_{\tilde{\Cc}'_0,\tilde{s}'_0}/ \m_{\tilde{\Cc}'_0,\tilde{s}'_0}^2)^{\dual}
\tensor
\m_{\tilde{\Cc}_0,\tilde{s}_0}/ \m_{\tilde{\Cc}_0,\tilde{s}_0}^2)
\end{eqnarray*}
The commutativity of diagram $(*)$ comes from the commutativity of
the diagram
$$
\xymatrix@R=2ex{
\text{$\tilde{\Cc}$} & \text{$\tilde{\Cc}_0$} \ar@{_{(}->}[l] \\
\text{$\tilde{\Cc}\setminus\tilde{s}$} \ar@{^{(}->}[u] &
\text{$\tilde{\Cc}_0\setminus\tilde{s}_0$} \ar@{^{(}->}[u] \ar@{_{(}->}[l]
}
$$ 
Let $A_{\vers}:=k[[t]]$ and let $S\to\Spec(A_{\vers})$ be a morphism
such that 
$\tilde{\Cc}'_0\isomorph \tilde{\Cc}'_{\vers}\times_{\Spec(A_{\vers})}S$.
Let 
$
\Md_{\vers}:=\Oo_{\tilde{\Cc}'_{\vers}}(\tilde{s}'_{\vers}))\tensor
             (f_{\vers})^*(\Oo_{\tilde{\Cc}_{\vers}}(-\tilde{s}_{\vers}))
$
and
$
N_{\vers}:=(f_{\vers})_*\Md_{\vers}
\quad.
$
By \cite{Knudsen}, Cor. 1.5, we have $N=N_{\vers}\tensor_{A_{\vers}}A$.
For the existence of the exact sequence
$
0 \to \Oo_{\tilde{\Cc}} \to 
\Oo_{\tilde{\Cc}}\oplus\Oo_{\tilde{\Cc}} \to
N \to 0
$
it suffices therefore to produce an exact sequence
$$
\vcenter{
\xymatrix@R=0.7ex{
0\ar[r] &
B_{A_{\vers}}\ar[r] &
B_{A_{\vers}}\oplus B_{A_{\vers}}\ar[r] & 
N_{\vers}\ar[r] &
0 \\
& 1\ar@{|->}[r] & (t,-z)
}}
\eqno(\dagger)
$$
where $B_{A_{\vers}}:=A_{\vers}[[z]]=k[[t]][[z]]$.
Observe that we have $\Md_{\vers}=\Oo_{\tilde{\Cc}'_{\vers}}(-R)$, where
$R$ is the exceptional divisor of the blowing-up morphism $f_{\vers}$.
We have an open affine covering $\tilde{\Cc}'_{\vers}=U_1\cup U_2$, where
$U_1=\Spec B_{A_{\vers}}[t/z]$,
$U_2=\Spec B_{A_{\vers}}[z/t]$ and
$U_1\cap U_2=\Spec B_{A_{\vers}}[t/z,z/t]$.
Over $U_1$ and $U_2$, the divisor $R$ is given by the equation
$z$ and $t$ respectively. Therefore we have
$$
N_{\vers}=\{(g,h)\in B_{A_{\vers}}[t/z]\times B_{A_{\vers}}[z/t]\ |\ \ 
             \text{$z\cdot g=t\cdot h$ in $B_{A_{\vers}}[t/z,z/t]$}\}
$$
It is easy to see that the morphism 
$B_{A_{\vers}}\oplus B_{A_{\vers}}\to N_{\vers}$
given by $(1,0)\mapsto (1,z/t)$ and $(0,1)\mapsto (t/z,1)$ 
is surjective whith kernel $B_{A_{\vers}}\cdot(t,-z)$.
The compatibility (in the sense of the lemma) 
of the exact sequence $(\dagger)$ with the
canonical isomorphisms 
$N_{\vers}|_{(\tilde{\Cc}_{\vers}\setminus s_{\vers)}}\isomto B_{A_{\vers}}[1/z]$
and 
$
N_{\vers}|_{\tilde{\Cc}_0} \isomto 
\m_{\tilde{\Cc}_0,\tilde{s}_0} \oplus 
((\m_{\tilde{\Cc}'_0,\tilde{s}'_0}/ \m_{\tilde{\Cc}'_0,\tilde{s}'_0}^2)^{\dual}
\tensor
\m_{\tilde{\Cc}_0,\tilde{s}_0}/ \m_{\tilde{\Cc}_0,\tilde{s}_0}^2)
$ 
follows by inspection.
\end{proof}

\begin{lemma}
\label{surjective}
In the situation of \ref{piMM to M} let 
$S=\Spec(A)$ for a local artinian $k$-algebra $A$ with residue
field $k$ and let $\Cc'_0:=\Cc'\tensor_A k$. Then the restricion map
$$
H^0(\Cc',\Md)\to H^0(\Cc'_0,\Md|_{\Cc'_0})
$$
is surjective.
\end{lemma}

\begin{proof}
Let $\Nc$ be the $\Oo_{\Cc}$-module $f_*(\Md)$.
By \cite{Knudsen}, cor 1.5, the formation of $\Nc$ commutes with
base change $S'\to S$. Therefore it suffices to show the surjectivity
of the restriction morphism
$$
H^0(\Cc,\Nc)\to H^0(\Cc_0,\Nc_0) \quad,
$$
where $\Cc_0:=\Cc\tensor_A k$ and $\Nc_0:=\Nc|_{\Cc_0}$.
Let $\tilde{\Cc}$ be the completion of $\Cc$ along the section $s$.
Since $\pi$ is smooth in the neighbourhood of $s$, we may identify
$\tilde{\Cc}$ with $\Spec(A[[z]])$ such that $z$ corresponds to a generator
of the $\Oo_{\Cc}$-Ideal defining $s$.
Let $\tilde{\Cc}':=\Cc'\times_{\Cc}\tilde{\Cc}$ and denote by
$\tilde{f}$ the projection $\tilde{\Cc}'\to\tilde{\Cc}$ onto
the second factor. Then we have a canonical isomorphism
$N:=\tilde{f}_*(\Oo_{\tilde{\Cc}'}(\tilde{s}'))
    \tensor\Oo_{\tilde{\Cc}}(-\tilde{s})
\isomto
\Nc\tensor_{\Oo_{\Cc}}\Oo_{\tilde{\Cc}}\quad.
$
Let $Q_A$ be the localization of $\Oo_{\tilde{\Cc}}$ with respect to $z$.
Since $\Nc$ is naturally trivialized on the complement of $s$, we have
canonical isomorphisms $\Nc\tensor_{\Oo_{\Cc}}Q_A\isomto Q_A$ and
$N\tensor_{\Oo_{\tilde{\Cc}}}Q_A\isomto Q_A$. By \cite{Beauville-Laszlo}
the canonical morphism
$$
H^0(\Cc,\Nc)\to N\times_{Q_A}H^0(\Cc\setminus\{s\},\Nc)
$$
is an isomorphism. Likewise, we have a canonical isomorphism
$$
H^0(\Cc_0,\Nc_0)\to N_0\times_{Q_k}H^0(\Cc_0\setminus\{s\},\Nc_0)
$$
where $N_0:=N\tensor_A k$ and $Q_k:=Q_A\tensor_Ak$.
By the first part of lemma \ref{presentation of N}
we have
\begin{eqnarray*}
N_0\times_{Q_k}H^0(\Cc_0\setminus\{s\},\Nc_0) 
&=&
H^0(\Cc_0,\Oo_{\Cc_0}(-s_0))\oplus
((\m_{\tilde{\Cc}'_0,\tilde{s}'_0}/ \m_{\tilde{\Cc}'_0,\tilde{s}'_0}^2)^{\dual}
\tensor
\m_{\tilde{\Cc}_0,\tilde{s}_0}/ \m_{\tilde{\Cc}_0,\tilde{s}_0}^2)
\\
&=&
(\m_{\tilde{\Cc}'_0,\tilde{s}'_0}/ \m_{\tilde{\Cc}'_0,\tilde{s}'_0}^2)^{\dual}
\tensor
\m_{\tilde{\Cc}_0,\tilde{s}_0}/ \m_{\tilde{\Cc}_0,\tilde{s}_0}^2
\quad.
\end{eqnarray*}
Let $\nu\in N$ be the image of 
$(0,1)\in \Oo_{\tilde{\Cc}}\oplus\Oo_{\tilde{\Cc}}$
by the surjection
$\Oo_{\tilde{\Cc}}\oplus\Oo_{\tilde{\Cc}}\to N$
from lemma \ref{presentation of N}.
By the second part of \ref{presentation of N} the pair
$(\nu,a)$ is an element of 
$N\times_{Q_A}H^0(\Cc\setminus\{s\},\Nc)$ and its image
under the restriction morphism 
$$
N\times_{Q_A}H^0(\Cc\setminus\{s\},\Nc)
\to
N_0\times_{Q_k}H^0(\Cc_0\setminus\{s\},\Nc_0)=
(\m_{\tilde{\Cc}'_0,\tilde{s}'_0}/ \m_{\tilde{\Cc}'_0,\tilde{s}'_0}^2)^{\dual}
\tensor
\m_{\tilde{\Cc}_0,\tilde{s}_0}/ \m_{\tilde{\Cc}_0,\tilde{s}_0}^2
$$
is $e\neq 0$.
\end{proof}

\begin{lemma}
\label{openness of simple modification}
Let $S$ be a locally noetherian scheme, $\pi:\Cc\to S$ a finitely presented
flat morphism whose geometric fibres are connected of dimension one.
Let $f:\Cc'\to\Cc$ be a proper morphism such that $\pi':=\pi\comp f$
is flat. Let $s':S\to\Cc'$ be a section of $\pi'$ and let $s:=f\comp s'$.
Assume that the restriction of $f$ to $\Cc'\setminus s'(S)$
is an isomorphism onto $\Cc\setminus s(S)$. 
Let $z\in S$ be a point. Denote by $\Cc_z$ and $\Cc'_z$ the fibres
of $\pi$ and of $\pi'$ over $z$ and by $\pi_z$, $\pi'_z$, $f_z$, $s_z$,
$s'_z$ the the restriction of $\pi$, $\pi'$, $f$, $s$, $s'$ on
the fibres. Assume that $(\Cc'_z,f_z,\pi'_z,s'_z)$ is a 
simple modification
of $(\Cc_z, \pi_z, s_z)$. Then there is an open neighbourhood $U\subseteq S$
of $z$, such that $(\Cc',f,\pi',s')|_U$ is a simple modification of 
$(\Cc,\pi,s)|_U$.
If we assume in addition that $H^0(\Cc_z,\Oo_{\Cc_z})=\kappa(z)$ and
$H^0(\Cc_z,\Oo_{\Cc_z}(-s(z)))=(0)$, then we
can choose $U$ such that the analogous property
holds for all $z'\in U$.
\end{lemma}

\begin{proof}
Admitting the first part of the lemma,
the last assertion in the lemma follows immediately by
\cite{EGA III}(7.8.7),(7.8.8) and (7.7.5)I.
By assumption, the $s(z)$ and $s'(z)$ are smooth points in their fibre.
Therefore there exists an open neighbourhood $V\subseteq\Cc$ of $s(z)$ and 
$V'\subseteq\Cc'$ of $s'(z)$ such that the restriction of $\pi$ to $V$
and the restriction of $\pi'$ to $V'$ is smooth (\cite{EGA IV}(12.1.7)(iv)). 
Replacing $S$ by $s^{-1}(V)\cap (s')^{-1}(V')$, we may assume 
that $s$ and $s'$ meet $\Cc$ and $\Cc'$ in the smooth locus of 
$\pi$ and $\pi'$.
If $f_z$ is an isomorphism, then by \cite{EGA IV}(17.8.3) and
\cite{EGA III}(4.4.5) there is a neighbourhood $W\subseteq\Cc$ of 
$s(z)$ such that the restriction of $f$ to $W':=f^{-1}(W)$ is
a finite etale morphism $W'\to W$. Since its degree over $s(z)$
is one, it is in fact an isomorphism. It is clear that in
this case $U:=s^{-1}(W)$ has the required property.
Now assume that $f_z$ is not an isomorphism.
Let $\widehat{S}$ be the completion of $S$ at $z$ and let
$\widehat{\Cc}$ (resp. $\widehat{\Cc}'$)
be the completion of 
$\Cc\times_S\widehat{S}$ (resp. of $\Cc'\times_S\widehat{S}$) 
along the subscheme 
$s(S)\times_S\widehat{S}$ (resp. $s'(S)\times_S\widehat{S}$).
Denote by 
$\widehat{\pi}$, $\widehat{\pi}'$, $\widehat{f}$, $\widehat{s}$, 
$\widehat{s}'$
the morphisms between $\widehat{S}$, $\widehat{\Cc}$ and $\widehat{\Cc}'$
induced by $\pi$, $\pi'$, $f$, $s$, $s'$.
From lemma \ref{universal} it follows that 
$(\widehat{\Cc'},\widehat{f},\widehat{\pi'},\widehat{s}')$ is a 
simple modification of 
$(\widehat{\Cc},\widehat{\pi},\widehat{s})$.
Since $\widehat{S}\to\Spec(\Oo_{S,z})$ is surjective 
(cf. \cite{Matsumura}\S 8), it follows
that $(\Cc',f,\pi',s')|_{\Spec(\Oo_{S,z})}$ is a simple modification of 
$(\Cc,\pi,s)|_{\Spec(\Oo_{S,z})}$. Thus it suffices to show
that the set $T$ of all $z'\in S$ such that 
$(\Cc'_{z'},f_{z'},\pi'_{z'},s'_{z'})$ is a 
simple modification
of $(\Cc_{z'}, \pi_{z'}, s_{z'})$ is a locally constructible 
subset of $S$. 
But this is a consequence of the following characterization
of $T$:
A point $z'\in S$ is contained in $T$ if and only if
(i) $f^{-1}(s(z'))$ is irreducible, smooth and of dimension $\leq 1$.
(ii) $H^1(f^{-1}(s(z')),\Oo)=0$.
(iii) $(\pi')^{-1}(z')\cap f^{-1}(V)$ is connected and is either smooth
      or contains precisely one ordinary double point.
\end{proof}

\begin{lemma}
\label{reduction to noetherian}
Let $S$ be an affine scheme, let $\pi:\Cc\to S$ be a proper
finitely presented
flat morphism whose geometric fibres are connected of dimension one
and let $s:S\to \Cc$ be a section of $\pi$ whose image consists of points
where $\pi$ is smooth. Let $(\Cc',f,\pi',s')$ be a simple modification of 
$(\Cc,\pi,s)$. Then there is a noetherian scheme $S_0$, a proper flat morphism
$\pi_0:\Cc_0\to S_0$ with connected geometric fibres of dimension one,
a section $s_0$ of $\pi_0$ meeting $\Cc_0$ in the smooth locus of $\pi_0$,
a simple modification 
$(\Cc'_0,f_0,\pi'_0,s'_0)$ of $(\Cc_0,\pi_0,s_0)$ and a morphism $S\to S_0$
such that the data $(\Cc',f,\pi',s',\Cc,\pi,s)$ are induced by pull back from
the data $(\Cc'_0,f_0,\pi'_0,s'_0,\Cc_0,\pi_0,s_0)$.
If we assume in addition that $H^0(\Cc_z,\Oo_{\Cc_z})=\kappa(z)$ and
$H^0(\Cc_z,\Oo_{\Cc_z}(-s(z)))=(0)$ for all $z\in S$, then we
can choose $(\Cc_0,\pi_0,s_0)$ such that the analogous property
holds for all $z_0\in S_0$.
\end{lemma}

\begin{proof}
Let $A$ be the coordinate ring of $S$ and let $(A_{\lambda})_\lambda$
be the inductive system of all noetherian subrings of $A$.
Let $S_\lambda:=\Spec(A_\lambda)$ for all $\lambda$.
Using the results from \cite{EGA IV} cited in the proof of
\ref{noetherian} one shows that there exists an index
$\lambda$, a proper flat morphism $\pi_\lambda:\Cc_\lambda\to S_\lambda$,
a proper morphism $f_\lambda:\Cc'_\lambda\to\Cc_\lambda$ such that 
the composition $\pi'_\lambda:=\pi_\lambda\comp f\lambda$ is flat,
and a morphism $s'_\lambda: S_\lambda\to\Cc'_\lambda$ such that
the data 
$(\Cc',f,\pi',s',\Cc,\pi,s)$ 
are induced from the data
$(\Cc'_\lambda,f_\lambda,\pi'_\lambda,s'_\lambda,
  \Cc_\lambda,\pi_\lambda,s_\lambda)$ 
where $s_\lambda:=f\comp s'_\lambda$.
By \cite{EGA IV}(8.8.2)(i) and (8.8.2.4) 
we may assume that $s'_\lambda$ is
a section of $\pi'_\lambda$ and that the restriction of
$f_\lambda$ to the open subset 
$f_\lambda^{-1}(\Cc_\lambda\setminus s_\lambda(S_\lambda))$
of $\Cc'_\lambda$ is an isomorphism onto 
$\Cc_\lambda\setminus s_\lambda(S_\lambda)$.
To conclude,
we now make use of \ref{openness of simple modification}
similarly as we employed \ref{open} in the proof of lemma 
\ref{noetherian}.
\end{proof}

\begin{proof}(of lemma \ref{piMM to M})
Consider first the case where $S$ is locally noetherian.
Then both $\pi'_*\Md$ and $(s')^*\Md$ are coherent
$S$-modules and $(s')^*\Md$ is locally free, so
lemma \ref{piMM to M} follows easily from lemma \ref{S=Spec k}, if
we can show that the formation of $\pi'_*\Md$ commutes
with base change by a morphism from any scheme to $S$.
But this base change property follows from lemma \ref{surjective}
and \cite{EGA III} 7.7.5 and 7.7.10.
In the general case we may assume $S$ to be affine.
By lemma \ref{reduction to noetherian} there exists a morphism
from $S$ to a locally noetherian scheme $S_0$ such that
the given data $(\Cc',f,\pi',s',\Cc,\pi,s)$ are induced by pull 
back from data $(\Cc'_0,f_0,\pi'_0,s'_0,\Cc_0,\pi_0,s_0)$ defined
over $S_0$. Let 
$\Md_0:=\Oo_{\Cc'_0}(s'_0)\tensor_{\Oo_{\Cc'_0}}
        f_0^*\Oo_{\Cc_0}(-s_0)$.
By the first case, we have $(\pi'_0)_*\Md_0\isomto (s'_0)^*\Md_0$.
Since furthermore our morphism $\pi'_*\Md\to (s')^*\Md$ is induced
from this morphism by base change, we are done.
\end{proof}

\section{Admissible bundles on simple modifications and bf-morphisms}
\label{section admissible bundles}

Throughout this section, $S$ will be a locally Noetherean scheme,
$\pi:\Cc\to S$ a proper and flat morphism  with connected geometric
fibers of dimension one such that $\pi_*\Oo_{\Cc}=\Oo_S$, 
and $s:S\to\Cc$ will be a section of $\pi$ whose
image consists of smooth points of $\pi$.
Furthermore we fix an $\Oo_{\Cc}$-module $\E$ and an $\Oo_S$-module
$E'$, both locally free of rank $n$. Let $E:=s^*\E$.

\begin{definition}
\label{admissible of degree d}
Let $(\Cc',f,\pi',s')$ be a simple modification of $(\Cc,\pi,s)$ 
(cf. \ref{def mpc}) and let $0\leq d\leq n$. 
A locally free $\Oo_{\Cc'}$-module $\E'$
of rank $n$ will be called {\em admissible of degree $d$ 
for $(\Cc',f,\pi',s')$}, if
for any point $x\in\Cc$ with $f^{-1}(x)\isomorph\Pp^1_{\kappa(x)}$
we have
$$
\E'|_{f^{-1}(x)} \ \isomorph\
\Oplus^d
\Oo_{\Pp^1_{\kappa(x)}}(1)
\ \oplus\
\Oplus^{n-d}
\Oo_{\Pp^1_{\kappa(x)}}
\quad.
$$
\end{definition}

\begin{construction}
\label{constr3}
Let $(\Cc',f,\pi',s')$ be a simple modification of $(\Cc,\pi,s)$ and let
$\E'$ be an admissible $\Oo_{\Cc'}$-module of degree $d$ together 
with isomorphisms 
$$
(f_*\E'(-s'))(s) \isomto \E 
\qquad \text{and} \qquad
(s')^*\E' \isomto E'
\quad.
$$

Since $R^1f_*\E'(-s')=0$ by \cite{Knudsen} 1.5,
we have the commutative diagram with exact rows:
$$
\xymatrix{
0 \ar[r] &
f_*\E'(-s') \ar[d] \ar[r] &
f_*\E' \ar[d] \ar[r] \ar@{.>}[dl]&
s_*(s')^*\E' \ar[d]^0 \ar[r] &
0 
\\
0 \ar[r] &
(f_*\E'(-s'))(s) \ar[r] &
(f_*\E')(s) \ar[r] &
(s_*(s')^*\E')(s) \ar[r] &
0
}
$$
The vanishing of the right vertical arrow implies that
the middle vertical arrow factorizes as indicated by the dotted arrow.
Composing this arrow with the isomorphism
$(f_*\E'(-s'))(s)\isomto \E$,
we obtain a morphism $f_*\E'\to \E$.
Now consider the diagram with exact rows
$$
\xymatrix{
0 \ar[r] &
f_*\E'(-s') \ar[d]^{\isomorph} \ar[r] &
f_*\E' \ar[d] \ar[r] &
s_*(s')^*\E' \ar@{.>}[d] \ar[r] &
0 
\\
0 \ar[r] &
\E(-s) \ar[r] &
\E \ar[r] &
s_*s^*\E \ar[r] &
0
}
$$
By the commutativity of the left square we obtain
a morphism as indicated by the dotted arrow.
This induces a morphism $E'\to E$.

Let $\M:=\Oo_{\Cc'}(-s')\tensor f^*\Oo_{\Cc}(s)$.
By the projection formula we have 
$f_*(\M\tensor\E')=(f_*\E'(-s'))(s)$.
Thus we obtain a morphism
$$
E=s^*\E\isomto s^*f_*(\M\tensor\E')=(s')^*f^*f_*(\M\tensor\E')
\to (s')^*(\M\tensor\E')\isomto M\tensor E'
\quad,
$$
where  $(M,\mu)$ is the invertible $\Oo_S$-module together with a global
section associated to the simple modification $(\Cc',f,\pi',s')$
(cf. proposition \ref{mod to Mmu}).
\end{construction}

\begin{remark}
\label{constr3 over K}
In case $S=\Spec(K)$ for a field $K$, 
the curve $\Cc$ can be identified with
a closed subscheme of $\Cc'$ and construction \ref{constr3}
can be specified more concretely as follows.
Assume that $f$ is not an isomorphism and let 
$R:=f^{-1}(s)\isomorph \Pp_K^1$.
Then $M=(\m/\m^2)^{\dual}\tensor(\m'/(\m')^2)$ and $\mu=0$, where
$\m$ (resp. $\m'$) is the maximal ideal belonging to the point $s\in\Cc$
(resp. to the point $s'\in\Cc'$).
Observe that we have 
$$
\Oo_{\Cc,s}\tensor f_*\E'(-s')=
(\E'|_{\Cc})_s \times_{\E'[s]} H^0(R,\E'(-s')|_R)=
(\E'|_{\Cc})_s \times_{\E'[s]} \F[s]
$$
where $0\to\F\to\E'|_R\to\G\to 0$ is the canonical exact sequence
with $\F\isomorph\oplus^d\Oo_R(1)$ and $\G\isomorph \oplus^{n-d}\Oo_R$.
Therefore we have a canonical exact sequence of $\Oo_{\Cc}$-modules:
$$
0 \to \E(-s) \to \E'|_{\Cc} \to \G[s] \to 0
\quad.
$$
Tensoring this exact sequence with $K=\kappa(s)$ leads to the exact
sequence 
$$
0\to(\m/\m^2)\tensor\G[s]\to
(\m/\m^2)\tensor E \to \E'[s] \to \G[s] \to 0
$$
and the image of $(\m/\m^2)\tensor E\to\E'[s]$ is the vector space 
$\F[s]$ which is canonically isomorphic to 
$(\m'/(\m')^2)\tensor\F[s']$.
Therefore we arrive at the two exact sequences
$$
\xymatrix@R=2ex{
0 \ar[r] &
\F[s'] \ar[r] &
E' \ar[r] &
\G[s']=\G[s] \ar[r] &
0
\\
0\ar[r] & 
\G[s] \ar[r] &
E \ar[r] &
M\tensor\F[s'] \ar[r] &
0
}
$$
The morphisms $E'\to E$ and $E\to M\tensor E'$ from
construction \ref{constr3} are just
the compositions $E'\to\G[s]\to E$ and 
$E\to M\tensor\F[s']\to M\tensor E'$.
\end{remark}

\begin{proposition}
\label{mod&E' to bf}
The construction \ref{constr3} defines an isomorphism of groupoids:
$$
\left\{
\begin{array}{lllllll}
\text{simple modifications}\\ 
\text{$(\Cc',f,\pi',s')$ of $(\Cc,\pi,s)$ } \\ 
\text{together with an}\\
\text{admissible $\Oo_{\Cc'}$-module $\E'$ of} \\
\text{degree $d$ and isomorphisms}\\
\text{$(f_*\E'(-s'))(s) \isomto \E$} \\
\text{and $(s')^*\E' \isomto E'$}
\end{array}
\right\}
\isomto
\left\{
\begin{array}{llll}
\text{bf-morphisms} \\
(M,\ \mu,\ E'\to E,\ M\tensor E' \ot E,\ n-d) \\
\text{of rank $n-d$ from $E'$ to $E$}\\
\text{in the sense of \cite{Kausz} 3.1}
\end{array}
\right\}
$$
\end{proposition}

In the  remaining of this section we will prove proposition
\ref{mod&E' to bf}.
Let $M$ be an invertible $\Oo_S$-module and $\mu$ a global section 
of $M$. Let $(\Cc',f,\pi',s')$ be the corresponding simple modification
of $(\Cc,\pi,s)$ (cf. proposition \ref{mod to Mmu}).
Let $\E_1$ and $\E_2$ be two locally free $\Oo_{\Cc}$-modules of
rank $d$ and $n-d$ respectively. To these data we associate 
a locally free $\Oo_{\Cc'}$-module $\E'(\E_1,\E_2)$ of rank $n$ as follows:
$$
\E'(\E_1,\E_2):=
\left(
\Md\tensor 
f^*\E_1
\right)
\ \oplus\  f^*\E_2
\quad,
$$
where $\M:=\Oo_{\Cc'}(-s')\tensor f^*\Oo_{\Cc}(s)$.

\begin{lemma}       
\label{special case}
Let $\E':=\E'(\E_1,\E_2)$ and $E_i:=s^*\E_i$ for $i=1,2$. Then 
$\E'$ is an admissible $\Oo_{\Cc'}$-module of degree $d$ and
there are natural isomorphisms
$$
(f_*\E'(-s'))(s)\isomto\E_1\oplus\E_2 
\quad \text{and} \quad 
(s')^*\E'\isomto (M^{\dual}\tensor E_1)\oplus E_2
\quad .
$$
Furthermore, the following diagram commutes:
$$
\xymatrix@R=2ex{
(M^{\dual}\tensor E_1)\oplus E_2 \ar[dd]_{\id} \ar[rr] & &
E_1\oplus E_2  \ar[dd]_{\id} \ar[rr] & &
E_1\oplus(M\tensor E_2) \ar[dd]_{\id}  & \\
& & & & &  \\
(M^{\dual}\tensor E_1) \oplus E_2 
\ar[rr]^(.6){
\left[
\begin{array}{ll}
\mu & 0 \\
0 & 1 
\end{array}
\right]
}
& &
E_1\oplus E_2
\ar[rr]^(.4){
\left[
\begin{array}{ll}
1 & 0 \\
0 & \mu 
\end{array}
\right]
}
& &
E_1\oplus(M\tensor E_2)
&
}
$$
where the upper horizontal arrows are defined as in 
construction \ref{constr3}.
\end{lemma}

\begin{proof}
It is clear that $\E'$ is admissible.
By lemma \ref{f_*} we have canonical isomorphisms
$f_*\Oo_{\Cc'}\isomto\Oo_{\Cc}$ and $f_*\M\isomto\Oo_{\Cc}$.
Therefore we get an isomorphism
$$
(f_*\E'(-s'))(s)=f_*(\M\tensor\E')=
(f_*\Oo_{\Cc'}\tensor\E_1)\oplus (f_*\M\tensor\E_2)
\isomto \E_1\oplus\E_2
\quad.
$$
The isomorphism $(s')^*\E'\isomto (M^{\dual}\tensor E_1)\oplus E_2$
is the obvious one.
For the commutativity of the diagram it 
suffices to consider the two cases $(\E_1,\E_2)=(\Oo_{\Cc},0)$
and $(\E_1,\E_2)=(0,\Oo_{\Cc})$.
In the first case $\E'=\Md$.
The main point in this case is the observation that the diagram
$$
\xymatrix@R=3ex{
0 \ar[r] &
\Oo_{\Cc}(-s) 
\ar[r] \ar[d] &
f_*\Md
\ar[r] \ar@{.>}[dl] \ar[d] &
s_*(s')^*\M 
\ar[r] \ar[d]^0 &
0 \\
0 \ar[r]&
\Oo_{\Cc} 
\ar[r] &
f_*\Oo_{\Cc'}(s') 
\ar[r] &
s_*(s')^*\Oo_{\Cc'}(s') 
\ar[r] & 
0 
}
$$
defining 
$f_*\Md=f_*\E'\to (f_*\E'(-s'))(s)=\Oo_{\Cc}$ 
(cf. construction \ref{constr3})
is obtained by applying $f_*$ to the diagram
$$
\xymatrix@R=3ex{
0 \ar[r] &
f^*\Oo_{\Cc}(-s) 
\ar[r] \ar[d] &
\Md
\ar[r] \ar@{.>}[dl] \ar[d] &
s'_*(s')^*\M 
\ar[r] \ar[d]^0 &
0 \\
0 \ar[r]&
\Oo_{\Cc'} 
\ar[r] &
\Oo_{\Cc'}(s') 
\ar[r] &
s'_*(s')^*\Oo_{\Cc'}(s') 
\ar[r] & 
0 
}
$$
and that the dotted arrow of this last diagram is dual to the dotted
arrow in construction \ref{constr1}.
If $(\E_1,\E_2)=(0,\Oo_{\Cc})$, then $\E'=\Oo_{\Cc'}$. Here the main
point is to show that the composition of the isomorphism
$\Oo_{\Cc'}\isomto f^*f_*\M$
with the adjunction morphism
$f^*f_*\M\to\M$
coincides with the dotted arrow in construction
\ref{constr1}. But this follows from the commutativity of the
following diagram
$$
\xymatrix{
f^*f_*\M \ar[r] \ar[d]^{\isomorph} &
\M \ar[d] \\
f^*\Oo_{\Cc} \ar[r] &
f^*\Oo_{\Cc}(s)
}
$$
which in turn is easy to see.
\end{proof}

\begin{lemma}       
\label{general case}
Let $\E'$ be an arbitrary admissible $\Oo_{\Cc'}$-module of degree $d$.
Then locally in a neighbourhood of $s(S)\subset\Cc$ there exist free
$\Oo_{\Cc}$-modules $\E_1$ and $\E_2$ of rank $d$ and $n-d$ respectively,
such that 
$\E'$ is isomorphic to $\E'(\E_1,\E_2)$.
\end{lemma}

\begin{proof}
We may assume that $S=\Spec(A)$, $\Cc=\Spec(B)$ for a noetherian ring $A$
and smooth $A$-algebra $B$. By proposition \ref{mod to Mmu} we may
assume furthermore that $(\Cc',f,\pi',s')$ is the simple modification of
$(\Cc,\pi,s)$ associated to a pair $(M,\mu)$ by construction
\ref{constr2}, where $M$ is a free $A$-module of rank one.
Let $A_0:=\coker(\mu:\Md\to A)$ and $S_0:=\Spec(A_0)$.
In other words, $S_0$ is the zero set of the section $\mu$.
If $S_0$ is empty, then $f$ is an isomorphism and the lemma becomes
trivial. So let us assume that $S_0$ is non-empty.
Let $\Rr$ be defined by the cartesian diagram
$$
\xymatrix@R=3ex{
\Rr \ar[rr] \ar[d] & &
\Cc' \ar[d] \\
S_0 \ar@{^{(}->}[r] &
S \ar@{^{(}->}[r]^s &
\Cc
}
$$
Then $\Rr=\Proj(\Sym(H))$ for a free $A_0$-module $H$ of rank two.
Let $\pi'_0:\Rr\to S_0$ be the projection,
$s'_0:S_0\to \Rr$ the section of $\pi'_0$ induced by $s'$,
and let $\E'_{\Rr}$ be the restriction of $\E'$ to $\Rr$.
The admissibility of $\E'$ implies the vanishing of
$H^1(\Rr_x,\E'_{\Rr}(-s')|_{\Rr_x})$ for every fiber $\Rr_x$ of $\pi'_0$.
From this and the flatness of $\pi'_0$ it follows that the sheaf 
$\E'_{\Rr,1}:=(\pi'_0)^*((\pi'_0)_*\E'_{\Rr}(-s'_0))(s'_0)$
is locally free and its formation commutes with base change.
The adjunction morphism gives us a canonical morphism
$\E'_{\Rr,1}\to \E'_{\Rr}$. 
Let $\E'_{\Rr,2}$ be the cokernel of this morphism.
By restricting to fibers one finds that the sequence
$$
0 \to
\E'_{\Rr,1}\to \E'_{\Rr} \to\E'_{\Rr,2}
\to 0
$$
is exact and that
$$
(\E'_{\Rr,1})|_{\Rr_x}\isomorph \bigoplus^d \Oo_{\Rr_x}(1)
\qquad,\qquad
(\E'_{\Rr,2})|_{\Rr_x}\isomorph \bigoplus^{n-d} \Oo_{\Rr_x}
\quad.
$$
for all $x\in S_0$.
Now we choose a direct sum decomposition 
$(s')^*\E'\isomto (M^{\dual}\tensor E_1)\oplus E_2$
compatible with the above exact sequence, i.e.
in such a way that there exist isomorphisms
$(s'_0)^*\E_{\Rr,1}\isomto (M^{\dual}\tensor E_1)|_{S_0}$ and 
$(s'_0)^*\E_{\Rr,2}\isomto E_2|_{S_0}$ 
which make the following diagram commute:
$$
\xymatrix{
0 \ar[r] &
(s'_0)^*\E'_{\Rr,1} \ar[r]\ar[d]^{\isomorph} &
(s'_0)^*\E'_{\Rr} \ar[r]\ar[d]^{\isomorph} &
(s'_0)^*\E'_{\Rr,2} \ar[r]\ar[d]^{\isomorph} & 
0 \\
0 \ar[r] &
(M^{\dual}\tensor E_1)|_{S_0} \ar[r] &
((M^{\dual}\tensor E_1)\oplus E_2)|_{S_0} \ar[r] &
E_2|_{S_0} \ar[r] & 0
}
$$
where the middle vertical arrow is the one induced
from $(s')^*\E'\isomto (M^{\dual}\tensor E_1)\oplus E_2$
by restriction to $S_0$.
Define the locally free $\Oo_{\Cc'}$-module $\F'$
by the exact sequence
$$
0\to\F'\to\E'\to s'_*(M^{\dual}\tensor E_1)\to 0
$$
and let $\F:=f_*\F'$.
It is easy to see that $\F'$ is trivial along the fibers
of $f$. 
By \cite{Knudsen}, Cor. 1.5 the sheaf $\F$ is flat over $S$ and
its formation commutes with base change. It is easy
to see that the restriction of $\F$ to fibers of $\pi$
is locally free of rank $n$. In particular, all the fibers
$\F[x]$ are of dimension $n$. 

I claim that $\F$ itself
is locally free of rank $n$. 
To prove this, let $I\in A$ be the nilradical
of $A$ and let $J:=IB$. Then 
$\overline{S}:=S_{\red}=\Spec(A/I)$
and $\overline{\Cc}:=\Cc\times_S S_{\red}=\Spec(B/J)$.
By \cite{Matsumura}, 22.3 (3) it suffices to show that
$J\tensor_B\F=J\cdot\F$ and that $\F|_{\overline{\Cc}}$ is a flat
$\Oo_{\overline{\Cc}}$-module. The first of these conditions
follows easily from the flatness of $B/A$ and $\F/A$.
Since $\overline{\Cc}\to \overline{S}$ is smooth and $\overline{S}$ 
is reduced, $\overline{\Cc}$ is also reduced 
(cf. \cite{SGA1} Exp. II, Prop 3.1). As observed above,
all fibers of $\F$ (and thus of $\F|_{\overline{\Cc}}$)
have the same dimension $n$.
But a coherent sheaf over a reduced locally noetherian scheme
all of whose fibres have the same dimension is locally free
(cf. \cite{Milne} 2.9, where the ground ring is assumed to
be integer, but the proof works as well for reduced ground
ring). This proves our claim.

I claim furthermore that the canonical morphism $f^*\F\to\F'$ 
is an isomorphism. Indeed, since $\F'$ is flat over $S$, it suffices
to check this on the fibers of $\pi'$, which is easy.

By applying $(s')^*$ to the defining sequence for $\F'$
we obtain the exact sequence
$$
\xymatrix@R=.5ex{
0 \ar[r] &
s^*\Oo_{\Cc}(-s)\tensor E_1 \ar[r] &
(s')^*\F' \ar[rr] \ar[dr] & &
(s')^*\E'\ar[r] &
(M^{\dual}\tensor E_1)\ar[r] & 
0 \\
& & & 
E_2 \ar[ur] \ar[dr] 
& & & \\ 
& & 0 \ar[ur] & & 0 & &
}
$$
Since we have $s^*\F=(s')^*f^*\F=(s')^*\F'$, we get the exact sequence
$$
0 \to s^*\Oo_{\Cc}(-s)\tensor E_1 \to s^*\F \to E_2 \to 0
\quad.
$$
After replacing $S$ and $\Cc$ by members of an open covering,
we can assume that there exists 
a decomposition $\F\isomto \E_1(-s)\oplus\E_2$
together with isomorphisms $E_1\isomto s^*\E_1$ and
$E_2\isomto s^*\E_2$ which is compatible with this exact sequence.
Observe that we have
\begin{eqnarray*}
\ker(\F'\to s'_*E_2)\tensor_{\Oo_{\Cc'}}\Oo_{\Cc'}(s') & \isomorph &
\ker(f^*\F\to s'_*E_2)\tensor_{\Oo_{\Cc'}}\Oo_{\Cc'}(s') \isomorph
\\
& \isomorph & 
(\Md\tensor f^*\E_1)\oplus f^*\E_2 =\E'(\E_1,\E_2)
\end{eqnarray*}
On the other hand, the exact diagram
$$
\xymatrix@R=2ex{
& & 0\ar[d] & & 
\\
& & \E'(-s') \ar[d] & & 
\\
0\ar[r] & \F'\ar[r]\ar[d] & 
\E'\ar[r] \ar[d] & s'_*(M^{\dual}\tensor E_1) \ar[d]^= \ar[r] & 0 
\\
0 \ar[r] & s'_*E_2 \ar[r]\ar[d] & 
s'_*E_2\oplus s'_*(M^{\dual}\tensor E_1) \ar[r]\ar[d] &
s'_*(M^{\dual}\tensor E_1) \ar[r] & 0
\\
& 0 & 0 & &
}
$$
shows that
$\ker(\F'\to s'_*E_2)\tensor_{\Oo_{\Cc'}}\Oo_{\Cc'}(s')\isomorph\E'$.
Therefore $\E'(\E_1,\E_2)\isomorph\E'$.
\end{proof}

From lemmas \ref{special case} and \ref{general case} it follows easily
that the tupel $(M,\ \mu,\ E'\to E,\ M\tensor E'\ot E,\ n-d)$ associated
in construction \ref{constr3}
to an admissible $\Oo_{\Cc'}$-module $\E'$ of degree $d$ 
and isomorphisms $(f_*\E'(-s'))(s)\isomto\E$ and $(s')^*\E'\isomto E'$
is indeed a bf-morphism of rank
$n-d$ from $E'$ to $E$ in the sense of \cite{Kausz} 3.1.
Thus the arrow in proposition \ref{mod&E' to bf} is well defined.

Our aim is now to define the inverse of construction \ref{constr3}.
Let
$$
(M,\ \mu,\ E'\to E,\ M\tensor E' \ot E,\ n-d)
$$
be an arbitrary bf-morphism of degree $n-d$ from $E'$ to $E$.

Assume first that $S$ and $\Cc$ are affine, that $M$ is trivial and
that there exist 
free $\Oo_{\Cc}$-modules $\E_1$ and $\E_2$ of 
rank $d$ and $n-d$ respectively,
and isomorphisms
\begin{eqnarray*}
\varphi &:& \E_1\oplus\E_2 \Isomto \E \\
{\varphi'} &:& (M^{\dual}\tensor E_1)\oplus E_2 \Isomto E'
\end{eqnarray*}
where $E_1:=s^*\E_1$ and $E_2:=s^*\E_2$,
such that the following diagram commutes:
$$
\label{star}
\xymatrix@R=2ex{
E' \ar[rr] & &
E  \ar[rr] & &
M\tensor E'  & \\
& & & & & (*)\quad. \\
(M^{\dual}\tensor E_1) \oplus E_2 
\ar[uu]_{\varphi'}
\ar[rr]^(.6){
\left[
\begin{array}{ll}
\mu & 0 \\
0 & 1 
\end{array}
\right]
}
& &
E_1\oplus E_2
\ar[uu]_{s^*\varphi}
\ar[rr]^(.4){
\left[
\begin{array}{ll}
1 & 0 \\
0 & \mu 
\end{array}
\right]
}
& &
E_1\oplus(M\tensor E_2)
\ar[uu]_{\varphi'}
&
}
$$
Let $\E':=\E'(\E_1,\E_2)$. Composing $\varphi$ and $\varphi'$ with
the natural isomorphisms of lemma \ref{special case} we get isomorphisms
which we denote again by the same symbols:
\begin{eqnarray*}
\varphi &:& (f_*\E'(-s'))(s) \Isomto \E \\
\varphi' &:& (s')^*\E' \Isomto E'
\end{eqnarray*}
If 
\begin{eqnarray*}
\psi &:& \E_1\oplus\E_2 \Isomto \E \\
\psi' &:& (M^{\dual}\tensor E_1)\oplus E_2 \Isomto E'
\end{eqnarray*}
is another pair of isomorphisms such that the diagram analogous to $(*)$
commutes, then we get in the same way a pair of isomorphisms
\begin{eqnarray*}
\psi &:& (f_*\E'(-s'))(s) \Isomto \E \\
\psi' &:& (s')^*\E' \Isomto E'
\end{eqnarray*}

\begin{lemma}                  
\label{unique iso}
In the above situation there exists a unique isomorphism
$
g:\E' \Isomto \E'
$
such that the following diagrams commute
$$
\xymatrix@C=2ex@R=2ex{
(f_*\E'(-s'))(s)
\ar[rr]^{\isomorph}
\ar[rd]_{\varphi}
& &
(f_*\E'(-s'))(s)
\ar[ld]^{\psi}
& &
(s')^*\E' 
\ar[rr]^{\isomorph}
\ar[dr]_{\varphi'}
& &
(s')^*\E'
\ar[dl]^{\psi'}
\\
& \E &
& &
& E' &
}
$$
where the horizontal arrows are the ones induced by $g$.
\end{lemma}

Before proving lemma \ref{unique iso}, we will establish
the following auxiliary result:

\begin{lemma}
\label{aux}
Let $a_{12}:\E_2\to\E_1$ be an $\Oo_{\Cc'}$-module
homomorphism such that the induced morphism $\alpha_{12}:E_2\to E_1$ 
factorizes as follows:
$$
\xymatrix{
E_2 \ar[r]^-{\beta_{12}} & 
M^{\dual} \ar[r]^-{\mu}\tensor E_1 & 
E_1
}
\quad.
$$
Then there exists a factorization 
$$
\xymatrix{
f^*\E_2\ar[r]^-{b_{12}} &
\Md\tensor f^*\E_1\ar[r] &
f^*\E_1
}
$$
of $f^*a_{12}$,
such that $(s')^*b_{12}=\beta_{12}$.
Here the morphism 
$\Md\tensor f^*\E_1\to f^*\E_1$ is the one induced by
the arrow $m:\Oo_{\Cc'}\to\M$ of construction \ref{constr1}.
\end{lemma}

\begin{proof}
After chosing a bases of $\E_1$ and of $\E_2$ one is easily
reduced to the case $\E_1=\E_2=\Oo_{\Cc'}$.
Let $\I$ be the $\Oo_{\Cc}$-module associated to the pair $(M,\mu)$
as in construction \ref{constr2}.
Consider the morphism 
$$
\pi^*M^{\dual} \oplus \Oo_{\Cc}(-s) \to \Oo_{\Cc}
$$
which is the sum of the morphism $\pi^*M^{\dual}\to\Oo_{\Cc}$ 
induced by $\mu$ and the inclusion $\Oo_{\Cc}(-s)\injto\Oo_{\Cc}$.
Obviously, this morphism factorizes through $\I$.
Thus we have a natural morphism $m:\I\to\Oo_{\Cc}$ and
this induces a morphism $\Sc(1)\to\Sc$ of graded $\Sc$-modules, 
where $\Sc$ is the graded $\Oo_{\Cc}$-algebra $\Sc=\Sym(\I)$.
Let $\M:=\Oo_{\Cc'}(-s')\tensor f^*\Oo_{\Cc}(s)$. 
We have seen on page \pageref{O(-1)=M} that 
$\Oo_{\Cc'}(1)=\Md$. It is not difficult to see that
under this identification 
the morphism $\Oo_{\Cc'}(1)\to\Oo_{\Cc'}$ induced by
$\Sc(1)\to\Sc$ corresponds to the morphism dual to 
$\Oo_{\Cc'}\to\M$
defined in construction \ref{constr1}.
Let $S=\Spec(A)$, $\Cc=\Spec(B)$ for a ring $A$ and an $A$-algebra $B$
(remember that we are still dealing with the affine case).
We are given a $B$-module morphism $a_{12}:B\to B$ such that the induced
morphism $\alpha_{12}:A\to A$ factorizes as 
$A\overset{\beta_{12}}{\to} 
M^{\dual}\overset{\mu}{\to} A$.
We have to establish a factorization 
$B\to\I\overset{m}{\to} B$ 
such that the
following diagram commutes:
$$
\xymatrix@R=3ex{
B \ar[r] \ar[d] &
\I \ar[r]^{m} \ar[d] &
B \ar[d] \\
A \ar[r]^{\beta_{12}} &
M^{\dual} \ar[r]^\mu &
A 
}
$$
where the middle vertical arrow is the one defining
$s'$ (cf. construction \ref{constr2}).
Let $b\in B$ be a generator of $\ker(B\to A)$ and let 
$l:M\isomto A$ be an isomorphism.
By assumption we have $a_{12}(1)=b_1l(\mu)+b_2b$ for some
$b_1,b_2\in B$. Furthermore we may choose $b_1, b_2$ in such a way
that $\overline{b_1}l=\beta_{12}(1)$, where $\overline{b_1}\in A$
is the image of $b_1$ under the surjection $B\to A$.
By definition, we may identify $\I$ with
the $B$-module 
$(B\cdot X \oplus B\cdot Y)/(bX-l(\mu)Y)$.
With this identification the morphism $\I\to M^{\dual}$ is given by
$X\mapsto l$, $Y\mapsto 0$ and the morphism
$m:\I\to B$ is given by $X\mapsto l(\mu)$, $Y\mapsto b$.
Therefore the morphism $B\to\I$ defined by $1\mapsto b_1X+b_2Y$
yields the required factorization.
\end{proof}

\begin{proof}(Of lemma \ref{unique iso}).
Let
\begin{eqnarray*}
\left[
\begin{array}{ll}
a_{11} & a_{12} \\
a_{21} & a_{22}
\end{array}
\right]
&:=&
\psi\comp\varphi^{-1}
:
\E_1\oplus\E_2 \isomto \E_1\oplus\E_2
\\
\left[
\begin{array}{ll}
\alpha_{11} & \alpha_{12} \\
\alpha_{21} & \alpha_{22}
\end{array}
\right]
&:=&
s^*(\psi\comp\varphi^{-1})
:
E_1\oplus E_2 \isomto E_1\oplus E_2
\\
\left[
\begin{array}{ll}
\beta_{11} & \beta_{12} \\
\beta_{21} & \beta_{22}
\end{array}
\right]
&:=&
\psi'\comp(\varphi')^{-1}
:
(M^{\dual}\tensor E_1)\oplus E_2 \isomto 
(M^{\dual}\tensor E_1)\oplus E_2 
\end{eqnarray*}
Let $\E'_1:=\Md\tensor f^*\E_1$ and $\E'_2:=f^*\E_2$.
Then $\E'=\E'_1\oplus\E'_2$ and a morphism $\E'\to\E'$
is of the form 
$$
\left[
\begin{array}{ll}
b_{11} & b_{12} \\
b_{21} & b_{22}
\end{array}
\right]
:
\E'_1\oplus\E'_2 \to \E'_1\oplus\E'_2
\quad.
$$
We have to show that there exist morphisms 
$b_{ij}\in\Hom(\E'_j,\E'_i)$
with
$$
f_*(\M\tensor b_{ij})= a_{ij}
\qquad \text{and} \qquad
(s')^*b_{ij} = \beta_{ij}
$$
for $i,j\in\{1,2\}$ and that they are uniquely determined by 
this property. 

We let 
$b_{11} := \Md\tensor f^*a_{11}$ and
$b_{22} := f^*a_{22}$
and define $b_{21}$ as the composition
$$
\xymatrix{
b_{21}\ :\ \Md\tensor f^*\E_1 \ar[r]^-{m} &
f^*\E_1 \ar[r]^-{f^*a_{21}} &
f^*\E_2
}
\quad.
$$
It is easy to check that the commutativity of the 
diagram $(*)$ on page \pageref{star} implies
$$
\xymatrix@R=1ex{
\alpha_{11} = \beta_{11} &
\alpha_{12} = \mu\beta_{12} \\
\mu\alpha_{21} = \beta_{21} &
\alpha_{22} = \beta_{22} 
}
$$
By lemma \ref{aux} there exists a morphism 
$b_{12}:f^*\E_2\to\Md\tensor f^*\E_1$
with $mb_{12}=f^*a_{12}$ and $(s')^*b_{12}=\beta_{12}$.
One checks without difficulty that the $b_{ij}$ thus defined have 
the required property. We will see below that the matrix
$[b_{ij}]$ defines an isomorphism $\E'\isomto\E'$.

The uniqueness of $b_{ij}$ is clear 
for $(i,j)\in\{(1,1),(2,1),(2,2)\}$, since in these cases
the morphism 
$\Hom(\E'_j,\E'_i)\to\Hom(\E_j,\E_i)$
given by $g\mapsto f_*(\M\tensor g)$
is an isomorphism.
For the uniqueness of $b_{12}$ we have to show the injectivity
of the morphism
$$
v':
\left\{
\begin{array}{ll}
\Hom(f^*\E_2,\Md\tensor f^*\E_1) &
\to\quad \Hom(\E_2,\E_1)\quad\oplus\quad\Hom(E_2,M^{\dual}\tensor E_1) \\
\qquad\qquad g &
\mapsto\quad (f_*(\M\tensor g)\ ,\ (s')^*g)
\end{array}
\right.
$$
Let $v_1:f_*\Md\to\Oo_{\Cc}$ and
$v_2:f_*\Md\to s_*M^{\dual}$    be the indicated arrows in the 
commutative diagram
with exact rows
$$
\xymatrix@R=3ex{
0 \ar[r] &
\Oo_{\Cc}(-s) \ar[r] \ar[d] &
f_*\Md \ar[r]^{v_2} \ar[d] \ar@{.>}[dl]_{v_1} &
s_*M^{\dual} \ar[r] \ar[d]^0 &
0 \\
0 \ar[r] &
\Oo_{\Cc} \ar[r]&
(f_*\Md)(s) \ar[r] &
(s_*M^{\dual})(s) \ar[r] &
0
}
$$
We have canonical identifications
\begin{eqnarray*}
& &\Hom(f^*\E_2,\Md\tensor f^*\E_1)= 
\Hom(\E_2,\E_1)\tensor f_*\Md \\
&\text{and}& 
\Hom(E_2,M^{\dual}\tensor E_1)=
\Hom(\E_2,\E_1)\tensor_{\Oo_{\Cc}}s_*M^{\dual}
\end{eqnarray*}
and it is not difficult to see that the morphism $v'$ coincides
with the morphism
$$
\Hom(\E_2,\E_1)\tensor f_*\Md \to
\Hom(\E_2,\E_1)\tensor(\Oo_{\Cc}\oplus s_*M^{\dual})
$$
obtained by tensoring 
$v:=v_1+v_2:f_*\Md\to\Oo_{\Cc}\oplus s_*M^{\dual} $ 
whith $\Hom(\E_2,\E_1)$.
Therefore the injectivity of $v'$ is equivalent to the injectivity
of $v$. But the injectivity of $v$ is clear, since
$$
\ker(v)=\ker(v_1)\cap\ker(v_2)=\ker(\Oo_{\Cc}(-s)\to\Oo_{\Cc})=0
\quad.
$$

It is now easy to see that the morphism $g:\E'\to\E'$ defined
by the matrix $[b_{ij}]$ is an isomorphism: Let $g':\E'\to\E'$ be
the morphism constructed in the same way as $g$, only with the roles of 
$(\varphi,\varphi')$ and $(\psi,\psi')$ interchanged. 
By uniqueness it follows that $g\comp g'=g'\comp g=\id$.
\end{proof}

\begin{construction}
\label{constr4}
We return now to the general case and let 
$$
(M,\ \mu,\ E'\to E,\ M\tensor E' \ot E,\ n-d)
$$
be a bf-morphism of rank $n-d$ from $E'$ to $E$.
With the help of lemma 4.1 in \cite{Kausz} it is easy to see
that there exists an open affine covering $\Cc=U_0\cup\bigcup_{i\in I} U_i$
with the following properties:
\begin{enumerate}
\item
$U_0$ is the complement of the section $s$ in $\Cc$.
\item
For each $i\in I$, $s$ induces a section of 
$\pi|_{U_i}:U_i\to\pi(U_i)$. In particular, $\pi(U_i)$ is affine.
\item
For each $i\in I$ the invertible sheaf $M$ is trivial over $\pi(U_i)$.
\item
For each $i\in I$ there exist free $\Oo_{U_i}$-modules $\E_1^i$ and $\E_2^i$
of rank $d$ and $n-d$ respectively, and isomorphisms
\begin{eqnarray*}
\varphi_i &:& \E_1^i\oplus\E_2^i \Isomto \E|_{U_i} \\
\varphi_i' &:& (M^{\dual}\tensor E_1^i)\oplus E_2^i 
\Isomto E'|_{s^{-1}(U_i)}
\end{eqnarray*}
where $E_1^i:=s^*\E_1^i$ and $E_2^i:=s^*\E_2^i$, such that
over $s^{-1}(U_i)$ the diagram analogous to $(*)$ (cf. page
\pageref{star}) commutes.
\end{enumerate}
For $i\in I$ we have over $V_i:=f^{-1}(U_i)$ the 
locally free $\Oo_{V_i}$-modules
$$
\E'_i:=\E'(\E_1^i,\E_2^i)=
(\Oo_{\Cc'}(s')\tensor f^*\E_1^i(-s))\oplus f^*\E_2^i
$$
together with the isomorphisms
\begin{eqnarray*}
\varphi_i &:& (f_*\E'_i(-s'))(s) \Isomto \E|{U_i} \\
\varphi'_i&:&
(s')^*\E'_i \Isomto E'|_{s^{-1}(U_i)}
\end{eqnarray*}
Let $U:=\bigcup_{i\in I}U_i$ and $V:=f^{-1}(U)$. Then
$(V,f|_V,\pi|_U,s)$ is a simple modification of $(U,\pi|_U,s)$
and by the lemmas \ref{special case} and \ref{unique iso} the data
$(\E'_i,\varphi_i,\varphi'_i)$ 
glue together to give an admissible $\Oo_{V}$-module $\E'_V$ 
of degree $d$ together with isomorphisms
$\varphi_V:((f|_V)_*\E'_V(-s'))(s)\Isomto\E|_U$ and 
$\varphi':(s')^*\E'_V\Isomto E'$.
Let $V_0:=f^{-1}(U_0)$ and $\E'_0:=(f|_{V_0})^*\E$.
Since the restriction of $f$ to $V_0$ is an isomorphism onto $U_0$,
we may identify $V_0$ with $U_0$,\ $\E'_0$ with $\E|_{U_0}$,
and $((f|_V)_*\E'_V(-s'))(s)|_{V_0\cap V}$ with 
$\E'_V|_{V_0\cap V}$. Thus $\E'_V$ and $\E'_0$ glue together
via $\varphi_V$ to give an admissible $\Oo_{\Cc'}$-module
$\E'$ of degree $d$ together with isomorphisms
$(f_*\E'(-s'))(s)\Isomto \E$ and $(s')^*\E'\Isomto E'$.
By lemma \ref{unique iso}, this construction does not depend on
the choice of the $U_i, \varphi_i, \varphi'_i$.
\end{construction}

\begin{remark}
Construction \ref{constr4} has the drawback that it involves a choice
of the data $U_i, \varphi_i, \varphi'_i$.
I have tried in vain to find a direct construction, which would produce a
canonical graded $\Sym(\I)$-module out of a bf-morphism from $E'$ to $E$ such 
that the associated $\Oo_{\Cc'}$-module would be $\E'$.
\end{remark}

It is an immediate consequence of lemmas \ref{special case},
\ref{general case} and \ref{unique iso} that the
constructions \ref{constr3} and \ref{constr4} are  
inverses of each other. This proves proposition 
\ref{mod&E' to bf}.

\section{Contractions}
\label{section contractions}

The following proposition is implicit in the
proof of 
theorem 2.1 in \cite{Knudsen}.
\begin{proposition}
\label{X^c}
Let $Y$ be a locally noetherian scheme, 
let $\pi:X\to Y$ be
a proper flat morphism with $\pi_*\Oo_X=\Oo_Y$
and let $\Ll$ be an invertible $\Oo_X$-module. 
For $y\in Y$ we denote by $X_y$ the fibre
of $\pi$ over $y$ and by $\Ll_y$ the restriction of $\Ll$
to $X_y$. Assume that for every $y\in Y$ the following holds:
\begin{enumerate}
\item
$H^1(X_y,\Ll_y^{\tensor i})=(0)$ for $i\geq 1$.
\item
The $\Oo_{X_y}$-module $\Ll_y$ is generated by global sections.
\item
For every $i\geq 1$ the canonical morphism 
$(H^0(X_y,\Ll_y))^{\tensor i}\to H^0(X_y,\Ll_y^{\tensor i})$
is surjective.
\end{enumerate}
Let $\Sc$ be the graded $\Oo_Y$-algebra 
$\oplus_{i\geq 0}\pi_*\Ll^{\tensor i}$ and let 
$X^c:=\Proj(\Sc)$.
Then the canonical morphism $\pi^*\pi_*\Ll\to\Ll$
(resp. $\Sym\pi_*\Ll\to\Sc$) 
is surjective
and therefore induces a morphism 
$p:X\to\Proj(\Sym\pi_*\Ll)$ 
(resp. a closed immersion $X^c\injto\Proj(\Sym\pi_*\Ll)$).
The scheme theoretic image of $p$ is $X^c$ and $X^c$ is
flat over $Y$.
\end{proposition}

\begin{proof}
We set $\Sc^i:=\pi_*\Ll^{\tensor i}$ for $i\geq 0$.
It is well known that
property 1. ensures that the sheaves $\Sc^i$ are locally
free and that their construction commutes with any base
change $Y'\to Y$. By definition, $\Sc$ is the graded algebra with
components $\Sc^i$.
An appplication of Nakayama's lemma shows that
property 2. and property 3. entails the surjectivity
of $\pi^*\Sc^1\to\Ll$ and of $\Sym\Sc^1\to\Sc$ respectively.
The first of these surjections induces a morphism
$p:X\to\Proj(\Sym\Sc^1)$ (cf. \cite{EGA II}(4.3.2.)),
the second induces a closed immersion $X^c\injto\Proj(\Sym\Sc^1)$
(cf. \cite{EGA II}(2.9.2.)(i)).

To prove that $X^c$ is the scheme-theoretic image of $p$, we
may assume that $Y$ is affine. Let $S^i:=\Gamma(Y,\Sc^i)$ and
$S:=\Gamma(Y,\Sc)$.
For $f\in S^1$ we consider the open subscheme 
$U_f:=\Spec(\Sym S^1)_{(f)}$ of $\Proj(\Sym S^1)$.
Its preimage under $p:X\to\Proj(\Sym S^1)$ is $X_f:=\{x\in X\ |\ f(x)\neq 0\}$
and its preimage under $X^c\injto\Proj(\Sym S^1)$ is $V_f:=\Spec(S_{(f)})$.
The restriction of $f$ to $X_f$ defines a trivialization
$\Ll|_{X_f}\isomto\Oo_{X_f}$ and thus a morphism
$S^1=\Gamma(X,\Ll)\to\Gamma(X_f,\Ll)\isomto\Gamma(X_f,\Oo_X)$
which maps $f$ to $1$.
The induced morphism $(\Sym S^1)_{(f)}\to\Gamma(X_f,\Oo_X)$ is
the defining one for the morphism 
$p|_{X_f}:X_f\to U_f\subset\Proj(\Sym S^1)$.
We have a canonical morphism 
$$
\left\{
\begin{array}{ll}
S_{(f)} &\To  \Gamma(X_f,\Oo_X) \\
a/f^n   &\To  \quad b
\end{array}
\right.
$$
where $a\in S^n=\Gamma(X,\Ll^{\tensor n})$
and $b\in\Gamma(X_f,\Oo_X)$ is the element
defined by the relation 
$a|_{X_f}=b\cdot f^n|_{X_f}\in\Gamma(X_f,\Ll^{\tensor n})$.
It is easy to see that this morphism is injective and that
the diagram
$$
\xymatrix@R=1.5ex{
(\Sym S^1)_{(f)} \ar[rr] \ar@{->>}[dr] & & \Gamma(X_f,\Oo_X) \\
& S_{(f)} \ar@{^(->}[ur] &
}
$$
commutes. This shows that $V_f$ is the scheme-theoretic image 
of $X_f$ in $U_f$ (cf. \cite{EGA I} (6.10)).

The flatness of $X^c\to Y$ is clear, since as we already
remarked, all the $\Sc^i$ ($i\geq 0$) are locally free $\Oo_Y$-modules.
\end{proof}

\begin{lemma}
\label{pi_*(O) and h_*(O)}
Let $S$ be a $k$-scheme, let $(\Cc,\pi,s_1,s_2,h)$ be a modification
of $(\tilde{C}_0,p_1,p_2)$ over $S$ (cf. \ref{modification/S of C0,p1,p2}).
Then the canonical morphisms $\Oo_S\to\pi_*\Oo_{\Cc}$ and 
$\Oo_{\tilde{C}\times S}\to h_*\Oo_{\Cc}$ are isomorphisms.
\end{lemma}

\begin{proof}
By \cite{Knudsen}, Cor. 1.5 the sheaf $h'_*\Oo_{\Cc'}$ is 
locally free and commutes with any base change $S'\to S$.
It is easy to see that when restricted to any geometric fibre
of the projection $\pr_2:\tilde{C}_0\times S\to S$,
the morphism
$\Oo_{\tilde{C}_0\times S}\to h'_*\Oo_{\Cc'}$ becomes an isomorphism.
Therefore $\Oo_{\tilde{C}_0\times S}\to h'_*\Oo_{\Cc'}$ is an
isomorphism. Clearly we have $(\pr_2)_*\Oo_{\tilde{C}_0\times S}\isomto\Oo_S$.
Consequently we have 
$\pi'_*\Oo_{\Cc'}=(\pr_2)_*( h'_*\Oo_{\Cc'})\isomto\Oo_S$.
\end{proof}

\begin{definition}
Let $Y$ be a locally noetherian connected scheme and let $\pi:X\to Y$
be a proper flat curve over $Y$. For a locally free $\Oo_X$-module
$\E$ we let
$\deg_{X/Y}\E:=\chi(\det\E_y)-\chi(\Oo_{X_y})$
for some $y\in Y$ where $\E_y$ is the restriction of $\E$ to the fibre
$X_y$ of $\pi$ over $y$. Note that by \cite{EGA III} (7.9.4) this number
is independent of the choice of $y$. In case $Y$ is the $\Spec$ of
a field, we recover definition \ref{definiton of degree}.    
\end{definition}

\begin{lemma}
\label{X_1+X_2}
Let $K$ be a field and let $X/K$ be a proper curve over $K$.
Let $x\in X(K)$ be an ordinary double point of $X$ such
that the blowing up of $X$ in $x$ is the disjoint union
of two curves $X_1$ and $X_2$. Let $\Ll$ be an invertible
$\Oo_X$-module and let $\Ll_\nu:=i^*_\nu\Ll$, where $i_\nu:X_\nu\to X$
is the canonical closed immersion ($\nu=1,2$).
Then the following holds:
\begin{enumerate}
\item
$\deg_{X/K}\Ll=\deg_{X_1/K}\Ll_1+\deg_{X_2/K}\Ll_2$
\item
Assume that either $\Ll_1$ or $\Ll_2$ 
is generated by global sections.
Then we have $H^1(X,\Ll)=H^1(X_1,\Ll_1)\oplus H^1(X_2,\Ll_2)$.
\item
If $\Ll_\nu$ is generated by global sections and
the canonical morphisms 
$\Gamma(X_\nu,\Ll_\nu)^{\tensor m}\to\Gamma(X_\nu,\Ll_\nu^{\tensor m})$
are surjective for all $m\geq 0$ ($\nu=1,2$),
then $\Ll$ is generated by global sections and the canonical morphisms
$\Gamma(X,\Ll)^{\tensor m}\to\Gamma(X,\Ll^{\tensor m})$
are surjective for all $m\geq 0$.
\item
If $\Ll_\nu$ is very ample for $\nu=1,2$, then also
$\Ll$ is very ample.
\end{enumerate}
\end{lemma}

\begin{proof}
The first two assertions follow easily from the following exact sequence
of $\Oo_X$-modules:
$$
0\To
\Ll\To
(i_1)_*\Ll_1\oplus(i_2)_*\Ll_2\To
\Ll[x]
\To 0
\quad.
$$
The proof of the third assertion is entirely similar to the proof of
theorem 1.8 c) in \cite{Knudsen}.
Recall that on a proper $K$-scheme $X$ an invertible $\Oo_X$-module
$\Ll$ is very ample if and only if the canonical morphisms
$\Gamma(X,\Ll)\to\Ll[x]\oplus\Ll[y]$ and $\Gamma(X,\Ll)\to\Ll_x/\m_x^2\Ll_x$
are surjective for all $x\neq y\in X$.
Using this criterion the fourth assertion is easy to see.
\end{proof}

\begin{lemma}
\label{Ll}
Let $S$ be a locally noetherian 
connected $k$-scheme, let $(\Cc',\pi',s'_1,s'_2,h')$
be a modification of $(\tilde{C}_0,p_1,p_2)$ over $S$ 
(cf. \ref{modification/S of C0,p1,p2})
and let $\E'$ be a locally free $\Oo_{\Cc'}$-module of rank $n$ which is
admissible for $(\Cc',\pi',s'_1,s'_2,h')$ 
(cf. \ref{admissible for Cc,pi,s1,s2,f}).
Let $\Ll_0$ be an invertible $\Oo_{\tilde{C}}$-module
of degree $\geq 2(2g+1)+2n-\deg_{\Cc'/S}\E'$ 
where $g$ is the genus of $\tilde{C}_0$
and let $d\in\{0,\dots,n\}$.
Then with respect to the morphism $\pi'$
the  invertible $\Oo_{\Cc'}$-module
$$
\Ll:=(det\E')(-ds_1')\tensor_{\Oo_{\Cc'}}  (h')^*(\Ll_0\extensor_k\Oo_S)
$$
has the properties 1-3 listed in proposition
\ref{X^c}. Furthermore, if $Z$ is the component of a fibre
$\Cc'_z$ of $\pi'$ which by $h'$ is mapped isomorphically
to $\tilde{C}_0\tensor_k\kappa(z)$, then $\Ll|_Z$ is very ample. 
\end{lemma}

\begin{proof}
It is clearly sufficient to consider the case where $S$ is the
spectrum of an algebraically closed field $K$. 
Let $\Rr_i:=(h')^{-1}(p_i)\subset\Cc'$
be the chain of rational curves mapped to to $p_i$ ($i=1,2$)
and let $Z\subseteq\Cc'$ be the component of
$\Cc'$ which is mapped isomorphically to $\tilde{C}_0\times S$.
Since $\E|_{\Rr_i}$ is strictly positive, it follows easily from
\ref{X_1+X_2}.2.-3. 
that $H^1(\Rr_i,\Ll|_{\Rr_i})=0$,
that $\Ll|_{\Rr_i}$ is generated by global sections
and that the morphisms
$\Gamma(\Rr_i,\Ll|_{\Rr_i})^{\tensor m}
 \to\Gamma(\Rr_i,(\Ll|_{\Rr_i})^{\tensor m})$
are surjective for $m\geq 0$ ($i=1,2$).
By \ref{X_1+X_2} it is clear therefore that it suffices to show 
that $H^1(Z,\Ll|_Z)=0$,
that $\Ll|_Z$ is very ample
and that the morphism 
$\Gamma(Z,\Ll|_Z)^{\tensor m}\to\Gamma(Z,(\Ll|_Z)^{\tensor m})$
is surjective for all $m\geq 0$.
By \ref{X_1+X_2}.1. we have 
$
\deg_{Z/K}(\E|_{Z})=\deg_{\Cc'/K}(\E)-\sum_{i=1}^2\deg_{\Rr_i/K}(\E|_{\Rr_i})
\geq \deg_{\Cc'/K}(\E)-n
$.
Therefore we have
$$
\deg_{Z/K}(\Ll|_{Z})\geq
\deg_{Z/K}(\E|_{Z})-d+\deg_{\tilde{C}_0/k}\Ll_0\geq 2(2g+1)
\quad.
$$
It is well known that this lower bound for the degree of $\Ll_Z:=\Ll|_Z$ 
implies very ampleness and the vanishing of the first cohomology 
of  $\Ll_Z$. The morphisms 
$\Gamma(Z,\Ll_Z)^{\tensor m}\to\Gamma(Z,\Ll_Z^{\tensor m})$
are trivially surjective for $m=0,1$.
Let $\Ll_1$, $\Ll_2$ be two invertible $\Oo_Z$-modules with
$\Ll_1\tensor\Ll_2=\Ll_Z$ and $\deg_{Z/K}\Ll_i\geq 2g+1$ for $i=1,2$.
For $m\geq 1$ consider the commutative diagram
$$
\xymatrix@R=3ex{
\Gamma(Z,\Ll_Z^{\tensor m})\tensor
\Gamma(Z,\Ll_1)\tensor\Gamma(Z,\Ll_2)
\ar[d]^\alpha \ar[r] 
&
\Gamma(Z,\Ll_Z^{\tensor m})\tensor
\Gamma(Z,\Ll_Z)
\ar[d]^\gamma
\\
\Gamma(Z,\Ll_Z^{\tensor m}\tensor\Ll_1)\tensor
\Gamma(Z,\Ll_2)
\ar[r]^\beta
&
\Gamma(Z,\Ll_Z^{\tensor m+1})
}
$$
By the generalized lemma of Castelnuovo (cf. \cite{Knudsen}, p.170)
the vanishing of
$H^1(Z,\Ll_Z^{\tensor m}\tensor\Ll_1^{-1})$ and of
$H^1(Z,\Ll_Z^{\tensor m}\tensor\Ll_1\tensor\Ll_2^{-1})$
implies the surjectivity of $\alpha$ and $\beta$ respectively.
The surjectivity of $\gamma$ now follows.
By induction this shows that 
$\Gamma(Z,\Ll_Z)^{\tensor m}\to\Gamma(Z,\Ll_Z^{\tensor m})$
is indeed surjective for all $m\geq 0$.
\end{proof}

\begin{proposition}
\label{contraction}
Let $S$ be a $k$-scheme, let $(\Cc',\pi',s'_1,s'_2,h')$ be a modification
of $(\tilde{C}_0,p_1,p_2)$ over $S$ (cf. \ref{modification/S of C0,p1,p2})
and let $\E'$ be a locally free $\Oo_{\Cc'}$-module of rank $n$ which is
admissible for $(\Cc',\pi',s'_1,s'_2,h')$ 
(cf. \ref{admissible for Cc,pi,s1,s2,f}).
Assume that $\E'$ is of extremal degree $\geq(d_1,d_2)$ 
(cf. \ref{admissible for Cc,pi,s1,s2,f}).
There exists a modification $(\Cc,\pi,s_1,s_2,h)$ of $(\tilde{C}_0,p_1,p_2)$
over $S$ and an $S$-morphism $f:\Cc'\to\Cc$ such that the
following holds:
\begin{enumerate}
\item
$h\comp f=h'$ and $f\comp s_i=s'_i$ for $i=1,2$.
\item
The diagram
$$
\xymatrix@R=1.5ex{
\Cc' \ar[rr]^{f} \ar[dr]^{\pi'} & & 
\Cc \ar[dl]_{\pi} \\
   & S \ar@/^/[ul]^{s_1'} \ar@/_/[ur]_{s_1} &
}
$$
defines a simple modification of $(\Cc,\pi,s_1)$ in the sense of 
\ref{def mpc} and the bundle $\E'$ is admissible of degree $d_1$ 
(in the sense of \ref{admissible of degree d})
for $(\Cc',f,\pi',s'_1)$.
\item
Let 
$$
\E:=(f_*\E'(-s'_1))(s_1)
\quad.
$$
Then $\E$ is a locally free 
$\Oo_{\Cc}$-module of rank $n$
and is admissible for $(\Cc,\pi,s_1,s_2,h)$. The extremal degree
of $\E$ is $\geq (d_1+1,d_2)$.     
\end{enumerate}
Furthermore, by these properties the data $(\Cc,\pi,s_1,s_2,h)$ 
and $f:\Cc'\to\Cc$ are uniquely determined in the following sense:
Let $(\Cc^0,\pi^0,s_1^0,s_2^0,h^0)$ and $f^0:\Cc'\to\Cc^0$ be another
set of data of the same kind. Then there exists a unique isomorphism
$\Cc\isomto\Cc^0$ such that the obvious diagrams commute.
\end{proposition}

\begin{proof}
We first prove the uniqueness of 
$(\Cc,\pi,s_1,s_2,h)$ and $f:\Cc'\to\Cc$.
Let $z\in S$. By properties \ref{contraction}.2. and \ref{contraction}.3.
the morphism $f_z:\Cc'_z\to\Cc_z$ induced
by $f$ between the fibres over $z$ is not an 
isomorphism if and only if $\E'|_{\Cc'_z}$ is
of extremal degree $(d_1,d)$ for some $d\geq d_2$.
The same holds for $f^0$, where
$(\Cc^0,\pi^0,s_1^0,s_2^0,h^0)$ and $f^0:\Cc'\to\Cc^0$ is another
set of data satisfying \ref{contraction}.1.-\ref{contraction}.3.
Therefore $f^0$ is constant on the fibres of $f$.
The morphism $f$ is proper by assumption and the morphism
$\Oo_{\Cc}\to f_*\Oo_{\Cc'}$ is an isomorphism by \ref{f_*}.
Therefore we may apply \cite{EGA II} (8.11.1) to conclude that
there is a unique morphism $u:\Cc\to\Cc^0$ such that $f^0=u\comp f$.
With the same argument, there exists a unique morphism 
$v:\Cc^0\to\Cc$ with $f=v\comp f^0$ and we have 
$u\comp v=\id_{\Cc^0}$ and $v\comp u=\id_{\Cc}$. In particular,
$u$ is an isomorphism.
It is clear that $u\comp s_i=s_i^0$ for $i=1,2$
and it follows immediately that $\pi = \pi^0\comp u$.
Since the morphism $h$ is clearly constant on the fibres of $h^0\comp u$
and furthermore the morphism $h^0$ is proper by assumption and we have 
$h^0_*\Oo_{\Cc^0}=\Oo_{\tilde{C}\times S}$ by \ref{pi_*(O) and h_*(O)},
we may again apply \cite{EGA II} (8.11.1) to conclude that
there is a unique morphism $w:\tilde{C}\times S\to\tilde{C}\times S$
such that $h=w\comp h^0\comp u$. Then we have 
$h'=h\comp f=w\comp h^0\comp u\comp f=w\comp h'$ and
yet another application of \cite{EGA II} (8.11.1) leads
to the conclusion that $w$ is in fact the identitiy on $\tilde{C}\times S$,
in particular that $h=h^0\comp u$.
This proves the uniqueness statement.

For the existence it suffices by \ref{GVBD/noetherian} to consider
the case where $S$ is locally noetherian. Furthermore we may assume
$S$ to be connected. 
By \ref{Ll} there is a very ample invertible $\Oo_{\tilde{C}_0}$-module
$\Ll_0$ such that the $\Oo_{\Cc'}$-module
$$
\Ll:=(det\E')(-d_1s_1')\tensor_{\Oo_{\Cc'}}  (h')^*(\Ll_0\extensor_k\Oo_S)
$$
has the properties 1.-3. listed in \ref{X^c} with respect to $\pi'$.
Let $\Cc\subset\Proj(\Sym \pi'_*\Ll)$ be the scheme-theoretic image
of the morphism $\Cc'\to\Proj(\Sym \pi'_*\Ll)$ induced by the surjection
$(\pi')^*\pi'_*\Ll\to\Ll$.
We have the commutative diagram
$$
\vcenter{
\xymatrix@R=1.5ex{
\Cc' \ar[rr]^{f} \ar[dr]^{\pi'} & & 
\Cc \ar[dl]_{\pi} \\
   & S \ar@/^/[ul]^{s_1'} \ar@/_/[ur]_{s_1} &
}}
\eqno(*)
$$
where $\pi:\Cc\to S$ is the structure morphism and
$s_i:=f\comp s'_i$ for $i=1,2$.
Since by \ref{pi_*(O) and h_*(O)} we have $\pi'_*\Oo_{\Cc'}=\Oo_S$,
we can apply proposition \ref{X^c} to conclude that $\pi$ is 
a flat morphism. Since $\pi$ and $\pi'$ are obviously proper,
the morphism $f$ is also proper.
Let $z\in S$ and let $\Cc'_z$, $\Cc_z$, $f_z$ etc. denote the
objects over $\Spec\kappa(z)$ induced by pull back.
By definition, the curve $\Cc'_z$ is of the following shape:
\vspace{3mm}
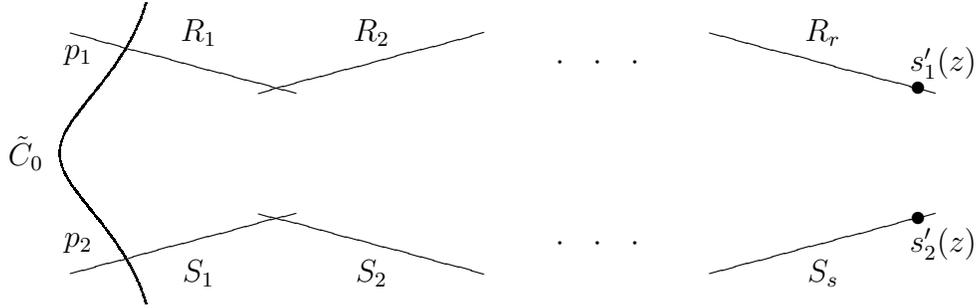
\begin{center}
\parbox{8cm}{
\xy <0mm,-20mm>; <1mm,-20mm>:
(10,20);(10,-20)**\crv{(7,10)&(-10,0)&(7,-10)};
(0,16);(30,8) **@{-}; (17,16) *{\text{$R_1$}};
(25,8);(55,16) **@{-}; (40,16) *{\text{$R_2$}};
(65,12)*{\cdot};
(70,12)*{\cdot};
(75,12)*{\cdot};
(85,16);(115,8) **@{-} ?>>>*{\bullet}; (100,16) *{\text{$R_{r}$}};
(0,-16);(30,-8) **@{-}; (17,-16) *{\text{$S_1$}};
(25,-8);(55,-16) **@{-}; (40,-16) *{\text{$S_{2}$}};
(65,-12)*{\cdot};
(70,-12)*{\cdot};
(75,-12)*{\cdot};
(85,-16);(115,-8) **@{-} ?>>>*{\bullet}; (100,-16) *{\text{$S_{s}$}};
(-6,0) *{\text{$\tilde{C}_0$}};
(1,13) *{\text{$p_1$}};
(1,-12) *{\text{$p_2$}};
(116,12) *{\text{$s'_1(z)$}};
(116,-12) *{\text{$s'_2(z)$}};
\endxy
}
\end{center}
\vspace{3mm}
for some $r,s\geq 0$, with $R_i\isomorph S_j\isomorph\Pp_{\kappa(z)}^1$.
Let 
$D\subseteq\Cc'_z$
be the closed subscheme of $\Cc'_z$ consisting of all components except
$R_r$. It is not difficult to check that $\Ll_z|_D$ is very ample
and that $\Ll_z$ is very ample if $d_1<\deg\E'_z|_{R_r}=:d$.
Furthermore it is clear that $\Ll_z|_{R_r}$ is trivial if 
$d_1=d$.
It follows that there are two possible cases:
\begin{enumerate}
\item
$d_1<d$. Then $f_z:\Cc'_z\to\Cc_z$ is an isomorphism.
\item
$d_1=d$. Then $\Cc_z\isomorph D$ and $f_z$ is the contraction of 
the rational curve $R_r$.
\end{enumerate}
It is now clear that $\Cc$ is a flat curve over $S$, that the
$s_i$ are sections of $\pi$ which meet $\Cc$ in points where $\pi$
is smooth, that the diagram $(*)$ defines a simple modification
of $(\Cc,\pi,s_1)$ and that $\E'$ is admissible of degree $d_1$ for 
$(\Cc',f,\pi',s'_1)$. 

In particular, the morphism $h'$ is constant
on the fibres of $f$. Since $f$ is proper and 
we have $f_*\Oo_{\Cc'}=\Oo_{\Cc}$ by \ref{f_*}, 
we can apply \cite{EGA II} (8.11.1) to show that there exists a 
unique morphism $h:\Cc\to\tilde{C}_0\times S$ with $h'=h\comp f$.
It is clear from the fibrewise consideration above 
that $(\Cc,\pi,s_1,s_2,h)$ is a modification of
$(\tilde{C}_0,p_1,p_2)$.

That $\E:=(f_*\Oo_{\Cc'}(-s'_1))(s_1)$ 
is locally free of rank $n$, follows from  
\ref{general case} and \ref{special case}. For the remaining part
of property 3. it suffices to consider the case $S=\Spec K$ for
some field $K$. Let $R'\subset \Cc'$ be the fibre $f^{-1}(s_1)$.
Since otherewise $f$ is an isomorphism, we may assume that
$R'$ is a chain of rational curves of length $r\geq 1$. Let
$R_i$ $(i=1,\dots,r)$ be its successive irreducible components.
For $i=1,\dots,r-1$ let $x_i\in R'(K)$ be the point where  $R_i$ and 
$R_{i+1}$ intersect. Let $x_r:=s'_1\in R_r(K)$ and let 
$x_0:=p_1\in R_1(K)$ be the 
point, where $R'$ meets the rest of $\Cc'$. We identify $\Cc$ with the 
closed subset of $\Cc'$ consisting of the union of all irreducible
components of $\Cc'$, except $R_r$. Therefore we consider
$R:=\cup_{i=1}^{r-1}R_i$ also as a closed subscheme of $\Cc$ and we have
$x_{r-1}=s_1$. For $i=1,\dots,r$ (resp. $i=1,\dots,r-1$) 
let $\E'_i:=\E'|_{R_i}$  and $\delta'_i:=\deg(\E'_i)$
(resp. $\E_i:=\E|_{R_i}$ and $\delta_i:=\deg(\E_i)$).
Since $\E'|_{R'}$ is strictly standard,
we have canonical exact sequences
$$
0\to\F'_i\to\E'_i\to\G'_i\to 0
$$ 
for all $i=1,\dots,r$, 
where $\F'_i\isomorph\oplus^{\delta'_i}\Oo_{R_i}(1)$ and 
$\G'_i\isomorph\oplus^{n-\delta'_i}\Oo_{R_i}$.
It is easy to see that 
$f_*\E'(-s'_1)=\ker(\E'\to\E'[x_{r-1}]\to\G_r[x_{r-1}])$.
Therefore we have an exact diagram of $\Oo_{R_{r-1}}$-modules
as follows:
$$
\xymatrix{
  &  0\ar[d]  &  0\ar[d]  &  ?\ar[d] &  
\\
0\ar[r] & 
\F'_{r-1}(-x_{r-1}) \ar[r] \ar[d] &
\F'_{r-1} \ar[r] \ar[d] &
\F'_{r-1}[x_{r-1}] \ar[r] \ar[d] & 0
\\
0\ar[r] & 
\E_{r-1}(-x_{r-1}) \ar[r] \ar[d] &
\E'_{r-1} \ar[r] \ar[d] &
\G'_{r}[x_{r-1}] \ar[r] \ar[d] & 0
\\
0\ar[r] & 
\Hc(-x_{r-1}) \ar[r] \ar[d] &
\G'_{r-1} \ar[r] \ar[d] &
Q \ar[r] \ar[d] & 0
\\
 & ? & 0 & 0 & 
}
$$
where $Q:=\coker(\F'_{r-1}[x_{r-1}]\to\G'_{r}[x_{r-1}])$ and
$\Hc:=\ker(\G'_{r-1}\to Q)(x_{r-1})$.
I claim that this diagram remains exact if we replace $''?''$ by
zero. Indeed, we have
$$
\ker(\F'_{r-1}[x_{r-1}]\to\G'_{r}[x_{r-1}])=
\F'_{r-1}[x_{r-1}]\cap\F'_{r}[x_{r-1}] =
H^0(R_{r-1}\cup R_r\ , \ \E'|_{R_{r-1}\cup R_r}(-x_{r-2}-x_r))
$$
which is zero, since $\E'$ and therefore also 
$\E'|_{R_{r-1}\cup R_r}$ is admissible.
By definition of $\Hc$ it is clear that it is strictly standard
of degree $\rk(\G'_{r-1})-\dim(Q)=\delta'_r$.
Being an extension of $\F'_{r-1}$ and $\Hc$, the $\Oo_{R_{r-1}}$-module
$\E_{r-1}$ is therefore strictly standard of degree 
$\delta_{r-1}=\delta'_r+\delta'_{r-1}$.
Furthermore, we have
$$
H^0(R,\E|_R(-x_0-x_{r-1}))=H^0(R,f_*\E'|_{R'}(-x_0-x_r))=(0)
\quad.
$$
This proves the admissibility of $\E$.
\end{proof}

\section{Gieseker vector bundle data and generalized isomorphisms}
\label{section GVBD and KGl}

Let $\VB(\tilde{C_0})$ be the moduli stack of vector bundles of
rank $n$ over $\tilde{C_0}$ (cf. \cite{LM} (4.6.2.1)).
Let $\E_{\univ}$ be the universal bundle over 
$\tilde{C_0}\times\VB(\tilde{C_0})$. By abuse of notation, we
denote by $p_i$ also the section of the projection
$\tilde{C_0}\times\VB(\tilde{C_0})\to\VB(\tilde{C_0})$
which is induced by the point $p_i\in \tilde{C_0}$ (i=1,2).
The purpose of this section is to show that the stack
$\GVBD(\tilde{C_0},p_1,p_2)$ is canonically isomorphic
to the stack $\KGln(p_1^*\E_{\univ},p_2^*\E_{\univ})$.
By definition, an $S$-valued point of 
$\KGln(p_1^*\E_{\univ},p_2^*\E_{\univ})$
is a vector bundle $\E$ over $C_0\times S$ together with a 
generalized isomorphism between the $\Oo_S$-modules 
$\E|_{\{p_1\}\times S}$ and $\E|_{\{p_2\}\times S}$ (cf. \cite{Kausz}).

\begin{construction}
\label{GVBD to VB}
Let $S$ be a $k$-scheme and let $(\Cc,\pi,s_1,s_2,h;\E,\varphi)$
an object in 
\linebreak[4]
$\GVBD_n(\tilde{C_0},p_1,p_2)(S)$ (cf. definition \ref{GVBD}).
We want to construct a vector bundle $h_{\bullet}(\E)$ of rank $n$ on 
$\tilde{C_0}\times_{\Spec(k)}S$ out of this data.

By setting 
$$
(\Cc^{(n)},\pi^{(n)},s_1^{(n)},s_2^{(n)},h^{(n)}):=
(\Cc,\pi,s_1,s_2,h)
$$
and $\E^{(n)}:=\E$
and
repeatedly applying \ref{contraction}, we obtain a sequence 
$$
(\Cc^{(i)},\pi^{(i)},s_1^{(i)},s_2^{(i)},h^{(i)})_{i=0,\dots,n}
$$
of modifications 
of $(\tilde{C_0},p_1,p_2)$ over $S$,
which are linked by morphisms as follows
$$
\Cc=\Cc^{(n)}\overset{f^{(n)}}{\To}
\Cc^{(n-1)}\overset{f^{(n-1)}}{\To}
\dots
\overset{f^{(1)}}{\To}\Cc^{(0)}
\quad,
$$
and for each $i\in[0,n]$ a
locally free $\Oo_{\Cc^{(i)}}$-module 
$\E^{(i)}$ of rank $n$, which is admissible for
$(\Cc^{(i)},\pi^{(i)},s_1^{(i)},s_2^{(i)},h^{(i)})$
and which is of extremal degree $\geq(n-i+1,1)$.

Now we set
$$
(\Cc^{[n]},\pi^{[n]},s_1^{[n]},s_2^{[n]},h^{[n]}):=
(\Cc^{(0)},\pi^{(0)},s_1^{(0)},s_2^{(0)},h^{(0)})
$$
and $\E^{[n]}:=\E^{(0)}$ and again successively apply
\ref{contraction} with the role of $s_1$ and $s_2$ interchanged,
thus obtaining a sequence
$$
(\Cc^{[i]},\pi^{[i]},s_1^{[i]},s_2^{[i]},h^{[i]})_{i=0,\dots,n}
$$
of modifications 
of $(\tilde{C_0},p_1,p_2)$ over $S$, a chain
$$
\Cc^{(0)}=\Cc^{[n]}\overset{f^{[n]}}{\To}
\Cc^{[n-1]}\overset{f^{[n-1]}}{\To}
\dots
\overset{f^{[1]}}{\To}\Cc^{[0]}
$$
of simple modifications, and for each $i\in[0,n]$
an $\Oo_{\Cc^{[i]}}$-module 
$\E^{[i]}$, which is admissible for
$(\Cc^{[i]},\pi^{[i]},s_1^{[i]},s_2^{[i]},h^{[i]})$
and of extremal degree $\geq(n+1,n-i+1)$.

In particular, $\E^{[0]}$ is of extremal degree $\geq(n+1,n+1)$
which is easily seen to imply that the morphism 
$h^{[0]}:\Cc^{[0]}\to\tilde{C_0}\times_{\Spec(k)}S$
is an isomorphism.
We set
$
h_{\bullet}(\E):=h^{[0]}_*(\E^{[0]})
$. 
\end{construction}

\begin{construction}
\label{GVBD to KGl}
Let $S$ be a $k$-scheme, $(\Cc,\pi,s_1,s_2,h;\E,\varphi)$
an object in 
\linebreak[4]
$\GVBD_n(\tilde{C_0},p_1,p_2)(S)$ (cf. definition \ref{GVBD})
and $h_{\bullet}(\E)$ the $\Oo_{\tilde{C_0}\times_{\Spec(k)}S}$-module
associated to this data (cf. \ref{GVBD to VB}).
We want to construct a generalized isomorphism 
$$
\Phi=
\left(
\xymatrix@C=1.3ex{
E
\ar@/^1.2pc/|{\tensor}[rr]
\ar @/_0.65pc/ @{} [rr]|{(M_0,\mu_0)}
& &
E_1 
\ar[ll]_0
\ar @/_0.65pc/ @{} [rr]|{(M_1,\mu_1)}
\ar@/^1.2pc/|{\tensor}[rr]
& &
E_2
\ar[ll]_1
& 
\dots
& 
E_{n-1}
\ar @/_0.65pc/ @{} [rr]|{(M_{n-1},\mu_{n-1})}
\ar@/^1.2pc/|{\tensor}[rr]
& & 
E_n
\ar[ll]_{n-1}
\ar[rr]^\sim
& & 
F_n
\ar[rr]^{n-1}
& & 
F_{n-1}
\ar@/_1.2pc/|{\tensor}[ll]
\ar @/^0.65pc/ @{} [ll]|{(L_{n-1},\lambda_{n-1})}
& 
\dots
&
F_2
\ar[rr]^1
& &
F_1
\ar[rr]^0
\ar @/^0.65pc/ @{} [ll]|{(L_1,\lambda_1)}
\ar@/_1.2pc/|{\tensor}[ll]
& &
F
\ar @/^0.65pc/ @{} [ll]|{(L_0,\lambda_0)}
\ar@/_1.2pc/|{\tensor}[ll]
}
\right)
$$
(cf. \cite{Kausz})
between the two
$\Oo_S$-modules $E:=p_1^*h_{\bullet}(\E)$ and $F:=p_2^*h_{\bullet}(\E)$,
where by abuse of notation we write 
$p_i:S\to\tilde{C_0}\times_{\Spec(k)}S$
for the section of $\tilde{C_0}\times_{\Spec(k)}S\to S$ induced
by the point $p_i\in\tilde{C_0}(k)$.

Let
$$
E_i:=(s_1^{(i)})^*\E^{(i)} 
\qquad , \qquad
F_i:=(s_2^{[i]})^*\E^{[i]} 
$$
for $i=0,\dots,n$
and let for $i=1,\dots,n$
$$
\xymatrix@C=2ex{
E_{i-1}
\ar@/^1.2pc/|{\tensor}[rr]
\ar @/_.8pc/ @{} [rr]|{(M_{i-1},\mu_{i-1})}
& & 
E_i
\ar[ll]_{i-1}
}
\qquad\text{resp.}\qquad
\xymatrix@C=2ex{
F_{i}
\ar[rr]^{i-1}
\ar @/_.8pc/ @{} [rr]|{(L_{i-1},\lambda_{i-1})}
& & 
F_{i-1}
\ar@/_1.2pc/|{\tensor}[ll]
}
$$
be the bf-morphism associated by construction \ref{constr3}
to the simple modification
$$
\xymatrix@R=1.5ex{
\Cc^{(i)} \ar[rr]^{f^{(i)}} \ar[dr]^{\pi^{(i)}} & & 
\Cc^{(i-1)} \ar[dl]_{\pi^{(i-1)}} \\
   & S \ar@{-<}@/^/[ul]^{s_1^{(i)}} \ar@/_/[ur]_{s_1^{(i-1)}} &
}
\qquad\text{resp.}\qquad
\xymatrix@R=1.5ex{
\Cc^{[i]} \ar[rr]^{f^{[i]}} \ar[dr]^{\pi^{[i]}} & & 
\Cc^{[i-1]} \ar[dl]_{\pi^{[i-1]}} \\
   & S \ar@/^/[ul]^{s_2^{[i]}} \ar@/_/[ur]_{s_2^{[i-1]}} &
}
$$
and the $\Oo_{\Cc^{(i)}}$-module $\E^{(i)}$
(resp. the $\Oo_{\Cc^{[i]}}$-module $\E^{[i]}$).
Observe that $E_0=E$, $F_0=F$, $E_n=s_1^*\E$ and $F_n=s_2^*\E$.
In particular, $\varphi$ defines an isomorphism $E_n\isomto F_n$
and thus we have constructed 
all the data $\Phi$, which make up a generalized isomorphism 
from $E$ to $F$.
\end{construction}

\begin{lemma}
\label{GVBD to KGl welldefined}
The data
$$
\Phi=
\left(
\xymatrix@C=1.3ex{
E
\ar@/^1.2pc/|{\tensor}[rr]
\ar @/_0.65pc/ @{} [rr]|{(M_0,\mu_0)}
& &
E_1 
\ar[ll]_0
\ar @/_0.65pc/ @{} [rr]|{(M_1,\mu_1)}
\ar@/^1.2pc/|{\tensor}[rr]
& &
E_2
\ar[ll]_1
& 
\dots
& 
E_{n-1}
\ar @/_0.65pc/ @{} [rr]|{(M_{n-1},\mu_{n-1})}
\ar@/^1.2pc/|{\tensor}[rr]
& & 
E_n
\ar[ll]_{n-1}
\ar[rr]^\sim
& & 
F_n
\ar[rr]^{n-1}
& & 
F_{n-1}
\ar@/_1.2pc/|{\tensor}[ll]
\ar @/^0.65pc/ @{} [ll]|{(L_{n-1},\lambda_{n-1})}
& 
\dots
&
F_2
\ar[rr]^1
& &
F_1
\ar[rr]^0
\ar @/^0.65pc/ @{} [ll]|{(L_1,\lambda_1)}
\ar@/_1.2pc/|{\tensor}[ll]
& &
F
\ar @/^0.65pc/ @{} [ll]|{(L_0,\lambda_0)}
\ar@/_1.2pc/|{\tensor}[ll]
}
\right)
$$
constructed in \ref{GVBD to KGl} is a generalized
isomorphism from $E$ to $F$. The association
$$
(\Cc,\pi,s_1,s_2,h;\E,\varphi)\mapsto (h_{\bullet}\E,\Phi)
$$
is functorial with respect to isomorphisms and commutes 
with base-change $S'\to S$,
and hence defines a morphism
$$
\GVBD(\tilde{C_0},p_1,p_2)\To
\KGl(p_1^*\E_{\univ},p_2^*\E_{\univ})
$$
of $k$-groupoids.
\end{lemma}

\begin{proof}
By construction it is clear that the association
$
(\Cc,\pi,s_1,s_2,h;\E,\varphi)\mapsto (h_{\bullet}\E,\Phi)
$
is functorial with respect to isomorphisms and commutes with base change.
All that remains to be shown is that $\Phi$ satisfies  
properties 1 and 2 in \cite{Kausz} definition 5.2. 
Since these are pointwise properties, we may assume that $S=\Spec(K)$ 
for some field $K$.
Consider the following situation:
$$
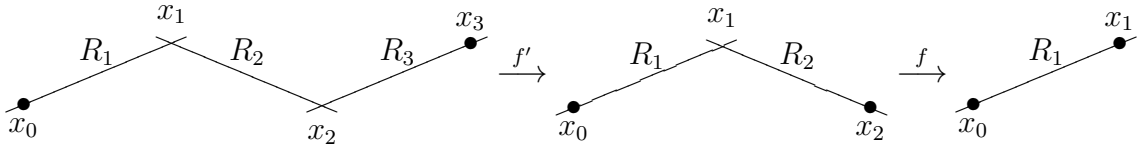
\vcenter{
\xy <1cm,0cm>:
(-0.2,0);
(2.2,1) **@{-} 
?<<< *{\bullet}
+(0,-0.3) *{x_0};
(2,1.3) *{x_1};
(1,0.8) *{R_1}
;
(1.8,1);
(4.2,0) **@{-};
(4,-0.3) *{x_2};
(3,0.8) *{R_2}
;
(3.8,0);
(6.2,1) **@{-}
?>>> *{\bullet}
+(0,0.3) *{x_3};
(5,0.8) *{R_3}
\endxy
}
\overset{f'}{\To}
\vcenter{
\xy <1cm,0cm>:
(-0.2,0);
(2.2,1) **@{-} 
?<<< *{\bullet}
+(0,-0.3) *{x_0};
(2,1.3) *{x_1};
(1,0.8) *{R_1}
;
(1.8,1);
(4.2,0) **@{-}
?>>> *{\bullet}
+(0,-0.3) *{x_2};
(3,0.8) *{R_2}
\endxy
}
\overset{f}{\To}
\vcenter{
\xy <1cm,0cm>:
(-0.2,0);
(2.2,1) **@{-} 
?<<< *{\bullet}
+(0,-0.3) *{x_0}
?>>> *{\bullet}
+(0,0.3) *{x_1};
(1,0.8) *{R_1}
\endxy
}
$$
where the $R_i$ are projective lines over $K$ and 
$f':R'':=\cup_{i=1}^3R_i\to R':=\cup_{i=1}^2R_i$
and
$f: R'\to R:=R_1$
are simple modifications over $K$ as indicated.
Let $\E''$ be an admissible 
$\Oo_{R''}$-module of rank $n$, let
$\E':=(f'_*\E''(-x_3))(x_2)$ and let 
$\E:=(f_*\E'(-x_2))(x_1)$.
Denote by $\E''_i$ (resp. $\E'_i$, $\E_i$) the restriction of $\E''$
(resp. of $\E'$, $\E$) to the component $R_i$.
We have canonical exact sequences 
$$
0\To\F''_i\To\E''_i\To\G''_i\To 0
$$
of $\Oo_{R_i}$-modules, where $\G''_i$ is trivial and $\F''_i$ is
a direct sum of copies of $\Oo_{R_i}(1)$. Analogously let
$\F'_i$ and $\G'_i$  be the canonical
subsheaf and quotient sheaf of $\E'_i$.
By \ref{constr3 over K} we have canonical diagrams
$$
\xymatrix@R=1.5ex{
\E[x_1] \ar[rr] \ar@{->>}[dr] &&
M\tensor\E'[x_2] \ar[rr] \ar@{->>}[dr] &&
M\tensor M'\tensor\E''[x_3] \\
& M\tensor\F'_2[x_2] \ar@{^(->}[ur] \ar[rr]^{u}  &&
M\tensor M'\tensor\F''_3[x_3] \ar@{^(->}[ur]  & 
\\
\E''[x_3] \ar[rr] \ar@{->>}[dr] &&
\E'[x_2] \ar[rr] \ar@{->>}[dr] &&
\E[x_1] \\
& \G''_3[x_2] \ar@{^(->}[ur] \ar[rr]^{v}  &&
\G'_2[x_1] \ar@{^(->}[ur]  & 
}
$$
where $M:=(\m_2/\m_2^2)\tensor(\m_1/\m_1^2)^{\dual}$,
$M':=(\m_3/\m_3^2)\tensor(\m_2/\m_2^2)^{\dual}$
and $\m_i$ is the maximal ideal of $\Oo_{R_i,x_i}$.
To show that $\Phi$ satisfies property 1 in \cite{Kausz}
definition 5.2 it is cleary sufficient to show that 
the morphisms $u$ and $v$ in these diagrams are surjective.
In fact by the exact diagram
$$
\xymatrix@R=1ex{
& & 0\ar[d] & & \\
& & \F'_2[x_2] \ar[d] \ar[dr]^u & & \\
0\ar[r] & \G''_3[x_2]\ar[r]\ar[dr]_v & 
\E'[x_2]\ar[r]\ar[d] & 
M'\tensor \F''_3[x_3]\ar[r] & 0 \\
& & \G'_2[x_1]\ar[d] & & \\
& & 0 & &
}
$$
it suffices to prove surjectivity of either $u$ or $v$.
From the canonical exact sequence of $\Oo_{R_2}$-modules 
(cf. \ref{constr3 over K}):
$$
0\To\E'_2(-x_2)\To\E''_2\To\G''_3[x_2]\To 0
$$
we deduce the exact diagram of $K$-vector spaces:
$$
\xymatrix@R=1.5ex{
0\ar[r] &
H^0(\E'_2(-x_2)) \ar[r] \ar[d]_{u'} &
H^0(\E''_2) \ar[r] \ar[d] &
\G''_3[x_2] \ar[r] \ar@{=}[d] &
0 \\
0\ar[r] &
\F''_3[x_2] \ar[r] \ar[d] &
\E''_2[x_2]\ar[r] \ar[d] &
\G''_3[x_2] \ar[r] &
0 \\
& 0 & 0 & & 
}
$$
Observe that we have canonical isomorphisms
$$
H^0(\E'_2(-x_2))=H^0(\F'_2(-x_2))=(\m_2/\m_2^2)\tensor\F'_2[x_2]
\quad\text{and}\quad
\F''_3[x_2]=(\m_3/\m_3^2)\tensor\F''_3[x_3]
$$
and that the morphism 
$$
(\m_2/\m_2^2)\tensor\F'_2[x_2]=H^0(\E'_2(-x_2))
\overset{u'}{\To}
\F''_3[x_2]=(\m_3/\m_3^2)\tensor\F''_3[x_3]
$$
coincides with the morphism $u$. This proves the surjectivity
of $u$ (and of $v$).

Now we want to show that $\Phi$ satisfies property 2.
in \cite{Kausz} definition 5.2., i.e. that 
$\varphi(\ker(E_0\ot E_n))\cap \ker(F_n\to F_0)=(0)$.
Let $R$ (resp. $R'$) be the chain of projective lines in $\Cc$,
which by $h$ is contracted to the point $p_1$ (resp. to $p_2$)
and let $r$ (resp. $r'$) be its length.
If $r\geq 1$ (resp. $r'\geq 1$), we denote by $V_{r}\subseteq\E[s_1]$ 
(resp. by $V'_{r'}\subseteq\E[s_2]$)
the image of $H^0(\E|_{R}(-p_1))$ in $\E[s_1]=E_n$
(resp. the image of $H^0(\E|_{R'}(-p_2))$ in $\E[s_2]=F_n$).
If $r=0$ (resp. $r'=0$), we set $V_{r}:=(0)$ (resp. $V'_{r'}=(0)$).
Admissibility of $\E$ implies that 
$\varphi(V_{r})\cap V'_{r'}=(0)$.
Therefore (and by reasons of symmetry) it suffices to show the equality 
$\ker(E_0\ot E_n)=V_{r}$. 
We will prove this by induction on $r$. The case $r=0$ is trivial.
If $r=1$, we have $\ker(E_0\ot E_n)=\F[s_1]=V_1$, where 
$\F\injto\E|_R$ is the canonical subsheaf of $\E|_R$.
Let $r\geq 2$ and let $f:R\to\Rb$ be the simple modification
of $(R,s_1)$ which contracts the component $R_r$ of $R$
to a point $x_{r-1}\in\Rb$. As usual, we identify $\Rb$ 
with a subchain of $R$. 
Let $\G$ be the canonical quotient
sheaf of $\E|_{R_r}$ and let $\Eb:=(f_*\E|_{R}(-s_1))(x_{r-1})$.
By \ref{constr3 over K} the morphism $E_n\to E_0$ factorizes
as follows: 
$
E_n=\E[s_1] \to \G[x_{r-1}] \to \Eb[x_{r-1}] \to E_0
$.
Let $\Vb_{r-1}\subseteq \Eb[x_{r-1}]$ 
(resp. $V_{r-1}\subseteq \E[x_{r-1}]$)
be the image of 
the map $H^0(\Eb|_{\Rb}(-p_1))\to\Eb[x_{r-1}]$
(resp. of the map $H^0(\E|_{\Rb}(-p_1))\to\E[x_{r-1}]$).
By induction hypothesis, we have 
$\Vb_{r-1}=\ker(\Eb[x_{r-1}]\to E_0)$ and by
lemma \ref{equivalents for admissible} the subspace
$V_r$ is the preimage under 
$
\E[s_1] \to \G[x_{r-1}]
$
of the image of $V_{r-1}$ by 
$
\E[x_{r-1}] \to \G[x_{r-1}]
$.
Therefore we have to show the following equality
of subspaces of $\G[x_{r-1}]$:
$$
\Ub:=
\left(
\begin{array}{ll}
\text{preimage of $\Vb_{r-1}$ under the}\\
\text{morphism $\G[x_{r-1}]\to\Eb[x_{r-1}]$}
\end{array}
\right)
=
\left(
\begin{array}{ll}
\text{image of $V_{r-1}$ by the}\\
\text{morphism $\E[x_{r-1}]\to\G[x_{r-1}]$}
\end{array}
\right)
=:U
\quad.
$$
Consider the exact diagram of $\Oo_{\Rb}$-modules (cf. \ref{constr3 over K}):
$$
\xymatrix@R=2ex{
& 0 \ar[d] & 0 \ar[d] & \G[x_{r-1}] \ar@{=}[d]& 
\\
0 \ar[r] &
\Eb(-p_1-x_{r-1}) \ar[d] \ar[r] &
\E|_{\Rb}(-p_1) \ar[d] \ar[r] &
\G[x_{r-1}] \ar[d]^0 \ar[r] & 0
\\
0 \ar[r] &
\Eb(-p_1) \ar[d] \ar[r] &
\E|_{\Rb}(-p_1+x_{r-1}) \ar[d] \ar[r] &
(\m/\m^2)^{\dual}\tensor\G[x_{r-1}] \ar@{=}[d] \ar[r] & 0
\\
 &
\Eb[x_{r-1}] \ar[d] \ar[r] &
(\m/\m^2)^{\dual}\tensor\E[x_{r-1}] \ar[d] \ar[r] &
(\m/\m^2)^{\dual}\tensor\G[x_{r-1}]  \ar[r] & 0
\\
& 0 & 0 & &
}
$$
Application of $H^0=H^0(\Rb,\cdot)$ leads to the exact diagram
of $K$-vector spaces
$$
\xymatrix@R=2ex{
& & 0 \ar[d] & 0 \ar[d] & U \ar@{=}[d]& 
\\
& 0 \ar[r] &
H^0(\Eb(-p_1-x_{r-1})) \ar[d] \ar[r] &
H^0(\E|_{\Rb}(-p_1)) \ar[d] \ar[r] &
U \ar[d]^0 \ar[r] & 0
\\
& 0 \ar[r] &
H^0(\Eb(-p_1)) \ar[d] \ar[r] &
H^0(\E|_{\Rb}(-p_1+x_{r-1})) \ar[d] \ar[r] &
\bullet \ar@{=}[d] \ar[r] & 0
\\
0 \ar[r] & 
\Ub \ar[r] &
V_{r-1} \ar[d] \ar[r] &
\bullet \ar[d] \ar[r] &
\bullet  \ar[r] & 0
\\
& & 0 & 0 & &
}
$$
which proves $\Ub=U$ as required.
\end{proof}

\begin{remark}
In \ref{GVBD to VB}, \ref{GVBD to KGl} and \ref{GVBD to KGl welldefined}   
we have  shown that starting with a Gieseker vector bundle data, one can
produce a sequence of simple modifications together with admissible
bundles on them, and that the corresponding (via \ref{mod&E' to bf})
bf-morphisms make up a generalized isomorphism.        
The following picture may help to keep track of which simple modification
corresponds to which bf-morphism:
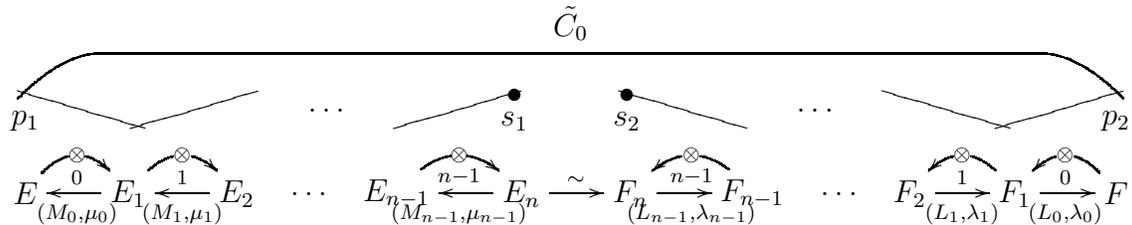
\samepage{
$$
\xy
0;<1cm,0cm>:
(-0.1,-0.1);(1,0.5) **\crv{(0.5,0.5)};
(1,0.5); (13.5,0.5) **@{-};                  
(13.5,0.5);(14.6,-0.1) **\crv{(14,0.5)};
(-0.1,0);(1.6,-0.5) **@{-}; 
(1.4,-0.5);(3.1,0) **@{-};
(4,-0.25) *{\dots};
(4.9,-0.5);(6.6,0) **@{-};              
(6.5,-0.06) *{\bullet};
(8,-0.06) *{\bullet};
(7.9,0);(9.6,-0.5) **@{-};
(10.5,-0.25) *{\dots};
(11.4,0);(13.1,-0.5) **@{-};
(12.9,-0.5);(14.6,0) **@{-};
(7.25,0.9 ) *{\text{$\tilde{C_0}$}};
(0,-0.4) *{p_1};
(14.5,-0.4) *{p_2};                    
(6.5,-0.4) *{s_1};
(8,-0.4) *{s_2};
\endxy
$$
$$
\xymatrix@C=1.3ex{
E
\ar@/^1.2pc/|{\tensor}[rr]
\ar @/_0.65pc/ @{} [rr]|{(M_0,\mu_0)}
& &
E_1 
\ar[ll]_0
\ar @/_0.65pc/ @{} [rr]|{(M_1,\mu_1)}
\ar@/^1.2pc/|{\tensor}[rr]
& &
E_2
\ar[ll]_1
& 
\dots
& 
E_{n-1}
\ar @/_0.65pc/ @{} [rr]|{(M_{n-1},\mu_{n-1})}
\ar@/^1.2pc/|{\tensor}[rr]
& & 
E_n
\ar[ll]_{n-1}
\ar[rr]^\sim
& & 
F_n
\ar[rr]^{n-1}
& & 
F_{n-1}
\ar@/_1.2pc/|{\tensor}[ll]
\ar @/^0.65pc/ @{} [ll]|{(L_{n-1},\lambda_{n-1})}
& 
\dots
&
F_2
\ar[rr]^1
& &
F_1
\ar[rr]^0
\ar @/^0.65pc/ @{} [ll]|{(L_1,\lambda_1)}
\ar@/_1.2pc/|{\tensor}[ll]
& &
F
\ar @/^0.65pc/ @{} [ll]|{(L_0,\lambda_0)}
\ar@/_1.2pc/|{\tensor}[ll]
}
$$
}

\vspace{3mm}
But one should keep in mind that a projective line corresponds only
to those bf-morphisms, where the section $\mu_i$ (resp. $\lambda_i$)
vanishes! In this sense the picture is slightly misleading.
\end{remark}

\begin{theorem}
\label{GVBD isomto KGl}
The morphism 
$$
\GVBD(\tilde{C_0},p_1,p_2)\To
\KGl(p_1^*\E_{\univ},p_2^*\E_{\univ})
$$
defined by construction \ref{GVBD to KGl} is an isomorphism 
of algebraic $k$-stacks.
\end{theorem}

The proof of theorem \ref{GVBD isomto KGl} will be given after
lemma \ref{KGl to GVBD welldefined} below.

\begin{construction}
\label{KGl to GVBD}
Let $S$ be a $k$-scheme,
let $\F$ be a locally free $\Oo_{\tilde{C_0}\times_{\Spec(k)}S}$-module
of rank $n$ and
let 
$$
\Phi=
\left(
\xymatrix@C=1.3ex{
E
\ar@/^1.2pc/|{\tensor}[rr]
\ar @/_0.65pc/ @{} [rr]|{(M_0,\mu_0)}
& &
E_1 
\ar[ll]_0
\ar @/_0.65pc/ @{} [rr]|{(M_1,\mu_1)}
\ar@/^1.2pc/|{\tensor}[rr]
& &
E_2
\ar[ll]_1
& 
\dots
& 
E_{n-1}
\ar @/_0.65pc/ @{} [rr]|{(M_{n-1},\mu_{n-1})}
\ar@/^1.2pc/|{\tensor}[rr]
& & 
E_n
\ar[ll]_{n-1}
\ar[rr]^\sim
& & 
F_n
\ar[rr]^{n-1}
& & 
F_{n-1}
\ar@/_1.2pc/|{\tensor}[ll]
\ar @/^0.65pc/ @{} [ll]|{(L_{n-1},\lambda_{n-1})}
& 
\dots
&
F_2
\ar[rr]^1
& &
F_1
\ar[rr]^0
\ar @/^0.65pc/ @{} [ll]|{(L_1,\lambda_1)}
\ar@/_1.2pc/|{\tensor}[ll]
& &
F
\ar @/^0.65pc/ @{} [ll]|{(L_0,\lambda_0)}
\ar@/_1.2pc/|{\tensor}[ll]
}
\right)
$$
be a generalized isomorphism from $E=E_0:=p_1^*\F$ to $F=F_0:=p_2^*\F$.
We want to construct a Gieseker vector bundle data
$
(\Cc,\pi,s_1,s_2,h;\E,\varphi)
$
on $(\tilde{C_0},p_1,p_2)$ over $S$ (cf. \ref{GVBD}).

Let 
$$
(\Cc^{[0]},\pi^{[0]},s_1^{[0]},s_2^{[0]},h^{[0]}):=
(\tilde{C_0}\times_{\Spec(k)}S,\pr_2,p_1,p_2,\id)
$$
and $\F^{[0]}:=\F$,
and for $i=1,\dots,n$ define inductively a modification
$
(\Cc^{[i]},\pi^{[i]},s_1^{[i]},s_2^{[i]},h^{[i]})
$
of
$
(\tilde{C_0},p_1,p_2)
$
over $S$
together with a locally free $\Oo_{\Cc^{[i]}}$-module
$\F^{[i]}$ and an isomorphism $(s_2^{[i]})^*\F^{[i]}\isomto F_i$
as follows:

Assume that 
$
(\Cc^{[i-1]},\pi^{[i-1]},s_1^{[i-1]},s_2^{[i-1]},h^{[i-1]})
$
together with 
$\F^{[i-1]}$ and an isomorphism $(s_2^{[i-1]})^*\F^{[i-1]}\isomto F_{i-1}$
has already been defined. 
By proposition \ref{mod&E' to bf}, the bf-morphism 
$$               
\xymatrix@C=2ex{
F_{i}
\ar[rr]^{i-1}
\ar @/_.8pc/ @{} [rr]|{(L_{i-1},\lambda_{i-1})}
& & 
F_{i-1}
\ar@/_1.2pc/|{\tensor}[ll]
}
$$
induces a canonical simple modification 
$
(\Cc^{[i]},f^{[i]},\pi^{[i]},s_2^{[i]})
$
of 
$
(\Cc^{[i-1]},\pi^{[i-1]},s_2^{[i-1]})
$
together with an $\Oo_{ \Cc^{[i]}}$-module $\F^{[i]}$
and an isomorphism $(s_2^{[i]})^*\F^{[i]}\isomto F_i$.
It is clear that there exists a unique section 
$s_1^{[i]}$ of $\pi^{[i]}$ such that 
$s_1^{[i-1]}=f^{[i]}\comp s_1^{[i]}$. Let 
$h^{[i]}:=h^{[i-1]}\comp f^{[i]}$.

Analogously, we define for $i=0,\dots,n$ a modification
$
(\Cc^{(i)},\pi^{(i)},s_1^{(i)},s_2^{(i)},h^{(i)})
$
of
$
(\tilde{C_0},p_1,p_2)
$
over $S$
together with a locally free $\Oo_{\Cc^{(i)}}$-module
$\F^{(i)}$ and an isomorphism $(s_1^{(i)})^*\F^{(i)}\isomto E_i$
by setting
$$
(\Cc^{(0)},\pi^{(0)},s_1^{(0)},s_2^{(0)},h^{(0)}):=
(\Cc^{[n]},\pi^{[n]},s_1^{[n]},s_2^{[n]},h^{[n]})
$$
and repeatedly applying proposition \ref{mod&E' to bf} with
the role of $s_1$ and $s_2$ interchanged where we now use
the bf-morphism 
$(M_{i-1},\mu_{i-1},E_i\to E_{i-1},M_{i-1}\tensor E_i\ot E_{i-1},i-1)$
instead of
$(L_{i-1},\lambda_{i-1},F_i\to F_{i-1},L_{i-1}\tensor F_i\ot F_{i-1},i-1)$
for the induction step.

Finally, let 
\begin{eqnarray*}
(\Cc,\pi,s_1,s_2,h)&:=&
(\Cc^{(n)},\pi^{(n)},s_1^{(n)},s_2^{(n)},h^{(n)})
\\
\E &:=& \F^{(n)}
\\
\varphi &:& s_1^*\E\isomto E_n\isomto F_n\isomto s_2^*\E
\end{eqnarray*}
This completes the construction.
\end{construction}

\begin{lemma}
\label{KGl to GVBD welldefined}
The tupel $(\Cc,\pi,s_1,s_2,h;\E,\varphi)$
constructed in \ref{KGl to GVBD} is a 
Gieseker vector bundle data
on $(\tilde{C_0},p_1,p_2)$ over $S$.
\end{lemma}

\begin{proof}
It is clear by construction that the diagram
$$
\xymatrix{
\Cc \ar@{->}[dr]^{\pi} \ar@{->}[rr]^(.4){h}    &   & 
\text{$\tilde{C_0}\times_{\Spec(k)}S$} \ar[dl]_{\pr_2}  \\
& S \ar@/^/[ul] \ar@<1ex>@/^/[ul]^{s_1,s_2} 
   \ar@/_/[ur] \ar@<-1ex>@/_/[ur]_{p_2\times\id,p_1\times\id} &   
}
$$
is commutative, that $\pi$ is flat, $h$ is proper and finitely
presented and that the $s_i$ are sections of $\pi$.
Since furthermore the construction \ref{KGl to GVBD} commutes
with base-change, we may assume that $S=Spec(K)$ for some field $K$.

Let $X$ be a curve over $K$, $\pi:X\to\Spec(K)$ the structure morphism
and $x_0\in X(K)$. 
Let $\G$ be a locally free $\Oo_X$-module of rank $n$ and let 
$G:=\G[x_0]$.
Let
$$
\xymatrix@C=4ex{
G
\ar@/^1.2pc/|{\tensor}[rr]
& & 
G' 
\ar[ll]_{n-d_1-d_2}^{(M,\mu)}
\ar@/^1.2pc/|{\tensor}[rr]
& & 
G''
\ar[ll]_{n-d_2}^{(M',\mu')}
}
$$
be a sequence of bf-morphisms of $K$-vectorspaces 
with $\mu=\mu'=0$ and such that
$$
\im(G\to\M\tensor G'\to\M\tensor\M'\tensor G'')=
\im(\M\tensor G'\to\M\tensor\M'\tensor G'')
$$
and
$$
\im(G\ot G')=\im(G\ot G'\ot G'')
\quad.
$$
Using proposition \ref{mod&E' to bf} we get a simple modification
$(X',f,\pi',x_1)$ of $(X,\pi,x_0)$ and an $\Oo_{X'}$-module
$\G'$ together with an isomorphism $\G'[x'_1]\isomto G'$.
A second application of \ref{mod&E' to bf} yields  a simple
modification $(X'',f',\pi'',x_2)$ of $(X',\pi',x_1)$ and an
$\Oo_{X''}$-module $\G''$. 
Thus we have the following situation:
$$
\xymatrix@C=2.7ex@R=1.5ex{
\vcenter{
\xy <1cm,0cm>:
(0.2,0); (2.2,1) 
**\crv{(0.5,0.6)&(1.5,0.4)} 
?>>> *{\bullet}
+(0,0.3) *{x_0};
(1,0.8) *{X}
\endxy
}
& & \ar[l]_{f} & 
\vcenter{
\xy <1cm,0cm>:
(0.2,0);
(2.2,1) **\crv{(0.5,0.6)&(1.5,0.4)}; 
(2,1.3) *{x_0};
(1,0.8) *{X}
;
(1.8,1);
(4.2,0) **@{-}
?>>> *{\bullet}
+(0,-0.3) *{x_1};
(3,0.8) *{R_1}
\endxy
}
& & \ar[l]_{f'} &
\vcenter{
\xy <1cm,0cm>:
(0.2,0);
(2.2,1) **\crv{(0.5,0.6)&(1.5,0.4)}; 
(2,1.3) *{x_0};
(1,0.8) *{X}
;
(1.8,1);
(4.2,0) **@{-};
(4,-0.3) *{x_1};
(3,0.8) *{R_1}
;
(3.8,0);
(6.2,1) **@{-}
?>>> *{\bullet}
+(0,0.3) *{x_2};
(5,0.8) *{R_2}
\endxy
}
\\
 &&& &&&
\\
\ar@{=}[u] &&& \ar@{=}[u] &&& \ar@{=}[u]
\\
X 
& & &
X' \ar[lll]_{f} 
& & &
X'' \ar[lll]_{f'} 
}
$$
and it is easy to see (by going through construction \ref{constr4}
for this case) that the restriction of $\G''$ to the chain of
rational curves $R_1\cup R_2$ is strictly standard.
This consideration 
shows that the restriction of $\E$ to the  
chains of rational curves $R':=h^{-1}(p_1)$  and $R'':=h^{-1}(p_2)$ 
is strictly standard. 

To prove the admissibility of $\E|_{R'}$, we have to show that
$H^0(R',\E|_{R'}(-s_1-p_1))=0$ (cf. \ref{equivalents for admissible}).
For this let $R'_i:=(h^{(i)})^{-1}(p_1)$ for $i=1,\dots,n$ and let 
$$
\xymatrix{
R'=R'_n \ar[r]^{f^{(n)}} &
R'_{n-1} \ar[r]^{f^{(n-1)}} &
\cdot & \cdot & \cdot & \ar[r]^{f^{(2)}} &
R'_1
}
$$
be the successive simple modifications intervening in construction
\ref{KGl to GVBD}. If all $R'_i$ are reduced to a point, there is
nothing to show. Otherwise there exists an $m\in[1,n]$ such that
$R'_m\isomorph\Pp^1$.
By construction we have
$
f^{(i)}_*(\F^{(i)}(-s_1^{(i)}))
=\F^{(i-1)}(-s_1^{(i-1)})
$
for $i=1,\dots,n$.
It follows that
$$
H^0(R',\E|_{R'}(-s_1-p_1)) =
H^0(R'_{1},\F^{(m)}|_{R'_{m}}(-s^{(m)}_1-p_1))
\quad,
$$
which is zero, since $\F^{(m)}|_{R'_{m}}$ is strictly standard.
Analogously, one shows the admissibility of $\E|_{R''}$.

Now let $R$ be the chain of rational curves which is induced from
$R'$ and $R''$ by identifying the points $s_1$ and $s_2$ and
let $\E_R$ be the $\Oo_R$-module induced by $\E|_{R'}$, $\E|_{R''}$
and $\varphi:\E|_{R'}[s_1]\isomto\E|_{R''}[s_2]$.
We have to show that $\E_R$ is admissible. This is clear, if
either $R'$ or $R''$ is a chain of length zero, so we may assume
that both are of length at least one.
In the last part of the proof of lemma \ref{GVBD to KGl welldefined}
we have seen that the image of $H^0(R',\E|_{R'}(-p_1))$
(resp. of $H^0(R'',\E|_{R''}(-p_2))$
in $\E|_{R'}[s_1]$ (resp. in $\E|_{R''}[s_2]$) is
$\ker(E_0\ot E_n)$ (resp. $\ker(F_n\to F_0)$).
Since by definition of a generalized isomorphism we have
$\varphi(\ker(E_0\ot E_n))\cap\ker(F_n\to F_0))=0$, it follows
that $H^0(R,\E_R(-p_1-p_2))=0$ as required.
\end{proof}

\begin{proof}(of theorem \ref{GVBD isomto KGl})
By lemma \ref{KGl to GVBD welldefined} construction \ref{KGl to GVBD}
is welldefined. By proposition \ref{mod&E' to bf}    
it is clear that 
construction \ref{KGl to GVBD} is inverse to construction
\ref{GVBD to KGl}.
\end{proof}

\section{Relationship with the stack of torsion free sheaves}
\label{section tfs}

Let $\TFS(C/B)$ be the algebraic $B$-stack which parametrizes
relatively
torsion-free sheaves of rank $n$ on $C$ over $B$ (cf. \cite{Faltings})
and let $\TFS(C_0/B_0)$ be its special fibre (the $B_0$-stack
which parametrizes torsion-free sheaves of rank $n$ on $C_0$).
The stack $\TFS(C/B)$ 
is known to be singular; its local structure has been studied in
\cite{Faltings}.
If $S$ is a $B$-scheme and $(h:\Cc\to C\times_BS\ ,\ \E)$ a Gieseker
vector bundle on $C$ over $S$, then the push-forward $h_*\E$ is
relatively torsion-free (cf. \ref{torsion-free}). Therefore we obtain
a morphism $\GVB(C/B)\to\TFS(C/B)$ which specializes to 
$\GVB(C_0/B_0)\to\TFS(C_0/B_0)$.
Denoting, as before, by $\E_{\univ}$ the universal vector bundle
on $\tilde{C_0}\times\VB(\tilde{C_0}/B_0)$,
there is also a morphism 
$
\Grass_n(p_1^*\E_{\univ}\oplus p_2^*\E_{\univ})
\to
\TFS(C_0/B_0)
$
which is defined as follows (cf. section 6 of \cite{Seshadri2}):
Let $S$ be a $k$-scheme, $\E$ a locally free 
$\Oo_{\tilde{C_0}\times S}$-module of rank $n$ and
$
\xymatrix{
p_1^*\E\oplus p_2^*\E \ar@{->>}[r] & Q
}
$ 
an epimorphism of $\Oo_S$-modules,
where $Q$ is locally free of rank $n$. To these data
we associate a relatively torsion-free 
$\Oo_{C_0\times S}$-module $\F$ by setting
$$
\F:= \ker(r_*\E \to p_*(p_1^*\E\oplus p_2^*\E)\to p_*Q)
\quad,
$$
where $r:\tilde{C_0}\times S\to C_0\times S$ is induced by the
normalization morphism and $p:S\to C_0\times S$ is the section
induced by the singular point $p\in C_0(k)$.
Finally, by section 10 of \cite{Kausz} 
there is a canonical morphism 
$
\KGln(p_1^*\E_{\univ},p_2^*\E_{\univ})\to
\Grass_n(p_1^*\E_{\univ}\oplus p_2^*\E_{\univ})
$.

\begin{proposition}
\label{GVB and TFS}
The following diagram is commutative:
$$
\xymatrix{
\KGl(p_1^*\E_{\univ},p_2^*\E_{\univ})
\ar[r]_(.55){\isomorph}^(.55){\text{cf. \ref{GVBD isomto KGl}}} \ar[d] 
&
\text{$\GVBD(\tilde{C_0},p_1,p_2)$}
\ar[r]^(.53){\text{cf. \ref{normalization}}}
&
\GVB(C_0/B_0)
\ar[d]
\\
\Grass_n(p_1^*\E_{\univ}\oplus p_2^*\E_{\univ})
\ar[rr]
& &
\TFS(C_0/B_0)
}
$$
\end{proposition}

\begin{proof}
Let $S$ be a $k$-scheme.
Let
$$
\xymatrix{
\text{$\tilde{C}$} \ar[r]^{f} \ar[d]_{\text{$\tilde{\pi}$}}
& 
\text{$\tilde{C_0}\times S$}
\\
S \ar@/_/[u] \ar@<-1ex>@/_/[u]_{s_1,s_2}         
&
}
$$
be a modification of $(\tilde{C_0},p_1,p_2)$ over $S$ and
let $(\tilde{\E},\varphi:s_1^*\tilde{\E}\isomto s_2^*\tilde{\E})$
be a locally free $\Oo_{\tilde{\Cc}}$-module of rank $n$ which is amissible
for $(\tilde{\Cc},\tilde{\pi},s_1,s_2,f)$.
Recall the definition of $f_{\bullet}\E$ from construction \ref{GVBD to VB}.
We have $f_{\bullet}\E=(f_*\tilde{\E}(-s_1-s_2))(p_1+p_2)$.
Let $(h:\Cc\to C_0\times S,\E)$ be the Gieseker vector bundle on $C_0$
over $S$ associated to the Gieseker vektor bundle data 
$(\tilde{\Cc},\tilde{\pi},s_1,s_2,f,\E,\varphi)$.
Then we have a commutative diagram:
$$
\xymatrix{
\text{$\tilde{\Cc}$} \ar[r]^q \ar[d]_f  & \Cc \ar[d]^h \\
\text{$\tilde{C_0}\times S$} \ar[r]^r &  C_0\times S
}
$$
where $q$ is the clutching morphism which maps the 
sections $s_1$ and $s_2$ of $\tilde{\pi}$ onto a section 
$s$ of $\pi$, which meets $\Cc$ in the
singular locus of $\pi$.
Let $s_1+s_2: S\amalg S\to \tilde{\Cc}$ and 
$p_1+p_2: S\amalg S\to \tilde{C_0}\times S$ be the morphisms induced by
the sections $s_1, s_2$ and $p_1, p_2$.
Similarly as in \ref{constr3}, we have a commutative diagram with exact rows:
$$
\xymatrix{
0 \ar[r] &
\text{$f_*\tilde{\E}(-s_1-s_2)$} \ar[r] \ar[d] &
\text{$f_*\tilde{\E}$} \ar[r]  \ar[d] \ar@{-->}[dl] &
\text{$(p_1+p_2)_*(s_1+s_2)^*\tilde{\E}$} \ar[r] \ar[d]^0 &
0 \\
0 \ar[r] &
\text{$f_{\bullet}\tilde{\E}$} \ar[r] &
\text{$(f_*\tilde{\E})(p_1+p_2)$} \ar[r]   &
\text{$((p_1+p_2)_*(s_1+s_2)^*\tilde{\E})(p_1+p_2)$} \ar[r]  &
0 
}
$$
which induces the arrow $f_*\tilde{\E}\to f_{\bullet}\tilde{\E}$.
Thus we obtain the following diagram with exact rows:
$$
\vcenter{
\xymatrix{
0 \ar[r] &
\text{$f_*\tilde{\E}(-s_1-s_2)$} \ar[r] \ar@{=}[d] &
\text{$f_*\tilde{\E}$} \ar[r]  \ar[d]  &
\text{$(p_1+p_2)_*(s_1+s_2)^*\tilde{\E}$} \ar[r] \ar@{-->}[d] &
0 \\
0 \ar[r] &
\text{$f_{\bullet}\tilde{\E}(-p_1-p_2)$} \ar[r] &
\text{$f_\bullet\tilde{\E}$} \ar[r]   &
\text{$(p_1+p_2)_*(s_1+s_2)^*(f_\bullet\tilde{\E})$} \ar[r]  &
0 
}}
\eqno(*)
$$
and in particular a canonical morphism
$\overline{\alpha}: s_1^*\tilde{\E}\oplus s_2^*\tilde{\E}
\to p_1^*(f_\bullet\tilde{\E})\oplus p_2^*(f_\bullet\tilde{\E})$.
Let $Q$ be the cokernel of the composed morphism
$$
\xymatrix{
\text{$s_1^*\tilde{\E}$} \ar@{^{(}->}[r]^(.4){(\id,\varphi)} &
\text{$s_1^*\tilde{\E}\oplus s_2^*\tilde{\E}$} 
\ar[r]^(.4){\text{$\overline{\alpha}$}} &
\text{$p_1^*(f_\bullet\tilde{\E})\oplus p_2^*(f_\bullet\tilde{\E})$}
}
$$
Then the data consisting in $f_\bullet\tilde{\E}$ together
with the quotient map 
$
\xymatrix{
\text{$p_1^*(f_\bullet\tilde{\E})\oplus p_2^*(f_\bullet\tilde{\E})$}
\ar@{->>}[r] & 
Q
}
$
is the image under 
$
\GVBD(\tilde{C_0},p_1,p_2)\isomto
\KGl(p_1^*\E_{\univ},p_2^*\E_{\univ})\to
\Grass_n(p_1^*\E_{\univ}\oplus p_2^*\E_{\univ})
$
of the Gieseker vector bundle data
$(\tilde{\Cc},\tilde{\pi},s_1,s_2,f,\E,\varphi)$.
Therefore it suffices to show that there exists a canonical exact
sequence of $\Oo_{\tilde{C_0}\times S}$-modules as follows:
$$
\xymatrix{
0 \ar[r] &
h_*\E \ar[r] &
\text{$r_*(f_\bullet\tilde\E)$} \ar[r]^\beta &
p_*Q \ar[r] &
0
}
$$
where 
$\beta:r_*(f_\bullet\tilde\E)\to p_*Q$ is 
the composed map 
$
r_*(f_\bullet\tilde\E) \to
r_*(p_1+p_2)_*(p_1+ p_2)^*(f_\bullet\tilde{\E}) =
p_*(p_1^*(f_\bullet\tilde{\E})\oplus p_2^*(f_\bullet\tilde{\E})) \to
p_*Q
$.
For abbreviation we set
$E_n:=s_1^*\tilde{\E}$, $F_n:=s_2^*\tilde{\E}$,
$E_0:=p_1^*f_\bullet\tilde{\E}$, $F_0:=p_2^*f_\bullet\tilde{\E}$.
Applying $r_*$ to the diagram $(*)$ gives the exact diagram
$$
\vcenter{
\xymatrix{
0 \ar[r] &
\text{$r_*f_*\tilde{\E}(-s_1-s_2)$} \ar[r] \ar@{=}[d] &
\text{$r_*f_*\tilde{\E}$} \ar[r]  \ar[d]^\alpha  &
p_*(E_n\oplus F_n) \ar[r] \ar[d]^{\text{$\overline{\alpha}$}} &
0 \\
0 \ar[r] &
\text{$r_*f_{\bullet}\tilde{\E}(-p_1-p_2)$} \ar[r] &
\text{$r_*f_\bullet\tilde{\E}$} \ar[r]   &
p_*(E_0\oplus F_0) \ar[r]  &
0 
}}
\eqno(**)
$$
We also have the following exact diagram:
$$
\vcenter{
\xymatrix{
0 \ar[r] &
E_n \ar[rr]^{(\id,\varphi)} \ar@{=}[d] & &
E_n\oplus F_n \ar[r]^(.55){\varphi-\id}  \ar[d]^{\text{$\overline{\alpha}$}}  &
F_n \ar[r] \ar[d]^{\text{$\alpha'$}} &
0 \\
0 \ar[r] &
E_n \ar[rr]^{\text{$\overline{\alpha}\comp(\id,\varphi)$}} & &
E_0\oplus F_0 \ar[r]   &
Q \ar[r]  &
0 
}}
\eqno(***)
$$
where $\alpha'$ is the morphism induced by $\overline{\alpha}$.
From diagrams $(**)$ and $(***)$ it follows that the following
diagram commutes:
$$
\vcenter{
\xymatrix{
\text{$r_*f_*\tilde{\E}$} \ar[r] \ar[d]^\alpha &
p_*(E_n\oplus F_n) \ar[r]^(.6){\varphi-\id} 
                   \ar[d]^{\text{$\overline{\alpha}$}} &
p_*F_n \ar[d]^{\alpha'}
\\
\text{$r_*f_\bullet\tilde{\E}$} \ar[r] \ar@/_2ex/[rr]_\beta &
p_*(E_0\oplus F_0) \ar[r] &
P_*Q
}}
\eqno(\dagger)
$$
Therefore the right square in the following diagram 
$(\dagger\dagger)$ commutes:
$$
\vcenter{
\xymatrix{
0 \ar[r] &
h_*\E \ar[r] \ar@{.>}[d] & 
\text{$r_*f_*\tilde{\E}$} \ar[r] \ar[d]^\alpha  &
p_*F_n \ar[r] \ar[d]^{\text{$\alpha'$}} &
0 \\
0 \ar[r] &
\ker(\beta) \ar[r] &
\text{$r_*f_\bullet\tilde{\E}$} \ar[r]^\beta   &
Q \ar[r]  &
0 
}}
\eqno(\dagger\dagger)
$$
The upper row in $(\dagger\dagger)$ is exact, since it comes from the
exact sequence
$$
\xymatrix{
0 \ar[r] &
\E \ar[r] &
\text{$f_*\tilde{\E}$} \ar[r] &
s_*F_n \ar[r] &
0
}
$$
by applying the functor $h_*$.
Furthermore, it follows from $(**)$ and $(***)$ that in  diagram $(\dagger)$
the horizontal arrows induce isomorphisms 
$\ker(\alpha)\isomto\ker(\alpha')$ and 
$\coker(\alpha)\isomto\coker(\alpha')$.
Therefore the dotted arrow in $(\dagger\dagger)$ is an isomorphism as
required.
\end{proof}

\begin{remark}
For $n>1$ the square in Proposition \ref{GVB and TFS} is
{\em not} cartesian. This can be seen as follows.
Let $K$ be an algebraically closed field and let
$(\G, \xymatrix@C=2.5ex{\G[p_1]\oplus\G[p_2]\ar@{->>}[r] & Q})$
be a $K$-valued point of 
$\Grass_n(p_1^*\E_{\univ}\oplus p_2^*\E_{\univ})$
such that
$\rk(\G[p_1]\to Q)\neq n\neq \rk(\G[p_2]\to Q)$ (here we need $n>1$).
Let $r:\tilde{C_0}\tensor_kK\to C_0\tensor_kK$ be the
normalization morphism and let
$\F:=\ker(r_*\G\to Q)$ be the torsion free sheaf
on $C_0\tensor_kK$, associated to 
$(\G, \xymatrix@C=2.5ex{\G[p_1]\oplus\G[p_2]\ar@{->>}[r] & Q})$.
By Lemma 2.1 (3) in \cite{Sun} there exists a $K$ valued point
$(\G', \xymatrix@C=2.5ex{\G'[p_1]\oplus\G'[p_2]\ar@{->>}[r] & Q'})$
of
$\Grass_n(p_1^*\E_{\univ}\oplus p_2^*\E_{\univ})$
such that 
$\ker(r_*\G'\to Q)\isomorph\F$ and
$\rk(\G'[p_1]\to Q')=n$. But then Proposition 10.1 in \cite{Kausz}
tells us that the respective fibres of 
$\KGl(p_1^*\E_{\univ},p_2^*\E_{\univ})\to
\Grass_n(p_1^*\E_{\univ}\oplus p_2^*\E_{\univ})$
over the points
$(\G, \xymatrix@C=2.5ex{\G[p_1]\oplus\G[p_2]\ar@{->>}[r] & Q})$
and
$(\G', \xymatrix@C=2.5ex{\G'[p_1]\oplus\G'[p_2]\ar@{->>}[r] & Q'})$
are not isomorphic (as would be the case, if the square in
Proposition  \ref{GVB and TFS} was cartesian).
\end{remark}

\end{document}